\documentclass[11pt,twoside,reqno]{amsart}
\usepackage{import}
\usepackage{newclude}
\usepackage{anysize}
\usepackage{amsmath}
\usepackage{amsthm}
\usepackage{amsmath,amscd}
\usepackage[utf8]{inputenc}
\usepackage{amssymb}
\usepackage{stmaryrd}
\usepackage{wasysym}
\usepackage[all]{xy}
\usepackage{tikz-cd}
\usepackage{mathrsfs}
\usepackage[pagebackref=true]{hyperref} 
\renewcommand*{\backref}[1]{}
\renewcommand*{\backrefalt}[4]{({\tiny%
   \ifcase #1 Not cited.%
         \or Cited on page~#2.%
         \else Cited on pages #2.%
   \fi%
   })}
\usepackage{cleveref}
\usepackage{enumitem}
\usepackage{relsize} 
\usepackage{dsfont}
\usepackage{soul}
\usepackage{filecontents}
\hypersetup{colorlinks=true,linkcolor=black,citecolor=black}
\usepackage{color}
\usepackage{xcolor}
\usepackage[all]{xy}
\usepackage{float}
\usepackage{bm}
\usepackage{mathtools}
\usepackage{thmtools}
\usepackage{setspace}
\usepackage{comment}
\usepackage{tikz}
\usepackage{ytableau}
\usepackage{mathdots}
\usepackage[myheadings]{fullpage}

\numberwithin{equation}{section}
\newcommand\mtop{1in}
\newcommand\mbottom{1in}
\newcommand\mleft{1in}
\newcommand\mright{1in}
\usepackage[top = \mtop, bottom = \mbottom, left = \mleft, right=\mright]{geometry}

\newtheorem{thm}{Theorem}[section]
\newtheorem{example}[thm]{Example}
\newtheorem{prop}[thm]{Proposition}

\newtheorem{lemma}[thm]{Lemma}
\newtheorem{cor}[thm]{Corollary}

\theoremstyle{definition}
\newtheorem{defi}{Definition}[section]

\newtheorem{rmk}{Remark}[section]

\usepackage{scalerel,stackengine}
\stackMath
\newcommand\reallywidehat[1]{%
\savestack{\tmpbox}{\stretchto{%
  \scaleto{%
    \scalerel*[\widthof{\ensuremath{#1}}]{\kern-.6pt\bigwedge\kern-.6pt}%
    {\rule[-\textheight/2]{1ex}{\textheight}}
  }{\textheight}%
}{0.5ex}}%
\stackon[1pt]{#1}{\tmpbox}%
}
\DeclareSymbolFont{bbold}{U}{bbold}{m}{n}
\DeclareSymbolFontAlphabet{\mathbbold}{bbold}

\makeatletter
\def\@tocline#1#2#3#4#5#6#7{\relax
  \ifnum #1>\c@tocdepth 
  \else
    \par \addpenalty\@secpenalty\addvspace{#2}%
    \begingroup \hyphenpenalty\@M
    \@ifempty{#4}{%
      \@tempdima\csname r@tocindent\number#1\endcsname\relax
    }{%
      \@tempdima#4\relax
    }%
    \parindent\z@ \leftskip#3\relax \advance\leftskip\@tempdima\relax
    \rightskip\@pnumwidth plus4em \parfillskip-\@pnumwidth
    #5\leavevmode\hskip-\@tempdima
      \ifcase #1
       \or\or \hskip 1em \or \hskip 2em \else \hskip 3em \fi%
      #6\nobreak\relax
    \hfill\hbox to\@pnumwidth{\@tocpagenum{#7}}\par
    \nobreak
    \endgroup
  \fi}
\makeatother

\makeatletter
\newcommand{\subalign}[1]{%
  \vcenter{%
    \Let@ \restore@math@cr \default@tag
    \baselineskip\fontdimen10 \scriptfont\tw@
    \advance\baselineskip\fontdimen12 \scriptfont\tw@
    \lineskip\thr@@\fontdimen8 \scriptfont\thr@@
    \lineskiplimit\lineskip
    \ialign{\hfil$\m@th\scriptstyle##$&$\m@th\scriptstyle{}##$\hfil\crcr
      #1\crcr
    }%
  }%
}
\makeatother

\newcommand{\R}{\mathbb{R}}
\newcommand{\Z}{\mathbb{Z}}
\newcommand{\Q}{\mathbb{Q}}
\newcommand{\N}{\mathbb{N}}
\newcommand{\C}{\mathbb{C}}
\newcommand{\F}{\mathbb{F}}

\newcommand{\mf}{\mathfrak}

\newcommand{\bbone}{\mathbbold{1}}
\newcommand{\la}{\lambda}

\newcommand{\tth}{^{th}}

\renewcommand{\L}{\Lambda}

\newcommand{\bd}{\mathbf{d}}
\newcommand{\bx}{\mathbf{x}}

\DeclareUnicodeCharacter{FFFD}{}

\DeclareMathOperator{\len}{len}

\DeclareMathOperator{\Span}{span}
\DeclareMathOperator{\Tr}{Tr}

\DeclareMathOperator{\Sig}{Sig}

\DeclareMathOperator{\Proj}{Proj}
\DeclareMathOperator{\GT}{GT}

\DeclareMathOperator{\Alt}{Alt}
\DeclareMathOperator{\Her}{Her}

\DeclareMathOperator{\SN}{SN}

\DeclareMathOperator{\diag}{diag}

\DeclareMathOperator{\Mat}{Mat}

\DeclareMathOperator{\rank}{rank}

\DeclareMathOperator{\alt}{alt}
\DeclareMathOperator{\her}{her}
\DeclareMathOperator{\Nm}{Nm}

\newcommand{\GL}{\mathrm{GL}}

\newcommand{\U}{\mathrm{U}}

\renewcommand{\l}{\lambda}
\newcommand{\bP}{\mathbf{P}}

\title{Non-Archimedean GUE corners and Hecke modules}

\author{Jiahe Shen}
\author{Roger Van Peski}

\date{\today}

\begin{document}

\begin{abstract}
    We compute the joint distribution of singular numbers for all principal corners of a $p$-adic Hermitian (resp. alternating) matrix with additive Haar distribution, the non-archimedean analogue of the GUE (resp. aGUE) corners process. In the alternating case we find that it is a Hall-Littlewood process, explaining---and recovering as a corollary---results of Fulman-Kaplan \cite{fulman2016hall}. In the Hermitian case we obtain a `marginal distribution' of a formal Hall-Littlewood process with both positive and negative transition `probabilities'. The proofs relate natural random matrix operations to structural results of Hironaka \cite{Hironaka_Hermitian,Hironaka_Hermitian_and_Symmetric_I} and Hironaka-Sato \cite{Hironaka} on modules over the spherical Hecke algebra, yielding other probabilistic statements of independent interest along the way.
\end{abstract}

\maketitle

\tableofcontents

\section{Introduction}

\subsection{Preface.} This work concerns non-archimedean analogues of the GUE corners process and its variants in classical random matrix theory, and the algebraic structures behind these objects. We first develop the combinatorics of certain modules over the spherical Hecke algebra, which encodes the action 
\begin{equation}\label{eq:main_operation_intro}
    A \mapsto B^*AB
\end{equation}
of $\GL_N$ on Hermitian and alternating matrices. From this we derive exact formulas for the distribution of these corners processes in terms of Hall-Littlewood polynomials, placing these objects into the Macdonald process framework of Borodin-Corwin \cite{borodin2014macdonald}.

The classical GUE corners process (of rank $k$) is defined as the joint distribution of eigenvalues of all top-left minors of a $k \times k$ Gaussian Hermitian matrix $X$ with density proportional to $\exp(-\Tr X^2/2)$ with respect to the Lebesgue measure\footnote{This means all off-diagonal entries are complex Gaussians and all diagonal entries real Gaussians, independent except for the restriction $a_{ij} = \bar{a}_{j,i}$.}. It is a basic universal object not just in random matrix theory, but also for discrete statistical mechanics models in appropriate limits, such as random tiling models and the six-vertex model, see for instance Johansson-Nordenstam \cite{johansson2006eigenvalues}, Okounkov-Reshetikhin \cite{okounkov2006birth}, Mkrtchyan-Petrov \cite{mkrtchyan2017gue}, Dimitrov \cite{dimitrov2020six}, or Aggarwal-Gorin \cite{aggarwal2022gaussian}. The analogue for alternating matrices is the so-called \emph{anti-symmetric GUE corners process (aGUE)} introduced by Forrester-Nordenstam \cite{forrester2009anti}, and shown to govern certain limits of sorting networks by Gorin-Xu \cite{gorin2024random}. 

Many papers have now been written on similar questions concerning joint distribution of the analogues of singular values, but in the discrete setting of matrices over the integers $\Z$ or $p$-adic integers $\Z_p$. This is motivated on the one hand by the appearance of such distributions in number theory, combinatorics and topology, where they encode universal distributions on abelian groups appearing in these contexts; see for instance Wood \cite{wood2023probability}. In particular, limits of alternating matrices conjecturally model the statistics of Tate-Shafarevich groups of elliptic curves, see Bhargava-Kane-Lenstra-Poonen-Rains \cite{bhargava2013modeling}. Furthermore, universality statements have been shown in the Hermitian one by Lee \cite{lee2022universality}, and in the alternating setting by Nguyen-Wood \cite{nguyen2022local} (which also makes extensive use of matrix corners operations of the type we consider here).

On the other hand, discrete random matrices are interesting probabilistic models in their own right, and analogues of classical questions in real/complex random matrix theory often yield nontrivial discrete analogues, see for instance Assiotis \cite{assiotis2022infinite}, M{\'e}sz{\'a}ros \cite{meszaros2024phase}, and Nguyen and the second author \cite{nguyen2024universality,nguyen2024rank} and \cite{van2020limits,vanpeski2021halllittlewood,van2023p,van2024local,van2023reflecting}. The $p$-adic setting encodes the essential features of the integer one, but simplifies it and provides additional tools. Proofs in this setting are also typically valid for extensions of $\Q_p$; since our results below are essentially algebraic, we will state them in this generality.

Given this history, it is somewhat surprising that the non-archimedean analogues of the GUE and aGUE corners processes have not (to our knowledge) been investigated until now. One reason for this may be that the right structural setting for such results was not at all clear, at least to us. For random matrices in $\GL_n(\C)$ with no symmetry constraints, so-called spherical Hecke algebra---functions on $\GL_n(\C)$ with appropriate invariance properties\footnote{See \Cref{subsec:hecke_ring} for precise definitions in the non-archimedean case.}---has an explicit basis of spherical functions given by degenerations of Macdonald polynomials. This can be parlayed into probabilistic tools through the general theory of Macdonald processes, see for instance Ahn \cite{ahn2019fluctuations}, Borodin-Gorin-Strahov \cite{borodin2018product}, or Gorin-Sun \cite{gorin2018gaussian}. 
This story on the archimedean side provided a very helpful guide for previous results on random non-archimedean matrices \cite{van2020limits}, which turned out to be exact structural analogues.

However, that was for matrices with no symmetry constraints. As soon as one restricts to alternating or Hermitian matrices as we do in this work, it was not clear whether an analogue existed. In the archimedean case, the GUE corners process is known to be a degeneration of a Macdonald process. However, we are not aware of any structural proof, from properties of spherical functions and Hecke algebras, which could have served as a guide for the non-archimedean setting. The arguments we are aware of simply compute explicit formulas for both sides and check that they are equal; this is not difficult because the formulas are simple, but all the same it does not provide any insight into the non-archimedean case, where the corresponding ones are more complicated.

The key to our results was the use of explicit modules over the spherical Hecke algebra as tools for random matrices with symmetry restrictions, just as explicit descriptions of the algebra itself have been useful for those without symmetry. As mentioned, the operation \eqref{eq:main_operation_intro} lets $\GL_n$ act on Hermitian matrices, and this gives an action of the Hecke algebra on appropriate spaces of functions on such matrices. In the $p$-adic groups literature these actions have already been understood explicitly and related to Hall-Littlewood polynomials in a manner similar to the classical Satake isomorphism, by Hironaka-Sato \cite{Hironaka} and Hironaka \cite{Hironaka_Hermitian,Hironaka_Hermitian_and_Symmetric_I} in the alternating and Hermitian cases respectively. 

Combining these structural results and symmetric function combinatorics, we place the non-archimedean versions of GUE and aGUE corners into the framework of Macdonald processes (\Cref{thm:corners_nomarkov_intro}), as well as proving probabilistic results regarding the operation \eqref{eq:main_operation_intro} (\Cref{thm:Product process}) along the way. This provides tools for future asymptotics, analogous to those used for matrices with no symmetry restrictions in the recent works \cite{nguyen2024universality,van2020limits,vanpeski2021halllittlewood,van2024local}, and we mention some directions in \Cref{subsec:applications}. It is also worth mentioning that beyond the Hermitian and alternating cases, many similar Hecke module computations have been done for other symmetry classes in the literature, e.g. \cite{corato2023spherical}. We treated these two cases because they had the most transparent existing results in terms of Hall-Littlewood polynomials, but the method of converting Hecke module computations to random matrix statements should be applicable more generally for other symmetry classes.


\subsection{Random matrix set-up.} Throughout the paper, we work in either one of the following two settings. We refer to \Cref{sec:prelim} for general background on non-archimedean local fields.

\begin{enumerate}[left=0pt]
\item (Alternating case) 

Let $F$ be a non-archimedean local field with characteristic zero\footnote{\label{footnote:pos_char} We believe that the characteristic assumption is merely technical and the results remain true in positive odd characteristic. However, our proofs rely on previous results of \cite{Hironaka_Hermitian}, \cite{Hironaka_Hermitian_and_Symmetric_I}, \cite{Hironaka} established there only in characteristic $0$. See \Cref{sec:char_appendix} for further discussion.};
    
$\mathfrak{o}$ be the ring of integers of $F$;

$\mathfrak{p}$ be the maximal ideal of $\mathfrak{o}$;

$\pi\in\mathfrak{p}$ be a generator of $\mathfrak{p}$.

$k=\mathfrak{o}/\mathfrak{p}$ be the residue field of $F$;

$q=|k|$ be the order of the residue field of $F$;

$|\cdot|:F^\times\rightarrow q^\Z,x\mapsto q^{-v(x)}$ be the absolute value defined over nonzero elements of $F$;

$\Alt_n(F),\Alt_n(\mathfrak{o})$ be the set of alternating\footnote{i.e. $A=-A^T$.} matrices with entries in $F,\mathfrak{o}$, respectively.

\quad

\item(Hermitian case)

Let $F$ be a non-archimedean local field with characteristic zero;

$\mathfrak{o}$ be the ring of integers of $F$;

$\mathfrak{p}$ be the maximal ideal of $\mathfrak{o}$;

$*$ be an involution of $F$, i.e., an isomorphism of order $2$;

$F_0$ be the fixed field of the involution. Also, suppose $F/F_0$ is unramified, so $F_0\cap\mathfrak{o}$ is the ring of integers of $F_0$;

$\pi\in F_0\cap\mathfrak{p}$ be a generator of $\mathfrak{p}$.

$k=\mathfrak{o}/\mathfrak{p}$ be the residue field of $F$;

$q^2=|k|$ be the order of the residue field of $F$;

$|\cdot|:F^\times\rightarrow q^\Z,x\mapsto q^{-2v(x)}$ be the absolute value defined over nonzero elements of $F$;

$B^* = *(B^T)$ where $*$ acts entrywise on the matrix;

$\Her_n(F),\Her_n(\mathfrak{o})$ be the set of Hermitian matrices with entries in $F,\mathfrak{o}$, respectively.

\end{enumerate}

\begin{example}
    In the alternating case, for any prime $p$ one may take $F = \Q_p, \mf{o} = \Z_p$ to be the $p$-adic numbers and integers respectively, in which case $\pi = p$ and $k = \F_p$. In the Hermitian case, one may take $F = \Q_p[\sqrt{d}]$ where $p$ is odd and $d \in \Z \subset \Q_p$ is a non-square in $(\Z/p\Z)^\times$, and $F_0=\Q_p$; then $*$ acts by $*(a+b\sqrt{d}) = a-b\sqrt{d}$. In particular if $p \equiv 3 \pmod{4}$ one may take $F=\Q_p[\sqrt{-1}]$.
\end{example}



We have the following results which are analogous to singular value decomposition for complex matrices, with the compact group $\GL_n(\mf{o})$ playing the role of $U(n)$ or $O(n)$. For any nonsingular non-archimedean Hermitian matrix $A\in \Her_n(F)$, there exists $U\in\GL_{n}(\mathfrak{o})$ such that
$$UAU^*=\diag_{n\times n}(\pi^{\lambda_1},\ldots,\pi^{\lambda_n})$$
for some integers $\infty>\lambda_1\ge\ldots\ge\lambda_n$, see e.g. \cite[page 567]{Hironaka_Hermitian}. We refer the integers $\lambda_i$ as the \textbf{singular numbers} of $A$ and write $\SN^{\her}(A)=(\lambda_1,\ldots,\lambda_n)=\lambda$ in the above case. 

For alternating matrices, singular numbers occur with even multiplicity, and matrices of odd and even size must be treated separately. Similarly to the above, by e.g. \cite[page 483]{Hironaka} one has that for any nonsingular non-archimedean alternating matrix $A\in\Alt_{2n}(F)$, there exists $U\in\GL_{2n}(\mathfrak{o})$ such that 
$$UAU^T=\diag_{2n\times 2n}\left(\begin{pmatrix}0&\pi^{\lambda_1}\\-\pi^{\lambda_1}&0\end{pmatrix},\ldots,\begin{pmatrix}0&\pi^{\lambda_n}\\-\pi^{\lambda_n}&0\end{pmatrix}\right) =: \pi_\lambda^{\alt}$$
for some integers $\la=(\la_1,\ldots,\la_n)$ with $\infty>\lambda_1\ge\ldots\ge\lambda_n$. In the odd case, for any alternating matrix $A\in \Alt_{2n+1}(F)$ with largest possible rank ($2n$), there exists $U\in \GL_{2n+1}(\mathfrak{o})$ such that 
$$UAU^T=\diag_{(2n+1)\times (2n+1)}\left(\begin{pmatrix}0&\pi^{\lambda_1}\\-\pi^{\lambda_1}&0\end{pmatrix},\ldots,\begin{pmatrix}0&\pi^{\lambda_n}\\-\pi^{\lambda_n}&0\end{pmatrix},0\right) $$
for $\la=(\la_1,\ldots,\la_n)$ as before. In both cases we refer the integers $\lambda_i$ as the \textbf{singular numbers} of $A$ and write $\SN^{\alt}(A)=(\lambda_1,\ldots,\lambda_n)=\lambda$; note that $\SN^{\alt}$ discards both the even multiplicity and the trailing $0$ in the odd case. In general, we refer to a weakly decreasing $n$-tuple of integers as an \emph{integer signature} of length $n$, and denote the set of these by
\begin{equation}
    \Sig_n := \{\la = (\la_1,\ldots,\la_n) \in \Z^n: \la_1 \geq \ldots \geq \la_n\}.
\end{equation}
These singular numbers are analogous to the negative logarithms of singular values in the usual setting of real or complex matrices, and they are our main object of study.

On $F$, the natural analogue of the Gaussian measure on $\R$ or $\C$ is the additive Haar measure on the subset $\mf{o}$. For instance, the product of these measures on $F^n$ is invariant under $\GL_n(\mf{o})$, and this property along with independence of coordinates characterizes it uniquely up to scaling, just as Gaussian vectors are the unique $O(n)$-invariant vectors with independent components. The analogues of the GUE and aGUE distributions are just random matrices with entries distributed by this additive Haar measure, which are i.i.d. except for the required symmetry constraints. These are the distributions in \Cref{thm:corners_nomarkov_intro}. 


\subsection{Hall-Littlewood polynomials and corners processes.} We find that the distributions of singular numbers are most neatly expressed in terms of the classical \emph{Hall-Littlewood (Laurent) polynomials}. Full background is given in \Cref{sec:prelim}, but for now let us introduce the facts needed for our results. They are a family of symmetric Laurent polynomials $P_\la(x_1,\ldots,x_n;t)$ in $n$ variables, indexed by integer signatures $\la$ of length $n$, which feature an extra parameter $t$ which we take to be real. Explicitly, they are given by 
\begin{equation}
    P_\la(x_1,\ldots,x_n;t) := \frac{1}{V_\la(t)} \sum_{\sigma \in S_n} \sigma\left(x_1^{\la_1}\cdots x_n^{\la_n} \prod_{1 \leq i < j \leq n} \frac{x_i-tx_j}{x_i-x_j}\right),
\end{equation}
where $\sigma$ acts by permuting the variables and $V_\la$ is the normalizing constant making the polynomial monic. These polynomials are well known in the representation theory of $\GL_n(F)$, see e.g. \cite{macdonald1971spherical}, often under the name of (type $A$) \emph{Macdonald spherical functions}. Also relevant are the skew Hall-Littlewood polynomials $P_{\la/\mu}$ defined by
\begin{equation}\label{eq:skew_def_intro}
    P_\lambda(x_1,\ldots,x_n;t) = \sum_{\mu \in \Sig_k} P_{\lambda/\mu}(x_1,\ldots,x_{n-k};t) P_\mu(x_{n-k+1},\ldots,x_n;t).
\end{equation}
Our first main result expresses the distribution of the non-archimedean versions of the GUE and aGUE corners processes in terms of these, together with the dual polynomials $Q_{\la/\mu}$ which are essentially constant multiples of $P_{\la/\mu}$ (see \eqref{eq:def_Q} for precise definitions).

\begin{thm}\label{thm:corners_nomarkov_intro}
Let $t=1/q$, and let $n\ge 1$ be an integer.
\begin{enumerate}[left=0pt]
\item(Alternating case) Let $A_{2n}$ be a random element of $\Alt_{2n}(\mf{o})$ with above-diagonal entries $a_{i,j}, 1 \leq i < j \leq 2n$ i.i.d. and distributed by the additive Haar measure on $\mf{o}$. For $k=1,\ldots,2n$, let $A_k$ be the top-left $k \times k$ submatrix of $A_{2n}$. Then the joint distribution of their singular numbers is given by
\begin{multline}\label{eq:alt_corner_proc_intro}
    \mathbf{P}(\SN^{\alt}(A_2)=\la^{(2)},\SN^{\alt}(A_3)=\nu^{(3)},\ldots,\SN^{\alt}(A_{2n})=\la^{(2n)}) \\ 
    = \frac{P_{\lambda^{(2)}}(t^{2n-2};t^2)Q_{\lambda^{(2)}/\nu^{(3)}}(t^{3-2n};t^2)P_{\lambda^{(4)}/\nu^{(3)}}(t^{2n-4};t^2)\cdots P_{\lambda^{(2n)}/\nu^{(2n-1)}}(1;t^2)Q_{\lambda^{(2n)}}(t,t^3,\ldots)}{\Pi_{t^2}(1,1,t^2,1,t^2,t^4,\ldots,1,t^2,\ldots,t^{2n-2};t,t^3,\ldots)},
\end{multline}
with the normalizing constant $\Pi_t(\cdots)$ defined in \eqref{eq: Cauchy kernel}.

\item(Hermitian case) Let $A_n\in\Her_n(\mathfrak{o})$ be random with i.i.d entries above the diagonal, distributed according to the additive Haar measure on $\mathfrak{o}$, and i.i.d entries on the diagonal independent of these, distributed according to the additive Haar measure on $\mathfrak{o}\cap F_0$. Then letting $A_1,\ldots,A_n$ be its top-left submatrices as before, the joint distribution of their singular numbers is given by
\begin{multline}\label{eq:herm_corner_proc_intro}
\mathbf{P}(\SN^{\her}(A_1)=\lambda^{(1)},\ldots,\SN^{\her}(A_n)=\lambda^{(n)})=\frac{1}{\Pi}\sum_{\nu^{(1)}\in\Sig_1^+,\ldots,\nu^{(n-1)}\in\Sig_{n-1}^+} P_{\lambda^{(1)}}((-t)^{n-1};-t)
\\\times \left(\prod_{i=1}^{n-1}Q_{\lambda^{(i)}/\nu^{(i)}}(-(-t)^{i+1-n};-t)P_{\lambda^{(i+1)}/\nu^{(i)}}((-t)^{n-i-1};-t)\right)Q_{\lambda^{(n)}}(t,-t^2,\ldots;-t)
\end{multline}
with the normalizing constant $\Pi$ given by 
$$\Pi=\Pi_{-t}(1,1,-t,1,-t,t^2,\ldots,1,-t,\ldots,(-t)^{n-1};t,-t^2,\ldots).$$
\end{enumerate}
\end{thm}

Probability measures on sequences of signatures of the same form as in \Cref{thm:corners_nomarkov_intro} are called \emph{Hall-Littlewood processes}. They are special cases of the Macdonald processes introduced by Borodin-Corwin \cite{borodin2014macdonald}, which have been applied in many contexts since then (some mentioned already above).

\Cref{thm:corners_nomarkov_intro} also yields attractive expressions for the marginal distributions of singular numbers of a single alternating or Hermitian matrix $A_k$, see \Cref{thm:The Hall-Littlewood measure of i.i.d entries}. In the even alternating case these read
\begin{align}\label{eq:haar_alt_intro}
\begin{split}
 \mathbf{P}(\SN^{\alt}(A_{2n}) = \la) &= \frac{P_\lambda(1,t^2,\ldots,t^{2n-2};t^2)Q_\lambda(t,t^3,\ldots;t^2)}{\Pi_{t^2}(1,t^2,\ldots,t^{2n-2};t,t^3,\ldots)} \\ 
 &= t^{\sum_{i \geq 1} (4i-3)\la_i}\frac{(t;t)^{2n}}{\prod_{i\ge 0}(t^2;t^2)_{m_i(\lambda)}},       
\end{split}
\end{align}
where $(a;t)_n := (1-a)(1-at) \cdots (1-at^{n-1})$ is the usual $q$-Pochhammer symbol and $m_i(\la) = \#\{j: \la_j = i\}$. In the Hermitian case, one obtains a similar formula
\begin{align}\label{eq:haar_her_intro}
    \begin{split}
\mathbf{P}(\SN^{\her}(A_n) = \la) &=\frac{P_\lambda(1,-t,t^2,\ldots,(-t)^{n-1};-t)Q_\lambda(t,-t^2,\ldots;-t)}{\Pi_{-t}(1,-t,t^2,\ldots,(-t)^{n-1};t,-t^2,t^3,\ldots)} \\ 
&= t^{\sum_{i \geq 1} (2i-1) \la_i}\frac{(t^2;t^2)_n}{\prod_{i\ge 0}(-t;-t)_{m_i(\lambda)}}.        
    \end{split}
\end{align}
Measures of this form are known as Hall-Littlewood measures.

The formula \eqref{eq:haar_alt_intro} recovers a result of Fulman-Kaplan \cite[Theorem 3.2]{fulman2018random}. They proved this by matching explicit formulas for Hall-Littlewood polynomials with an explicit expression for the probability computed earlier by Bhargava-Kane-Lenstra-Poonen-Rains \cite{bhargava2013modeling}. However, the reasons behind this coincidence of formulas remained mysterious, to us at least. Our methods give not just an alternate proof, but a conceptual reason why Hall-Littlewood polynomials should appear in such formulas: they are spherical functions on the relevant groups and homogeneous spaces, as we see below.

In the Hermitian case, a formula for the left hand side of \eqref{eq:haar_her_intro} in terms of module automorphisms---without reference to Hall-Littlewood polynomials---was given by Lee \cite{lee2022universality}. It may be easily checked to coincide with the right hand side of \eqref{eq:haar_her_intro}, so our proof gives an alternate derivation.

\subsection{Matrix operations and Hall-Littlewood polynomials.} \Cref{thm:corners_nomarkov_intro} is closely related to the following result, \Cref{thm:Corner process}, which considers a single step of the corners dynamics.

\begin{defi}
    \label{def:invariant_dist}
    We say that a random matrix $A \in \Alt_n(F)$ (resp. $A \in \Her_n(F)$) is \emph{$\GL_n(\mf{o})$-invariant} if $BAB^T = A$ (resp. $BAB^*=A$) in distribution for any fixed $B \in \GL_n(\mf{o})$. 
\end{defi}

Explicitly, the unique invariant measure on $\Alt_{2n}(F)$ with singular numbers $\la \in \Sig_n$ is given by $U\pi_\lambda^{\alt}U^T$ where $U$ is distributed by the Haar probability measure on $\GL_N(\mf{o})$, and the analogous result is true for odd alternating and for Hermitian matrices. 





\begin{thm}\label{thm:Corner process}
Let $t=1/q$, $n\ge 1$ be an integer.
\begin{enumerate}[left=0pt]
\item(Alternating case)\begin{enumerate}[left=0pt]
\item Let $\nu \in \Sig_n$ and let $A\in\Alt_{2n+1}(F)$ be random, with the unique $GL_{2n+1}(\mathfrak{o})$-invariant distribution such that $\SN^{\alt}(A)=\la$. Let $A_{2n}$ be the $2n\times2n$ corner of $A$ on the top left. Then $\SN^{\alt}(A_{2n})\in\Sig_n$ has distribution $\mathbf{P}_{2n<2n+1}^{\alt}(\cdot\mid\nu)$ defined by

\begin{equation}\label{eq:odd_alt_corner}
\mathbf{P}_{2n<2n+1}^{\alt}(\la \mid \nu)=\frac{Q_{\la/\nu}(t;t^2)P_\la(1,t^2,\ldots,t^{2n-2};t^2)}{P_\nu(1,t^2,\ldots,t^{2n-2};t^2)\prod_{t^2}(t;1,t^2,\ldots,t^{2n-2})}
\end{equation}

\item Let $\la \in \Sig_n$, and let $A\in \Alt_{2n}(F)$ be random, with the unique $GL_{2n}(\mathfrak{o})$-invariant distribution such that $\SN^{\alt}(A)=\lambda$. Let $A_{2n-1}$ be the top left $(2n-1)\times(2n-1)$ corner of $A$. Then $\SN^{\alt}(A_{2n-1})\in\Sig_{n-1}$ has distribution $\mathbf{P}_{2n-1<2n}^{\alt}(\cdot\mid\lambda)$ given by

\begin{equation}\label{eq:even_alt_corner}
\mathbf{P}_{2n-1<2n}^{\alt}(\nu\mid\lambda)=\frac{P_{\lambda/\nu}(1;t^2)P_\nu(t^2,t^4,\ldots,t^{2n-2};t^2)}{P_\lambda(1,t^2,\ldots,t^{2n-2};t^2)}
\end{equation}
\end{enumerate}\label{item:corner_alternating}

\item(Hermitian case) Let $\la \in \Sig_n$ and let $A\in\Her_n(F)$ be random, with the unique $GL_n(\mathfrak{o})$-invariant distribution such that $\SN^{\her}(A)=\lambda$. Let $A_{n-1}$ be the top left $(n-1)\times(n-1)$ corner of $A$. Then $\SN^{\her}(A_{n-1})\in\Sig_{n-1}$ has distribution $\mathbf{P}_{n-1<n}^{\her}(\cdot\mid\lambda)$ given by
\begin{equation}\label{eq:herm_corner}
\mathbf{P}_{n-1<n}^{\her}(\nu\mid\lambda)=\sum_{\kappa\in\Sig_{n-1}} \frac{P_{\la/\kappa}(1;-t) Q_{\nu/\kappa}(-1;-t)P_\nu(-t,\ldots,(-t)^{n-1};-t)}{P_\la(1,\ldots,(-t)^{n-1};-t)\Pi_{-t}(1;t,\ldots,-(-t)^{n-1})}  
\end{equation}\label{item:hermitian_corner}
\end{enumerate}
\end{thm}

\begin{rmk}\label{rmk:formal_macdonald_process}
It is interesting that formally speaking, the transition probabilities \eqref{eq:herm_corner} are given by a composition of two Markov maps of the same form as in \eqref{eq:odd_alt_corner} and \eqref{eq:even_alt_corner}, since the summands in \eqref{eq:herm_corner} are given by
    \begin{equation}\label{eq:branch_prod_herm_rmk}
        \frac{P_{\la/\kappa}(1;-t) P_\kappa(-t,\ldots,(-t)^{n-1};-t)}{P_\la(1,\ldots, (-t)^{n-1};-t)} \\ \times \frac{Q_{\nu/\kappa}(-1;-t)P_\nu(-t,\ldots,(-t)^{n-1};-t)}{P_\kappa(-t,\ldots,(-t)^{n-1};-t)\Pi_{-t}(1;t,\ldots,-(-t)^{n-1})}.
    \end{equation}
However, there is a crucial difference: depending on $\kappa$, the terms \eqref{eq:branch_prod_herm_rmk} may be either positive or negative (as may be checked by \Cref{prop:hl_principal_formulas}), though the sum over $\kappa$ is always positive. Possibly-negative weights defined on sequences of partitions in this manner are often called \emph{formal Macdonald process} after Borodin-Corwin-Gorin-Shakirov \cite{borodin2016observables}, and these `stochastic processes' have been used in intermediate steps of proofs such as \cite[Theorem 1.1]{borodin2016observables}. However, to the best of our knowledge, this work is the first time that `marginal distributions' of multiple partitions in a genuine formal Macdonald process have arisen as the answer to a probabilistic question. 

We note also that the appearance of the parameter $-1/q$ is similar to what has been observed in the context of representation theory of unitary groups over finite fields under the name of Ennola duality \cite{ennola1963characters}, see also e.g. Kawanaka \cite{kawanaka1985generalized} or Cuenca-Olshanski \cite{cuenca2022infinite}. We expect them to be related, see \Cref{rmk:ennola}, and it would be interesting to understand the details. A seemingly related probabilistic result is the appearance of Hall-Littlewood measures at $t=-1/q$ in asymptotics of unipotent Jordan blocks of random unitary matrices over a finite field, see Fulman \cite[Section 4.4]{fulman_thesis}.
\end{rmk}

\subsection{Matrix products and Hecke modules.} In previous work \cite{van2020limits}, the second author showed results analogous to \Cref{thm:Corner process} for matrices without symmetry restrictions. There, the key was to understand the singular numbers of matrix products, as then corners could be understood through the limiting case where one matrix was a projection. The singular numbers of matrix products were then parametrized via the so-called spherical Hecke algebra $H(\GL_n(F),\GL_n(\mf{o}))$ of $\GL_n(\mf{o})$-invariant functions on $\GL_n(F)$, which is classically known to be governed by Hall-Littlewood polynomials \cite[Chapter V]{Macdonald}. 

At first, it was very unclear to us whether there was any analogue of this strategy for matrices with extra symmetry constraints. The product of two alternating matrices is not an alternating matrix, for instance. The main novelty of our methods here is to find an analogue of this strategy, which requires an understanding of the structure not just of the Hecke algebra, but also certain modules over it.

One has a natural action of $B \in \GL_{2n}(F)$ on $\Alt_{2n}(F)$ by 
$$A \mapsto B^TAB,$$
and similarly with $\Her_n(F)$, rendering suitable spaces of functions on $\Alt_{2n}(F)$ and $\Her_n(F)$ as modules over the spherical Hecke algebra; see \Cref{subsec:hecke_ring} for definitions of the Hecke algebra, and \Cref{sec:alt_hecke} and \Cref{sec:her_hecke} for modules over it. Fortuitously, a precise explicit description of this action was given by Hironaka-Sato \cite{Hironaka} (in the alternating case) and Hironaka \cite{Hironaka_Hermitian,Hironaka_Hermitian_and_Symmetric_I} (in the Hermitian case). Perhaps surprisingly it is also parametrized by the same (type $A$) Hall-Littlewood polynomials. 

We will give the precise module structure statements extracted from those works in \Cref{thm: alt_Hironaka} and \Cref{thm: her_Hironaka} after introducing the necessary notation. However, they are equivalent---with some computation---to a probabilistic description of the above $\GL_n(F)$-actions with random inputs, which is easy to state now. The probabilities are described in terms certain modifications of the classical Littlewood-Richardson coefficients, $c_{\mu,\nu}^{\alt,\lambda}(t)$ and $c_{\mu,\nu}^{\her,\lambda}(t)$, defined by
\begin{equation}\label{eq: alt_Littlewood Richardson}
P_\mu(x_1,x_1t,\ldots,x_n,x_nt;t)P_\nu(x_1,\ldots,x_n;t^2)=\sum_\lambda c_{\mu,\nu}^{\alt,\lambda}(t)P_\lambda(x_1,\ldots,x_n;t^2)
\end{equation}
and
\begin{equation}\label{eq: her_Littlewood Richardson}
P_\mu(x_1^2,\ldots,x_n^2;t^2)P_{\nu}(x_1,\ldots,x_n;-t)=\sum_\lambda c_{\mu,\nu}^{\her,\lambda}(t)P_{\lambda}(x_1,\ldots,x_n;-t)\end{equation}
respectively. 


\begin{thm}\label{thm:Product process} Let $t=1/q$ and $n\ge 1$ be an integer.

\begin{enumerate}[left=0pt]
\item(Alternating case) Fix $\nu \in \Sig_n, \mu \in \Sig_{2n}$. Let $A\in\Alt_{2n}(F)$ be a random matrix with the unique $GL_n(\mathfrak{o})$-invariant distribution such that $\SN^{\alt}(A)=\nu $, and let $B\in \GL_{2n}(F)$ be any (fixed or random) matrix with singular numbers $\SN(B)=\mu $. Then $\SN^{\alt}(B^TAB)$ has distribution 

$$
\mathbf{P}(\SN^{\alt}(B^TAB)=\lambda)=c_{\mu,\nu}^{\alt,\lambda}(t)\frac{P_\lambda(1,t^2,\ldots,t^{2n-2};t^2)}{P_\mu(1,t,\ldots,t^{2n-2},t^{2n-1};t)P_\nu(1,t^2,\ldots,t^{2n-2};t^2)}.$$

\item(Hermitian case) Fix $\nu,\mu \in \Sig_n$. Let $A\in\Her_n(F)$ be a random matrix with the unique $\GL_n(\mathfrak{o})$-invariant distribution such that $\SN^{\her}(A) = \nu$, and let $B\in \GL_{2n}(F)$ be any (fixed or random) matrix with singular numbers $\SN(B)=\mu $. Then $\SN^{\her}(B^*AB)$ has distribution 

$$
\mathbf{P}(\SN^{\her}(B^*AB)=\lambda)=c_{\mu,\nu}^{\her,\lambda}(t)\frac{P_\lambda(1,-t,\ldots,(-t)^{n-1};-t)}{P_\mu(1,t^2,\ldots,t^{2n-2};t^2)P_\nu(1,-t,\ldots,(-t)^{n-1};-t)}.$$
\end{enumerate}
\end{thm}



\begin{rmk}
In the case of matrices without symmetry restrictions, the distribution of singular numbers of a product of two such matrices is encoded by the usual structure coefficients\footnote{When $t=0$ these are the Littlewood-Richardson coefficients, and we will sometimes use this name to refer to them at general $t$ as well.} $c_{\mu,\nu}^\la(t)$ of the Hall-Littlewood polynomials, rather than the nonstandard versions in \eqref{eq: alt_Littlewood Richardson} and \eqref{eq: her_Littlewood Richardson}, see \cite[Theorem 1.3, Part 3]{van2020limits}. However, we have not seen the modified Littlewood-Richardson-type coefficients of \eqref{eq: alt_Littlewood Richardson} and \eqref{eq: her_Littlewood Richardson} studied explicitly in the literature, and it was quite surprising to us that they arose naturally in these random matrix problems. 
\end{rmk}


\begin{rmk}
    The classical Hall-Littlewood structure coefficients are rich combinatorial objects, see e.g. Kirillov \cite{kirillov1998new} or Schwer \cite{schwer2006galleries}, and it would be interesting to understand the combinatorics and positivity properties of $c_{\mu,\nu}^{\alt,\la}$ and $c_{\mu,\nu}^{\her,\la}$.
\end{rmk}

\subsection{Applications.} \label{subsec:applications} In addition to the structural interest of our results, the connection to Macdonald processes provides tools for which a few immediate applications present themselves:
\begin{enumerate}
    \item Previous results of Olshanski-Vershik \cite{olshanski1996ergodic} classified all ergodic, unitarily-invariant measures on the space of infinite Hermitian matrices (over $\C$). Bufetov-Qiu \cite{bufetov2017ergodic} solved the analogous problem for invariant measures on the space of infinite matrices with no symmetry constraints over a non-archimedean local field $F$. Later, \cite{vanpeski2021halllittlewood} gave an alternate proof of their result using Hall-Littlewood formulas from \cite{van2020limits}, which are analogous to the ones we derive in \Cref{thm:Corner process}. The same arguments as in \cite{vanpeski2021halllittlewood}, combined with \Cref{thm:Corner process}, should yield generalizations of \cite{bufetov2017ergodic} to invariant measures on infinite Hermitian matrices (analogous to \cite{olshanski1996ergodic}) and infinite alternating matrices.
    \item Recent work \cite{van2024local} described the local limits of singular numbers of products $A_k A_{k-1} \cdots A_1$ of random matrices with no symmetry constraints, in terms of a new family of discrete distributions. Our \Cref{thm: product convolution and HL process}, combined with the techniques in those works, should find the analogous distributions for `Hermitized products' $A_k^* \cdots A_1^* B A_1 \cdots A_k$, and similarly in the alternating case. We believe that the distributions in the alternating case will be the same as in \cite{van2024local}, up to changing the Hall-Littlewood parameter. However, the negative $t$ parameter in the Hermitian case should lead to different analysis and a different final answer (the formula for the limit distribution in \cite[Proposition 5.1]{van2024local} does not make sense for $t \in (-1,0)$). It would be interesting to investigate this.  
    \item A related classical appearance of the usual Hall-Littlewood structure constants is in the Hall algebra of the category of finite $\mf{o}$-modules \cite[Chapter II]{Macdonald}. In upcoming work, we use the structural results here to study modules over this algebra corresponding to $\mf{o}$-modules with alternating (resp. Hermitian) forms, studied for instance in \cite{bhargava2013modeling,delaunay} and \cite{lee2022universality} respectively. 
\end{enumerate}

\subsection{Plan of paper.} In \Cref{sec:prelim} we give preliminaries on non-archimedean random matrices and symmetric polynomials. We deduce \Cref{thm:Product process} from the structure of Hecke modules in \Cref{sec:alt_hecke} and \Cref{sec:her_hecke} in the alternating and Hermitian cases respectively; these two sections are almost identical, but there are enough minor differences that we feel it better to write everything out fully in both, even though the proofs and the wording are mostly quite similar. In \Cref{sec:HL} we use \Cref{thm:Product process} and Hall-Littlewood combinatorics to deduce \Cref{thm:Corner process} and \Cref{thm:corners_nomarkov_intro}. \Cref{sec:char_appendix} is a brief discussion of positive characteristic issues in the literature.

\subsection{Acknowledgments.} We thank Amol Aggarwal, Alexei Borodin, Bill Casselman, Ivan Corwin, Cesar Cuenca, Vadim Gorin, Yumiko Hironaka, Chao Li, Omer Offen, Grigori Olshanski, and Yiannis Sakellaridis for helpful conversations and comments. RVP was partially supported by the European Research Council (ERC), Grant Agreement No. 101002013, and JS was partially supported by NSF grant DMS-2246576 and Simons Investigator grant 929852.

\section{Preliminaries}\label{sec:prelim}

\subsection{Non-archimedean matrix background}

We begin with a few paragraphs of background which are essentially quoted from \cite{van2020limits}, and are a condensed version of the exposition in Evans \cite[Section 2]{evans2002elementary}. Many of these are modified versions of results in \cite[Section 2]{van2023p} which were proven for matrices without symmetry restrictions, which we adapt here to the setting of alternating and Hermitian non-archimedean matrices.

Fix a non-archimedean local field $F$ with $\mathfrak{o}$ its ring of integers, $\pi$ a generator of the maximal ideal $\mathfrak{p}$, and $q=|\mathfrak{o}/\mathfrak{p}|$ the order of the residue field. Any nonzero element $x\in F^\times$ may be written as $x=\pi^my$ with $m\in\Z$ and $y\in\mathfrak{o}^\times$. Define $|\cdot|: F \to \R_{\ge 0}$ by setting $|x| = q^{-m}$ for $x$ as before, and $|0|=0$. Then $|\cdot|$ defines a norm on $F$ and $d(y_1,y_2) :=|y_1-y_2|$ defines a metric. We additionally define $v(x)=m$ for $x$ as above and $v(0)=\infty$, so $|x|=q^{-v(x)}$.

$F$ is noncompact but is equipped with a left- and right-invariant (additive) Haar measure; this measure is unique if we normalize so that the compact subgroup $\mathfrak{o}$ has measure $1$. The restriction of this measure to $\mathfrak{o}$ is the unique Haar probability measure on $\mathfrak{o}$, and is explicitly characterized by the fact that its pushforward under any map $r_n:\mathfrak{o}\to\mathfrak{o}/\mathfrak{p}^n$ is the uniform probability measure. For concreteness, it is often useful to view elements of $\mathfrak{o}$ as `power series in $\pi$' $a_0 + a_1 \pi + a_2 \pi^2 + \cdots$, with $a_i\in\{b_0,\ldots,b_{q-1}\}$, a set of representatives of the residue field $\mathfrak{o}/\mathfrak{p}$. Clearly these specify a coherent sequence of elements of $\mathfrak{o}/\mathfrak{p}^n$ for each $n$. The Haar probability measure then has the alternate explicit description that each $a_i$ is iid uniformly random from $\{b_0,\ldots,b_{q-1}\}$. Additionally, $F$ is isomorphic to the ring of Laurent series in $\pi$, defined in exactly the same way.

Similarly, $\GL_N(F)$ has a unique left- and right-invariant measure for which the total mass of the maximal compact subgroup $\GL_N(\mathfrak{o})$ is $1$. The restriction of this measure to $\GL_N(\mathfrak{o})$, which we denote by $M_{Haar}(\GL_N(\mathfrak{o}))$, pushes forward to $\GL_N(\mathfrak{o}/\mathfrak{p}^n)$ and is the uniform measure on these finite groups. 

The following result provides the standard form of the matrices we study throughout the paper. 

\begin{prop}\label{prop:smith}
\begin{enumerate}[left=0pt]
\item(Alternating case)\begin{enumerate}[left=0pt]
\item For any nonsingular alternating matrix $A\in\Alt_{2n}(F)$, there exists $U\in\GL_{2n}(\mathfrak{o})$ such that

$$UAU^T=\diag_{2n\times 2n}(\begin{pmatrix}0&\pi^{\lambda_1}\\-\pi^{\lambda_1}&0\end{pmatrix},\ldots,\begin{pmatrix}0&\pi^{\lambda_n}\\-\pi^{\lambda_n}&0\end{pmatrix})$$
for some integers $\infty>\lambda_1\ge\ldots\ge\lambda_n$. 

\item For any alternating matrix $A\in \Alt_{2n+1}(F)$, there exists $U\in \GL_{2n+1}(\mathfrak{o})$ such that

$$UAU^T=\diag_{(2n+1)\times(2n+1)}(\begin{pmatrix}0&\pi^{\lambda_1}\\-\pi^{\lambda_1}&0\end{pmatrix},\ldots,\begin{pmatrix}0&\pi^{\lambda_n}\\-\pi^{\lambda_n}&0\end{pmatrix},0)$$
for some integers $\infty\ge\lambda_1\ge\ldots\ge\lambda_n$. Also, $\lambda_1<\infty$ if and only if $A$ has corank $1$.
\end{enumerate}

\item(Hermitian case)For any nonsingular non-archimedean Hermitian matrix $A\in \Her_n(F)$, there exists $U\in\GL_{n}(\mathfrak{o})$ such that

$$UAU^*=\diag_{n\times n}(\pi^{\lambda_1},\ldots,\pi^{\lambda_n})$$
for some integers $\infty>\lambda_1\ge\ldots\ge\lambda_n$.
\end{enumerate}
\end{prop}

For the alternating case see page 483 of \cite{Hironaka}, and for the Hermitian case see page 567 of \cite{Hironaka_Hermitian}. We will sometimes omit the dimensions $2n\times 2n,(2n+1)\times(2n+1),n\times n$ in the $\diag$ notation when they are clear from context.

Similarly to eigenvalues and singular values, singular numbers have a variational characterization. We first recall the version for singular values, one version of which states that for $A\in\Alt_{2n}(\C)$ with singular values $a_1 \geq \ldots \geq a_n\ge 0$,
\begin{equation}\label{eq:alternating_minmax}
\prod_{i=1}^k a_i^2 = \sup_{\substack{V \subset \C^{2n}: \dim(V) = 2k}} |\det(\Proj_V \circ A|_V)|
\end{equation}
where $\Proj$ is the orthogonal projection and $A|_V$ is the restriction of the linear operator $A$ to the subspace $V$. Likewise, for $A\in\Her_n(\C)$ with singular values $|a_1| \geq \ldots \geq |a_n|$,
\begin{equation}\label{eq:Hermitian_minmax}
\prod_{i=1}^k |a_i| = \sup_{\substack{V \subset \C^n: \dim(V) = k}} |\det(\Proj_V \circ A|_V)|
\end{equation}  \eqref{eq:alternating_minmax} and \eqref{eq:Hermitian_minmax} hold because the right hand side is unchanged by multiplying $A$ by unitary matrices, hence $A$ may be taken to be diagonal with singular values on the diagonal by singular value decomposition, at which point the result is easy to see. 

For non-archimedean matrices, we state the result differently without referring to orthogonal projection, since  $\GL_n(\mathfrak{o})$ does not preserve a reasonable inner product as $\U(n)$. It turns out that one does not have to work with arbitrary projections and subspaces, but may instead consider only minors of the matrix $A$. Here by $k \times k$ minor, we mean any $k \times k$ matrix obtained by deleting rows and columns of the original matrix. Also, by $k\times k$ principal minor, we mean a $k\times k$ minor that is obtained by taking the same rows and columns.

\begin{prop}\label{thm:submatrices_suffice}
\begin{enumerate}[left=0pt]
\item(Alternating case) Let $A\in\Alt_{2n}(F)$ or $A\in\Alt_{2n+1}(F)$ with $\SN^{\alt}(A)=(\lambda_1,\ldots,\lambda_n)$. Then for any $1\le k\le n$,
\begin{equation}\label{eq:alternating_minor}
2(\la_n+\cdots+\la_{n-k+1}) = \inf_{A' \text{ $2k \times 2k$ principal minor of $A$}} v(\det(A')).
\end{equation}

\item(Hermitian case) Let $A\in\Her_n(F)$ with $\SN^{\her}(A)=(\lambda_1,\ldots,\lambda_n)$. Then for any $1\le k\le n$,
\begin{equation}\label{eq:Hermitian_minor}
\la_n+\cdots+\la_{n-k+1} = \inf_{A' \text{ $k \times k$ minor of $A$}} v(\det(A')).
\end{equation}
\end{enumerate}
\end{prop}
\begin{proof}

The Hermitian case is a simple corollary of Proposition 2.3 of \cite{van2023p}, and the alternating case closely follows the proof of that result. We will only prove the alternating even order case here, and the odd case is similar. Clearly, the statement holds when 

$$A=\diag(\begin{pmatrix}0&\pi^{\lambda_1}\\-\pi^{\lambda_1}&0\end{pmatrix},\ldots,\begin{pmatrix}0&\pi^{\lambda_n}\\-\pi^{\lambda_n}&0\end{pmatrix})$$

Since $UAU^T = \diag(\begin{pmatrix}0&\pi^{\lambda_1}\\-\pi^{\lambda_1}&0\end{pmatrix},\ldots,\begin{pmatrix}0&\pi^{\lambda_n}\\-\pi^{\lambda_n}&0\end{pmatrix})$ for some $U\in \GL_{2n}(\mathfrak{o})$, it suffices to show that the right hand side of \eqref{eq:alternating_minor} is invariant under the action of $\GL_{2n}(\mathfrak{o})$ on both sides, i.e.,
\begin{equation}
\label{eq:minor_invariant}
\inf_{A' \text{ $2k \times 2k$ principal minor of $B$}} v(\det(A')) = \inf_{A' \text{ $2k \times 2k$ minor of $UBU^T$}} v(\det(A'))
\end{equation}
for any $B = (b_{i,j})_{\substack{1 \leq i \leq 2n \\ 1 \leq j \leq 2n}} \in \Alt_{2n}(F)$ and $U\in \GL_{2n}(\mathfrak{o})$. 

First note that since $\GL_{2n}(\mathfrak{o})$ is generated by the three elementary row operations:
\begin{enumerate}[left=0pt]
\item[(i)] elementary transposition matrices, 
\item[(ii)] unit multiple matrices $\diag(1[i-1],u,1[2n-i])$ for $u \in \mathfrak{o}^\times, 1 \leq i \leq 2n$, and
\item[(iii)] matrices $(\bbone(i=j) + \bbone(i=x,j=y))_{1 \leq i,j \leq 2n}$ for some $x \neq y$,
\end{enumerate}
it suffices to prove \eqref{eq:minor_invariant} when $U$ is each one of the above types. This is clear for types (i) and (ii). Suppose $U$ is of type (iii). Now we need the following lemma:

\begin{lemma}\label{lemma:judge the singular numbers}
Let $D\in\Alt_{2k-1}(F)$, and $\alpha_1=(c_{1,1},\ldots,c_{(2k-1),1})^t,\alpha_2=(c_{1,2},\ldots,c_{(2k-1),2})^t\in M_{(2k-1)\times 1}(F)$. Also, let 

$$D_1=\begin{pmatrix}0 & \alpha_1^t\\ -\alpha_1 & D\end{pmatrix},D_2=\begin{pmatrix}0 & \alpha_2^t\\ -\alpha_2 & D\end{pmatrix},D_{12}=\begin{pmatrix}0 & \alpha_1^t+\alpha_2^t\\ -\alpha_1-\alpha_2 & D\end{pmatrix}$$

Then $\det (D_{12})=\det D_1+\det D_2\pm2\sqrt{\det D_1\det D_2}$. Here $\pm\sqrt{\det D_1\det D_2}$ is one of two elements in $F$ such that its square equals $\det D_1\det D_2$.
\end{lemma}

\begin{proof}
After we take action of $GL_{2k-1}(\mathfrak{o})$ on both sides, there is no loss for us to assume $D$ has the form

$$D=\diag(0,\begin{pmatrix}0&\pi^{\nu_1}\\-\pi^{\nu_1}&0\end{pmatrix},\ldots,\begin{pmatrix}0&\pi^{\nu_{k-1}}\\-\pi^{\nu_{k-1}}&0\end{pmatrix})$$

Then $\det(D_1)=\alpha_{11}^2\det D,\det(D_2)=\alpha_{12}^2\det D$, and $\det(D_{12})=(\alpha_{11}+\alpha_{12})^2\det D$, hence the lemma holds.
\end{proof}

Turning back to the proof, we have
\begin{equation}
UBU^T = (b_{i,j} + \bbone(i=x)b_{y,j}+\bbone(j=x)b_{i,y})_{\substack{1 \leq i \leq 2n \\ 1 \leq j \leq 2n}}
\end{equation}
differs from $B$ only in the $x\tth$ row and column. For any set of indices $I_x=\{x,i_1,\ldots,i_{k-1}\}$ which include the row and column $x$, let $B_{I_x}$ be the corresponding minor. Then \Cref{lemma:judge the singular numbers} shows that
\begin{equation}
\det (UBU^T)_{I_x} = \det B_{I_x} + \det B_{I_y}\pm 2\sqrt{\det B_{I_x}\det B_{I_y}}
\end{equation}
so by the ultrametric inequality, one can verify that
\begin{equation}\label{eq:ultrametric}
\min(v(\det (UBU^T)_{I_x}),v(\det B_{I_y}))=\min(v(\det B_{I_x}),v(\det B_{I_y})).
\end{equation}
Indeed, if $v(\det B_{I_x})>v(\det(B_{I_y}))$, then $v(\det B_{I_x})>v(\sqrt{\det B_{I_x}\det B_{I_y}})>v(\det B_{I_y})$, and LHS of \eqref{eq:ultrametric} $=v(\det B_{I_y})=$ RHS of\eqref{eq:ultrametric}; If $v(\det B_{I_x})\le v(\det(B_{I_y}))$, then $v(\det B_{I_x})\le v(\sqrt{\det B_{I_x}\det B_{I_y}})\le v(\det B_{I_y})$, and LHS of \eqref{eq:ultrametric} $=v(\det B_{I_x})=$ RHS of\eqref{eq:ultrametric}. This ends the proof.
\end{proof}

\begin{rmk}
As we have seen in \eqref{eq:alternating_minor}, the result of the LHS only depends on the principal minors. Nevertheless, this statement does not hold for the Hermitian case. As a counterexample, the minimal value of the entries in the matrix $\begin{pmatrix} \pi & 1 \\ 1 & \pi\end{pmatrix}$ does not appear on the diagonal.
\end{rmk}

The following theorem will be useful in the coming sections.

\begin{thm}\label{thm:invertible probability}
\begin{enumerate}[left=0pt]
\item(Alternating case)Let $A\in\Alt_{2n}(\mathfrak{o})$ and $\mathbf{P}_{2n}^{\alt}$ be the same as \Cref{thm:The Hall-Littlewood measure of i.i.d entries}. Then we have

$$\mathbf{P}_{2n}^{\alt}(0)=\frac{|\Alt_{2n}(k)\cap\GL_{2n}(k)|}{|\Alt_{2n}(k)|}=(1-q^{1-2n})(1-q^{3-2n})\cdots(1-q^{-1})$$

\item(Hermitian case)Let $A\in\Her_n(\mathfrak{o})$ and $\mathbf{P}_n^{\her}$ be the same as \Cref{thm:The Hall-Littlewood measure of i.i.d entries}. Then we have

$$\mathbf{P}_n^{\her}(0)=\frac{|\Her_n(k)\cap\GL_n(k)|}{|\Her_n(k)|}=(1+(-q)^{-n})(1+(-q)^{1-n})\cdots(1-q^{-1})$$
\end{enumerate}
\end{thm}

\begin{proof}
\begin{enumerate}[left=0pt]
\item See Lemma 3.6 of \cite{bhargava2013modeling}.
\item We prove the equation by induction. It is clear that the equation holds when $n=1$. Suppose we already know that the equation holds for $\le n-1$. In this case, to make the matrix $A\in\Her_n(\mathfrak{o})$ invertible, one of the following cases must occur:
\begin{enumerate}[left=0pt]
\item The left upper corner of $A$ belongs to $\mathfrak{o}^\times\cap F_0$, which has probability $1-q^{-1}$. Then, we can use this element to eliminate the first row and column and deal with the right lower $(n-1)\times (n-1)$ corner. The measure of this corner is still i.i.d. uniform, since our only operation was to add a matrix that comes from the product of the first row and column, and the measure is additive invariant;

\item The left upper corner of $A$ belongs to $\mathfrak{p}\cap F_0$, which has probability $q^{-1}$; In this case, at least one of the elements in the first row has to be in $\mathfrak{o}^\times$, which has probability $1-q^{2-2n}$. Suppose, without loss of generality, the element on the second column belongs to $\mathfrak{o}^\times$, then we can use this element to eliminate the first two rows and columns, and deal with the right lower $(n-2)\times (n-2)$ corner. Again, the measure of this corner is still i.i.d. uniform, since our only operation was to add a matrix that comes from the product of the first two rows and columns, and the measure is additive invariant.
\end{enumerate}
Therefore by induction, we know 

\begin{align}\begin{split}\mathbf{P}_n^{\her}(0)&=(1-q^{-1})\mathbf{P}_{n-1}^{\her}(0)+q^{-1}(1-q^{2-2n})\mathbf{P}_{n-2}^{\her}(0)\\
&=(1-q^{-1})(1+(-q)^{1-n})\cdots(1-q^{-1})+q^{-1}(1-q^{2-2n})(1+(-q)^{2-n})\cdots(1-q^{-1})\\
&=(1+(-q)^{-n})\cdots(1-q^{-1})
\end{split}\end{align}
\end{enumerate}
This ends our proof.
\end{proof}

\subsection{Hecke ring} \label{subsec:hecke_ring}
The Hecke ring structure plays an essential role throughout our method. See \cite[Chapter V]{Macdonald} for further details. We use the term ``Hecke ring'' in this subsection since it is the original terminology from \cite{Macdonald}, though it is also often called the spherical Hecke algebra in modern works.

Let $G=\GL_n(F)$ be the group of all invertible $n\times n$ matrices over $F$. Also, let $G^+=G\cap M_{n\times n}(\mathfrak{o})$ be the subsemigroup of $G$ consisting of all matrices $x\in G$ with entries $x_{ij}\in\mathfrak{o}$, and let $K=\GL_{n}(\mathfrak{o})=G^+\cap(G^+)^{-1}$ so that $K$ consists of all $x\in G$ with entries $x_{ij}\in\mathfrak{o}$ and $\det(x)$ a unit in $\mathfrak{o}$.
\begin{defi}\label{defi: Hecke algebra}
Let $L(G,K)$ (resp. $L(G^+,K)$) denote the space of all complex-valued continuous functions of compact support on $G$ (resp. $G^+$) which are bivariant with respect to $K$, i.e., such that $f(k_1xk_2)=f(x)$ for all $x\in G$ (resp. $G^+$) and $k_1,k_2\in K$. We define a multiplication on $L(G,K)$ as follows: for $f,g\in L(G,K)$,
$$(f*g)(x)=\int_Gf(xy^{-1})g(y)dy$$
where $dx$ is the unique Haar measure on $G$ such that $K$ has measure $1$. This product is associative and commutative, which makes $L(G^+,K)$ a subring of $L(G,K)$. Also, let $H(G,K)$ (resp. $H(G^+,K)$) denote the subspace of $L(G,K)$ (resp. $L(G^+,K)$) consisting of functions with integer values, which shall be called the \emph{Hecke ring} of $G$ (resp. $G^+$). 
\end{defi}
\begin{prop}
Every $f\in L(G,K)$ (resp. $f\in L(G^+,K)$) could be written as a finite linear combination of functions of the form $c_\mu$, where $\mu=(\mu_1,\ldots,\mu_n)\in\Sig_n$ (resp. $\mu=(\mu_1,\ldots,\mu_n)\in\Sig_n^+$), and $c_\mu$ is the characteristic function of the double coset
$$K\diag_{n\times n}(\pi^{\mu_1},\ldots,\pi^{\mu_n})K.$$
The similar statements for $H(G,K),H(G^+,K)$ translate mutatis mutandis.
\end{prop}

\subsection{Hall-Littlewood polynomials.} For a more complete introduction see \cite[Chapter III]{Macdonald}; some of the condensed treatment below is adapted from \cite[Section 2]{van2020limits}.

\begin{defi}
Let $\Sig_n:=\{\lambda=(\lambda_1,\ldots,\lambda_n)\in\Z^n\mid\lambda_1\ge\ldots\ge\lambda_n\}$ denote the set of integer signatures of length $n$, and let $\Sig_n^+\subset\Sig_n$ denote set of signatures with all parts nonnegative. Given $\l = (\l_1,\ldots,\l_n) \in \Sig_n$, we refer to the integers $\l_i$ as the \emph{parts} of $\l$. We set $|\l| := \sum_{i=1}^n \l_i,n(\lambda):=\sum_{i=1}^n (i-1)\lambda_i$, and $m_k(\l) = \#\{i\mid \l_i = k\}$. If $\lambda\in\Sig_n^+$ only has non-negative parts, we also say $\lambda$ is a \emph{partition} and write $l=l(\lambda)$ as the \emph{length} of $\lambda$, i.e., the number of positive parts. For $\l \in \Sig_n$ and $\mu \in \Sig_{n-1}$, write $\mu \prec_P \l$ if $\l_i \geq \mu_i$ and $\mu_i \geq \l_{i+1}$ for $1 \leq i \leq n-1$. For $\nu \in \Sig_n$ write $\nu \subset \l$ if $\l_i \geq \nu_i$ for $1 \leq i \leq n$, and $\nu \prec_Q \l$ if furthermore $\nu_i \geq \l_{i+1}$ for $1 \leq i \leq n-1$. We write $c[k]$ for the signature $(c,\ldots,c)$ of length $k$, and $(\lambda,\mu)$ for the tuple $(\lambda_1,\ldots,\lambda_n,\mu_1,\ldots,\mu_m)$ when $\lambda\in\Sig_n,\mu\in\Sig_m$. We additionally write $-\lambda=(-\lambda_n,\ldots,-\lambda_1)\in\Sig_n$ for $\lambda\in\Sig_n$.
\end{defi}

\begin{defi}
Let $\lambda=(\l_1,\ldots,\l_n)\in\Sig_n$ be an integer signature of length $n$. For each integer $m\ge 0$, denote $V_m(t)=\frac{(1-t)(1-t^2)\cdots(1-t^m)}{(1-t)^m}$. Then, the \emph{Hall-Littlewood polynomial} $P_\lambda(x_1,\ldots,x_n;t)$ is defined by

\begin{equation}\label{eq: Hall Littlewood polynomial}
P_\l(x_1,\ldots,x_n;t) = \frac{1}{V_\l(t)} \sum_{\sigma \in S_n} \sigma\left(x_1^{\lambda_1}\cdots x_n^{\lambda_n} \prod_{1 \leq i < j \leq n} \frac{x_i-tx_j}{x_i-x_j}\right)
\end{equation}
where $V_\l(t) = \prod_{i \in \Z} V_{m_i(\lambda)}(t)=\frac{1}{(1-t)^n}\prod_{i\in\Z}(t;t)_{m_i(\lambda)}$. Here the notation $(a;q)_n := (1-a)(1-aq) \cdots (1-aq^{n-1})$ for $n \geq 0$, with $(a;q)_0 =1$ and $(a;q)_\infty$ defined in the obvious way. 
\end{defi}

\begin{prop}
The Hall-Littlewood polynomials $P_\l(x_1,\ldots,x_n;t)$ satisfy the following properties:

\begin{enumerate}[left=0pt]
\item They are monic and have the form
$$P_\l(x_1,\ldots,x_n;t) = x_1^{\l_1}x_2^{\l_2}\cdots x_n^{\l_n} + \text{(lower-order monomials in the lexicographic order)};$$

\item When $\lambda\in\Sig_n^+$ ranges over all nonnegative integer signatures of length $n$, the Hall-Littlewood polynomials $P_\l(x_1,\cdots,x_n;t)$ form a $\Z[t]$-basis of $\L_n[t]$, where $\L_n[t]:=\Z[t][x_1,\ldots,x_n]^{S_n}$ is the ring of symmetric polynomials in $n$ variables $x_1,\ldots,x_n$ with coefficients in $\Z[t]$;

\item When $\lambda\in\Sig_n$ ranges over all integer signatures of length $n$, the Hall-Littlewood polynomials $P_\l(x_1,\cdots,x_n;t)$ form a $\Z[t]$ basis of $\Z[t][x_1^{\pm 1},\ldots,x_n^{\pm 1}]^{S_n}$, the ring of symmetric Laurent polynomials in $n$ variables $x_1,\ldots,x_n$ with coefficients in $\Z[t]$.
\end{enumerate}
\end{prop}

Let $\l\in\Sig_n^+$, the dual basis $Q_\l(x_1,\ldots,x_n;t)$ (of the basis $P_\l(x_1,\ldots,x_n;t)$ for the ring $\L_n[t]$ under the natural inner product defined in \cite[Chapter III.4]{Macdonald}) is given by 
\begin{equation}\label{eq:def_Q}
    Q_\l(x_1,\ldots,x_n;t) = \prod_{i\ge 1}(t;t)_{m_i(\lambda)}P_\l(x_1,\ldots,x_n;t).
\end{equation}
Because the $P_\l,\lambda\in\Sig_n^+$ form a basis for the vector space of symmetric polynomials in $n$ variables, there exist symmetric polynomials $P_{\l/\mu}(x_1,\ldots,x_{n-k};t) \in\Lambda_{n-k}[t]$ indexed by $\l \in \Sig_n^+, \mu \in \Sig_k^+$ which are defined by
\begin{equation*}\label{eq: skew Hall Littlewood polynomial}
    P_\l(x_1,\ldots,x_n;t) = \sum_{\mu \in \Sig_k} P_{\l/\mu}(x_{k+1},\ldots,x_n;t) P_\mu(x_1,\ldots,x_k;t).
\end{equation*}
Similarly, for $\l,\nu \in \Sig_n^+$ and $k \geq 1$ arbitrary, define $Q_{\l/\nu}(x_1,\ldots,x_k;t)$ by
\begin{equation*}
    Q_{(\l,0[k])}(x_1,\ldots,x_{n+k};t) = \sum_{\nu \in \Sig_n^+} Q_{\l/\nu}(x_{n+1},\ldots,x_{n+k};t) Q_\nu(x_1,\ldots,x_n;t).
\end{equation*}
In particular, 
\begin{equation}\label{eq:Qs_agree}
    Q_{\l/(0[n])}(x_1,\ldots,x_{n+k};t) = Q_{(\l,0[k])}(x_1,\ldots,x_{n+k};t)
\end{equation}
where the polynomials $P_{\l/\mu},Q_{\l/\mu}$ are named the \emph{skew Hall-Littlewood polynomial}. The following lemma is deduced from (5.11') of \cite[Chapter V]{Macdonald}, which gives an explicit form of these skew Hall-Littlewood polynomials: 
\begin{lemma}\label{lem:branching_formulas}
(Branching rule) For $\l=(\l_1,\ldots,\l_n)\in \Sig_n^+, \mu=(\mu_1,\ldots,\mu_{n-1}) \in \Sig_{n-1}^+$ with $\mu \prec_P \l$, let
\begin{equation}\label{eq: branching coefficient psi}
    \psi_{\l/\mu}(t):= \prod_{\substack{i \in \Z \\ m_i(\mu) = m_i(\l)+1}} (1-t^{m_i(\mu)}).
\end{equation}
For $\nu \in \Sig_n^+$ with $\nu \prec_Q \l$, let
\begin{equation}\label{eq: branching coefficient varphi}\varphi_{\l/\nu}(t)=\prod_{\substack{i \in \Z \\ m_i(\l) = m_i(\nu)+1}} (1-t^{m_i(\l)})\end{equation}
Then for $\l,\nu\in\Sig_n^+,\mu\in\Sig_{n-k}^+$, we have
\begin{equation}\label{eq: branching rule for P}
P_{\l/\mu}(x_1,\ldots,x_k;t)=\sum_{\mu = \l^{(0)} \prec_P \l^{(1)} \prec_P \cdots \prec_P \l^{(k)}= \l} \prod_{i=1}^k x_i^{|\l^{(i)}|-|\l^{(i-1)}|}\psi_{\l^{(i)}/\l^{(i-1)}}(t).
\end{equation} 
and
\begin{equation}\label{eq: branching rule for Q}
 Q_{\l/\nu}(x_1,\ldots,x_k;t) = \sum_{\nu = \l^{(0)} \prec_Q \l^{(1)} \prec_Q \cdots \prec_Q \l^{(k)}=\l} \prod_{i=1}^{k} x_i^{|\l^{(i)}|-|\l^{(i-1)}|}\varphi_{\l^{(i)}/\l^{(i-1)}}(t).
\end{equation}
\end{lemma}

The above formulas inspire us to extend the definition of $P$ and $Q$ to possibly negative signatures. 

\begin{defi}\label{def:skew_gt_patterns}
For $\mu \in \Sig_n, \l \in \Sig_{n+k}$, we define $\GT_P(\l/\mu)$ to be the set of sequences of interlacing signatures $\mu = \l^{(0)} \prec_P \l^{(1)} \prec_P \cdots \prec_P \l^{(k)}= \l$. 

For $\l,\nu \in \Sig_n$, we define $\GT_{Q,k}(\l/\nu)$ to be the set of sequences of length $n$ interlacing signatures $\nu = \l^{(0)} \prec_Q \l^{(1)} \prec_Q \cdots \prec_Q \l^{(k)}=\l$. We refer to elements of either $\GT_P$ or $\GT_{Q,k}$ as \emph{Gelfand-Tsetlin patterns}.

For $T \in \GT_P(\l/\mu)$ with $\len(\l) = \len(\mu)+k$, set $\psi(T) := \prod_{i=1}^{k} \psi_{\l^{(i)}/\l^{(i-1)}}$. For $T \in \GT_{Q,k}(\l/\nu)$, set $\varphi(T) := \prod_{i=1}^{k} \varphi_{\l^{(i)}/\l^{(i-1)}}$. In both cases, let $wt(T) := (|\l^{(1)}|-|\l^{(0)}|,\ldots,|\l^{(k)}|-|\l^{(k-1)}|) \in \Z^k$. 
\end{defi}

From now on, we often write $\bx$ for the collection of variables $x_1,\ldots,x_n$ when $n$ is clear from context. For $\bd \in \Z^n$ we write $\bx^\bd := x_1^{d_1} \cdots x_n^{d_n}$.

\begin{defi}\label{def:skew_signature_functions}
For any $n \geq 0, k \geq 1$ and $\l \in \Sig_{n+k}, \mu \in \Sig_{n}$, we let 
\begin{equation}\label{eq:Pbranch_with_GT}
    P_{\l/\mu}(x_1,\ldots,x_k;t) = \sum_{T \in \GT_P(\l/\mu)} \psi(T) \bx^{wt(T)}
\end{equation}
For any $\nu, \kappa \in \Sig_n$ we let 
\begin{equation}\label{eq:Qbranch_with_GT}
    Q_{\kappa/\nu}(x_1,\ldots,x_k;t) = \sum_{T \in \GT_{Q,k}(\kappa/\nu)} \varphi(T) \bx^{wt(T)}. 
\end{equation}
\end{defi}

Note that these are just the formulas in \Cref{lem:branching_formulas}, with the only change being that we do not require the signatures to be nonnegative.

\begin{rmk}
Notice that for $\lambda,\nu\in\Sig_n$, the skew Hall-Littlewood polynomial $Q_{\lambda/\nu}(x_1,\ldots,x_r;t)$ is a polynomial in the power sums $p_k=x_1^k+\cdots+x_r^k$, which is independent of $r$ for all sufficiently large $r$. This enables us to consider $Q_{\lambda/\nu}(a_1,a_2,\ldots;t)$ with an infinite sequence of real numbers satisfying $\sum_{i\ge 1}|a_i|<\infty$ by setting $p_k(a_1,a_2,\ldots)=\sum_{i\ge 1}a_i^k$, and the definition of $Q_{\lambda/\nu}(a_1,a_2,\ldots;t)$ simultaneously follows. From now on, given $\mu\in\Sig_n^+$, we sometimes briefly write $Q_\mu(a_1,a_2,\ldots;t)$ for $Q_{(\mu/0[n])}(a_1,a_2,\ldots;t)$.
\end{rmk}

\begin{lemma}\label{lem:basic_signature_func_properties}
Let $\l,\nu \in \Sig_{n+k}, \mu \in \Sig_{n}$. Then
\begin{align}
    P_{-\l/-\mu}(x_1,\ldots,x_k;t) &= P_{\l/\mu}(x_1^{-1},\ldots,x_k^{-1};t) \label{eq:p_invert_vars}\\
    Q_{-\l/-\nu}(x_1,\ldots,x_r;t) &= Q_{\l/\nu}(x_1^{-1},\ldots,x_r^{-1};t) \label{eq:q_invert_vars}\\
    P_{(\l + (d[n+k]))/(\mu+d[n])}(x_1,\ldots,x_k;t) &= (x_1\cdots x_k)^d P_{\l/\mu}(x_1,\ldots,x_k;t) \label{eq:p_add_constant}\\
    Q_{(\l+(d[n+k]))/(\nu + (d[n+k]))}(x_1,\ldots,x_r;t) &=  Q_{\l/\nu}(x_1,\ldots,x_r;t). \label{eq:q_add_constant}
\end{align}
\end{lemma}
\begin{proof}
See \cite[Lemma 2.2]{van2020limits}.
\end{proof}

\begin{lemma}\label{lem:asym_cauchy}
(Skew Cauchy identity) Let $\nu \in \Sig_N, \mu \in \Sig_{N+k}$, and $x_1,\ldots,x_k, y_1,\ldots,y_r$ be indeterminates. Then
\begin{multline}\label{eq:asym_cauchy}
    \sum_{\kappa \in \Sig_{N+k}} Q_{\kappa/\mu}(y_1,\ldots,y_r;t) P_{\kappa/\nu}(x_1,\ldots,x_k;t) \\
    = \Pi_t(x_1,\ldots,x_k; y_1,\ldots,y_r)\sum_{\l\in \Sig_N} Q_{\nu/\l}(y_1,\ldots,y_r;t) P_{\mu/\l}(x_1,\ldots,x_k;t) 
\end{multline}
where the Cauchy kernel $\Pi_t(x_1,\ldots,x_k; y_1,\ldots,y_r)$ has the form 
\begin{equation}\label{eq: Cauchy kernel}
    \Pi_t(x_1,\ldots,x_k; y_1,\ldots,y_r) = \prod_{\substack{1 \leq i \leq k \\ 1 \leq j \leq r}} \frac{1-tx_iy_j}{1-x_iy_j}=\exp(\sum_{n\ge1}\frac{1-t^n}{n}\sum_{i=1}^kx_i^n\cdot\sum_{j=1}^ry_j^n)
\end{equation}
and \eqref{eq:asym_cauchy} is interpreted as an equality of formal power series in the variables.
\end{lemma}
\begin{proof}
See \cite[Lemma 2.3]{van2020limits}.
\end{proof}

Recall that the set $\{P_\l(\bx;t): \l \in \Sig_n\}$ forms a basis for the ring of symmetric Laurent polynomials $\L_n[x_1^{\pm 1},\ldots,x_n^{\pm 1}]$. Hence for any $\mu,\nu \in \Sig_n$ one has
\begin{equation}\label{eq:mac_func_struct_coefs}
    P_\mu(x_1,\ldots,x_n;t) \cdot P_\nu(x_1,\ldots,x_n;t) = \sum_{\l \in \Sig_n} c_{\mu,\nu}^\l(t) P_\l(x_1,\ldots,x_n;t)
\end{equation}
for the \emph{Littlewood-Richardson coefficients} $c_{\mu,\nu}^\l(t)\in\Z[t]$. By comparing the degrees on both sides, it is clear that $c_{\mu,\nu}^\l(t)\ne 0$ only when $|\l|+|\mu|=|\nu|$. The Littlewood-Richardson coefficients are related to the skew-$Q$ polynomials as follows.

\begin{prop}\label{prop:q_coproduct_coefs}
Let $m,n \in \N$ and $\l,\nu \in \Sig_n$. Then
\begin{equation*}
    \sum_{\mu \in \Sig_n^+} c_{\nu,\mu}^\l(t) Q_{(\mu/0[n])}(x_1,\ldots,x_m;t)=Q_{\l/\nu}(x_1,\ldots,x_m;t) .
\end{equation*}
\end{prop}

When all signatures are nonnegative this follows by specializing the corresponding statement for symmetric functions, see \cite[Chapter III.5]{Macdonald} where \Cref{prop:q_coproduct_coefs} is taken as the definition of the skew $Q$ polynomials. The case of general signatures follows by adding $d[n],d>0$ sufficiently large to both $\lambda$ and $\nu$, for which feasibility has already proved in \Cref{lem:basic_signature_func_properties}.

In the end, we give the value of the Hall-Littlewood polynomial under principal specialization, which serves as a corollary of \eqref{eq: Hall Littlewood polynomial}: 

\begin{prop}[Principal specialization formulas]\label{prop:hl_principal_formulas}
For $J,n \geq 1$ and $\l \in \Sig_n^+$,
\begin{align}\label{eq: Principal specialization formulas}
\begin{split}
    P_\l(x,xt,\ldots,xt^{n-1};t) &= x^{|\l|} t^{n(\l)} \frac{(t;t)_n}{\prod_{i \in \Z} (t;t)_{m_i(\l)}}=\frac{V_n(t)}{V_\lambda(t)}x^{|\l|} t^{n(\l)}\\
    Q_{\l/(0[n])}(x,xt,\ldots,xt^{J-1};t) &= x^{|\l|} t^{n(\l)} \frac{(t;t)_J}{(t;t)_{m_0(\l)+J-n}} \bbone(m_0(\l)+J-n \geq 0)
\end{split}
\end{align}
\end{prop}

Note that the principal specialization formula for $Q$ differs from the statement in \cite[Ch. III.2, Ex. 1]{Macdonald} due to our conventions on signatures, but it may be derived directly from that statement using \eqref{eq:Qs_agree} to translate between skew and non-skew $Q$ polynomials.

\begin{defi}\label{defi: specialization}
Fix the parameter $t$, we say the specialization $\theta=(a_1,a_2,\ldots)$ is \emph{non-negative} if it takes nonnegative values on the skew Hall-Littlewood symmetric functions $Q_{\lambda/\nu}(a_1,a_2,\ldots;t)\ge 0$ for any $n\ge 1$ and integer signatures $\nu\prec_Q\lambda\in\Sig_n$.
\end{defi}

\begin{defi}
Let $\theta=(a_1,\ldots,a_n),\psi=(b_1,b_2,\ldots)$ be non-negative specializations that satisfy 
$$\Pi_t(\theta;\psi)=\sum_{\lambda\in\Sig_n^+} P_\lambda(\theta;t)Q_\lambda(\psi;t)<\infty.$$ The \emph{Hall-Littlewood measure} with specializations $\theta,\psi$ is the measure on $\Sig_n^+$ given by
$$\mathbf{P}(\lambda)=\frac{P_\lambda(\theta;t)Q_\lambda(\psi;t)}{\Pi_t(\theta;\psi)}.$$
\end{defi}

The \emph{Hall-Littlewood process}, originally defined in \cite{borodin2014macdonald}, describes measures over sequences of partitions by Hall-Littlewood polynomials. Due to our conventions with integer signatures, we must define two slightly different versions:

\begin{defi}
\begin{enumerate}
\item (Corresponding to the corners process) Let $n\ge 1$. For $i=1,2,\ldots,n$, let $b_i$ be a real number, $\theta_i$ be a sequence (finite or infinite) of real numbers $a_j^{(i)}$. The \emph{Hall-Littlewood process} is the probability measure on sequence $$(\lambda,\nu)=\lambda^{(1)}\succ _Q\nu^{(1)}\prec_P\lambda^{(2)}\succ _Q\cdots\succ _Q\nu^{(n-1)}\prec_P\lambda^{(n)},\lambda^{(i)}\in\Sig_i^+,\nu^{(j)}\in\Sig_j^+$$
satisfying
$$\textbf{P}(\lambda,\nu)=\frac{P_{\lambda^{(1)}}(b_0;t)Q_{\lambda^{(1)}/\mu^{(1)}}(\theta_1;t)\cdots P_{\lambda^{(N)}/\mu^{(N-1)}}(b_{n-1};t)Q_{\lambda^{(N)}}(\theta_n;t)}{\prod_{0\le i<j\le n}\Pi_t(b_i;\theta_j)}$$
Here the specializations $b_i,\theta_j$ have to make sense, i.e., the denominator $\prod_{0\le i<j\le n}\Pi_t(b_i;\theta_j)<\infty$ is finite, and the measure is always non-negative.

\item (Corresponding to the product process) Let $n,k\ge 1$. For $i=0,1,\ldots,k$, let $\psi_i$ be a sequence (finite or infinite) of real numbers $a_j^{(i)}$, and $\theta=(a_1,\ldots,a_n)$. The \emph{Hall-Littlewood process} is the probability measure on the sequence
$$\lambda=\lambda^{(0)}\subset\lambda^{(1)}\subset\cdots\subset\lambda^{(k)},\lambda^{(i)}\in\Sig_n^+$$
satisfying
$$\textbf{P}(\lambda)=\frac{Q_{\lambda^{(0)}}(\psi_0;t)Q_{\lambda^{(1)}/\lambda^{(0)}}(\psi_1;t)\cdots Q_{\lambda^{(k)}/\lambda^{(k-1)}}(\psi_k;t)P_{\lambda^{(k)}}(\theta;t)}{\prod_{0\le i\le k}\Pi_t(\psi_i;\theta)}$$
Here the specializations have to make sense, i.e., the denominator $\prod_{0\le i\le k}\Pi_t(\psi_i;\theta)<\infty$ is finite, and the measure is always non-negative.\end{enumerate}
\end{defi}

\subsection{Mixed Littlewood-Richardson coefficients and their relation to Hall-Littlewood processes.} While the coefficients $c_{\mu,\nu}^\l(t)$ in \eqref{eq:mac_func_struct_coefs} are standard, it will be useful for us to define two different versions which are not. These are relevant to the alternating and Hermitian cases respectively in random matrix theory, which explains our notation below, although the definitions themselves do not require any reference to random matrix problems.

\begin{defi}\label{def:nonstandard_LR_coefs}
\begin{enumerate}
\item (Alternating case) Given $\mu\in\Sig_{2n},\nu\in\Sig_n$, we define the \emph{alternating Littlewood-Richardson coefficient} $c_{\mu,\nu}^{\alt,\lambda}(t)\in\Z[t]$ given by
\begin{equation}\label{eq: alt_LR coef}
P_\mu(x_1,x_1t,\ldots,x_n,x_nt;t)P_\nu(x_1,\ldots,x_n;t^2)=\sum_{\lambda\in\Sig_n} c_{\mu,\nu}^{\alt,\lambda}(t)P_\lambda(x_1,\ldots,x_n;t^2)
\end{equation}
\item (Hermitian case) Given $\mu,\nu\in\Sig_n$, we define the \emph{Hermitian Littlewood-Richardson coefficient} $c_{\mu,\nu}^{\her,\lambda}(t)\in\Z[t]$ given by
\begin{equation}\label{eq: her_LR coef}
P_\mu(x_1^2,\ldots,x_n^2;t^2)P_{\nu}(x_1,\ldots,x_n;-t)=\sum_{\lambda\in\Sig_n} c_{\mu,\nu}^{\her,\lambda}(t)P_{\lambda}(x_1,\ldots,x_n;-t)
\end{equation}
\end{enumerate}
\end{defi}

\begin{rmk}\label{rmk:lr_poly}
The coefficient $c_{\mu,\nu}^{\alt,\lambda}(t)$ is a polynomial in $\Z[t]$, because the left hand side of \eqref{eq: alt_LR coef} is in $\Lambda_n[t]$, and $P_\lambda\in\Lambda_n[t]$ are symmetric Laurent polynomials that are monic under the lexicographical ordering. Similarly, we can prove that $c_{\mu,\nu}^{\her,\lambda}(t)\in \Z[t]$ has integer coefficients.
\end{rmk}

Similarly to classical Littlewood-Richardson coefficients, these also appear in coproduct operations similarly to \Cref{prop:q_coproduct_coefs}, as the result below shows.

\begin{prop}\label{prop: skew and LR coef}
Given $m,n\ge 0$ and $\l,\nu\in\Sig_n$, we have
\begin{equation}\label{eq: alt_skew and LR coef}\sum_{\mu\in\Sig_{2n}^+} c_{\mu,\nu}^{\alt,\lambda}(t)Q_{(\mu/0[2n])}(x_1,\ldots,x_m;t)=Q_{\lambda/\nu}(x_1,\ldots,x_m;t^2),\end{equation}
and
\begin{equation}\label{eq: her_skew and LR coef}\sum_{\mu\in\Sig_n^+} c_{\mu,\nu}^{\her,\lambda}(t)Q_{(\mu/0[n])}(x_1^2,\ldots,x_m^2;t^2)=Q_{\lambda/\nu}(x_1,-x_1,\ldots,x_m,-x_m;-t).\end{equation}
\end{prop}
\begin{proof}
For the alternating case, we only need to prove the equality
\begin{multline}\label{eq: alt_y assistant}\sum_{\mu\in\Sig_{2n}^+,\lambda\in\Sig_n} c_{\mu,\nu}^{\alt,\lambda}(t)P_\lambda(y_1,\ldots,y_n;t^2)Q_{(\mu/0[2n])}(x_1,\ldots,x_m;t)=\\
\sum_{\lambda\in\Sig_n} P_\lambda(y_1,\ldots,y_n;t^2)Q_{\lambda/\nu}(x_1,\ldots,x_m;t^2)\end{multline}
must hold, since the $P_\la$ form a basis for symmetric Laurent polynomials in $y_1,\ldots,y_n$ and we may thus equate coefficients on both sides. Now, applying the skew Cauchy identity from \Cref{lem:asym_cauchy} in the second equality below, we have

\begin{align}\begin{split}
\text{LHS of \eqref{eq: alt_y assistant}}&=\sum_{\mu\in\Sig_{2n}^+}P_\mu(y_1,y_1t,\ldots,y_n,y_nt;t)Q_{(\mu/0[2n])}(x_1,\ldots,x_m;t)P_\nu(y_1,\ldots,y_n;t^2)\\
&=\exp(\sum_{n\ge 1}\frac{1-t^n}{n}(1+t^n)(y_1^n+\cdots+y_n^n)(x_1^n+\cdots+x_m^n))P_\nu(y_1,\ldots,y_n;t^2)\\
&=\exp(\sum_{n\ge 1}\frac{1-t^{2n}}{n}(y_1^n+\cdots+y_n^n)(x_1^n+\cdots+x_m^n))P_\nu(y_1,\ldots,y_n;t^2)\\
&=\text{RHS of \eqref{eq: alt_y assistant}}.
\end{split}\end{align}
Likewise, for the Hermitian case, we only need to prove the equality
\begin{multline}\label{eq: her_y assistant}\sum_{\mu\in\Sig_n^+,\lambda\in\Sig_n} c_{\mu,\nu}^{\her,\lambda}(t)P_\lambda(y_1,\ldots,y_n;-t)Q_{(\mu/0[n])}(x_1^2,\ldots,x_m^2;t^2)=\\
\sum_{\lambda\in\Sig_n}P_\lambda(y_1,\ldots,y_n;-t)Q_{\lambda/\nu}(x_1,-x_1,\ldots,x_m,-x_m;-t)\end{multline}
must hold. In fact, apply the skew Cauchy identity from \Cref{lem:asym_cauchy}, and we have

\begin{align}\begin{split}
\text{LHS of \eqref{eq: her_y assistant}}&=\sum_{\mu\in\Sig_n^+}P_\mu(y_1^2,\ldots,y_n^2;t^2)Q_{(\mu/0[n])}(x_1^2,\ldots,x_m^2;t^2)P_\nu(y_1,\ldots,y_n;-t)\\
&=\exp(\sum_{n\ge 1}\frac{1-t^{2n}}{n}(y_1^{2n}+\cdots+y_n^{2n})(x_1^{2n}+\cdots+x_m^{2n}))P_\nu(y_1,\ldots,y_n;-t)\\
&=\text{RHS of \eqref{eq: her_y assistant}}.
\end{split}\end{align}
This ends the proof.
\end{proof}

The following definition is motivated by \Cref{prop: skew and LR coef}, which will be useful in further discussion.

\begin{defi}\label{defi: change of specialization}Let $\theta=(a_1,\ldots,a_n),\psi=(b_1^2,b_2^2,\ldots)$ be specializations ($\psi$ is a finite or infinite sequence) and $t$ a parameter. Also, we write $p_k$ as the $k$th power sum.
\begin{enumerate}
\item(Alternating case) We denote the specialization $\theta'=(a_1,a_1t,\ldots,a_n,a_nt)$. This makes the below equality hold:
$$p_k(\theta')=(1+t^k)p_k(\theta),\quad\forall k\ge 1.$$
\item(Hermitian case) We denote the specializations $\theta'=(a_1^2,\ldots,a_n^2),\psi^*=(b_1,-b_1,b_2,-b_2,\ldots)$. This makes the below equality hold: 
$$p_k(\theta')=p_{2k}(\theta),p_{2k-1}(\psi^*)=0,p_{2k}(\psi^*)=2p_k(\psi),\quad\forall k\ge 1.$$
\end{enumerate}
\end{defi}

The below product convolution is related to the random matrix products but is more general. For simplicity, we give only the definition for specializations of the form $\theta=(a_1,\ldots,a_n)$.

\begin{defi}\label{defi: product convolution}
\begin{enumerate}[left=0pt]
\item (Alternating case) Given $\mu\in\Sig_{2n},\nu\in\Sig_n$, we define a random signature $\mu\boxtimes_\theta^{\alt}\nu\in\Sig_n$ by

$$\mathbf{P}(\mu\boxtimes_\theta^{\alt}\nu=\lambda)=\frac{c_{\mu,\nu}^{\alt,\lambda}(t)P_\lambda(\theta;t^2)}{P_\mu(\theta';t)P_\nu(\theta;t^2)}$$
where $c_{\mu,\nu}^{\alt,\lambda}(t)\in\Z[t]$ is as in \Cref{def:nonstandard_LR_coefs}.



\item (Hermitian case) Given $\mu\in\Sig_n,\nu\in\Sig_n$, we define a random signature $\mu\boxtimes_\theta^{\her}\nu\in\Sig_n$ by

$$\mathbf{P}(\mu\boxtimes_\theta^{\her}\nu=\lambda)=\frac{c_{\mu,\nu}^{\her,\lambda}(t)P_\lambda(\theta;-t)}{P_\mu(\theta';t^2)P_\nu(\theta;-t)}$$
where $c_{\mu,\nu}^{\her,\lambda}(t)\in\Z[t]$ is as in \Cref{def:nonstandard_LR_coefs}.

\end{enumerate}
\end{defi}

The product convolution operations are related to the Hall-Littlewood process as follows. 

\begin{thm}\label{thm: product convolution and HL process}
\begin{enumerate}[left=0pt]
\item(Alternating case)Let $\theta=(a_1,\ldots,a_n), \psi,\psi_1,\ldots,\psi_k$ be sequences (finite or infinite) of real numbers $b_j,b_j^{(1)},\ldots,b_j^{(k)}$ such that 

\begin{enumerate}[left=0pt]
\item For all integer signatures $\nu\prec_Q\lambda\in\Sig_n$, we have 

$$Q_{\lambda/\nu}(\psi;t^2),Q_{\lambda/\nu}(\psi_1;t^2),\ldots,Q_{\lambda/\nu}(\psi_k;t^2)\ge 0;$$

\item$\Pi_t(\theta';\psi)=\Pi_{t^2}(\theta;\psi),\Pi_t(\theta';\psi_1)=\Pi_{t^2}(\theta;\psi_1),\ldots,\Pi_t(\theta';\psi_k)=\Pi_{t^2}(\theta;\psi_k)<\infty$. Here $\theta'$ follows from \Cref{defi: change of specialization}.
\end{enumerate} Suppose $\nu\in\Sig_n^+$ is distributed by the Hall-Littlewood measure of indeterminant $t^2$ with specializations $\theta,\psi$, and for all $1\le i\le k$, $\mu_i\in\Sig_{2n}^+$ is distributed by the Hall-Littlewood measure of indeterminant $t$ with specializations $\theta',\psi_i$ for each $i=1,2,\ldots,k$. Then for every fixed $\lambda,\lambda^{(1)},\ldots,\lambda^{(k)}\in\Sig_n^+$, we have the probability

$$\mathbf{P}(\mu_\tau\boxtimes_\theta\ldots\boxtimes_\theta\mu_1\boxtimes_\theta^{\alt}\nu=\lambda^{(\tau)},\quad\forall 0\le\tau\le k)$$
(here $\lambda^{(0)}=\lambda$) is equal to

\begin{equation}\label{eq: alt_joint distribution of product convolution}
\frac{Q_{\lambda}(\psi;t^2)Q_{\lambda^{(1)}/\lambda}(\psi_1;t^2)\cdots Q_{\lambda^{(k)}/\lambda^{(k-1)}}(\psi_k;t^2)P_{\lambda^{(k)}}(\theta;t^2)}{\Pi_{t^2}(\theta;\psi)\Pi_{t^2}(\theta;\psi_1)\cdots\Pi_{t^2}(\theta;\psi_k)}
\end{equation}

\item(Hermitian case)Let $\theta=(a_1,\ldots,a_n),\psi,\psi_1,\ldots,\psi_k$ be sequences (finite or infinite) of square of real numbers $b_j^2,(b_j^{(1)})^2,\ldots,(b_j^{(k)})^2$ such that
\begin{enumerate}
\item For all integer signatures $\nu\prec_Q\lambda\in\Sig_n$ such that $|\lambda|-|\nu|$ an even integer, we have $$Q_{\lambda/\nu}(\psi;-t),Q_{\lambda/\nu}(\psi_1^*;-t),\ldots,Q_{\lambda/\nu}(\psi_k^*;-t)\ge 0;$$

\item $\Pi_{-t}(\theta;\psi),\Pi_{t^2}(\theta';\psi_1)=\Pi_{-t}(\theta;\psi_1^*),\ldots,\Pi_{t^2}(\theta';\psi_k)=\Pi_{-t}(\theta;\psi_k^*)<\infty$. Here $\theta',\psi_1^*,\ldots,\psi_k^*$ follows from \Cref{defi: change of specialization}.
\end{enumerate} Suppose $\nu\in\Sig_n^+$ is distributed by the Hall-Littlewood measure of indeterminant $-t$ with specializations $\theta,\psi$, and for all $1\le i\le k$, $\mu_i\in\Sig_n^+$ is distributed by the Hall-Littlewood measure of indeterminant $t^2$ with specializations $\theta',\psi_i$ for each $i=1,2,\ldots,k$. Then for every fixed $\lambda,\lambda^{(1)},\ldots,\lambda^{(k)}\in\Sig_n^+$, we have the probability

$$\mathbf{P}(\mu_\tau\boxtimes_\theta\ldots\boxtimes_\theta\mu_1\boxtimes_\theta^{\her}\nu=\lambda^{(\tau)},\quad\forall 0\le\tau\le k)$$
(here $\lambda^{(0)}=\lambda$) is equal to

\begin{equation}\label{eq: her_joint distribution of product convolution}
\frac{Q_\lambda(\psi;-t)Q_{\lambda^{(1)}/\lambda}(\psi_1^*;-t)\cdots Q_{\lambda^{(k)}/\lambda^{(k-1)}}(\psi_k^*;-t)P_{\lambda^{(k)}}(\theta;-t)}{\Pi_{-t}(\theta;\psi)\Pi_{-t}(\theta;\psi_1^*)\cdots\Pi_{-t}(\theta;\psi_k^*)}
\end{equation}
\end{enumerate}
\end{thm}

\begin{proof}
We begin with the alternating case. We only need to prove the case $k=1$, and the larger case is similar. Let $\lambda\in\Sig_n$ be any fixed integer signature. By the definition of product convolution, we have
\begin{align}\label{eq: alt_conditional distribution of HL process}\begin{split}
\mathbf{P}(\mu_1\boxtimes_\theta^{\alt}\lambda=\lambda^{(1)})&=\sum_{\mu}\frac{P_\mu(\theta';t)Q_\mu(\psi_1;t)}{\Pi_{t}(\theta';\psi_1)}\frac{c_{\mu,\lambda}^{\alt,\lambda^{(1)}}(t)P_{\lambda^{(1)}}(\theta;t^2)}{P_\mu(\theta';t)P_{\lambda}(\theta;t^2)}\\
&=\frac{P_{\lambda^{(1)}}(\theta;t^2)\sum_\mu c_{\mu,\lambda}^{\alt,\lambda^{(1)}}(t)Q_\mu(\psi_1;t)} {P_\lambda(\theta;t^2)\Pi_{t^2}(\theta;\psi_1)}\\
&=\frac{Q_{\lambda^{(1)}/\lambda}(\psi_1;t^2)P_{\lambda^{(1)}}(\theta;t^2)}{P_\lambda(\theta;t^2)\Pi_{t^2}(\theta;\psi_1)}.
\end{split}\end{align}
Here the last row is deduced from \eqref{eq: alt_skew and LR coef}. This implies
\begin{align*}
\mathbf{P}(\nu=\lambda,\mu_1\boxtimes_\theta^{\alt}\nu=\lambda^{(1)})&=\frac{P_\lambda(\theta;t^2)Q_\lambda(\psi;t^2)}{\Pi_{t^2}(\theta;\psi)}\mathbf{P}(\mu_1\boxtimes_\theta^{\alt}\lambda=\lambda^{(1)})\\
&=\frac{Q_\lambda(\psi;t^2)Q_{\lambda^{(1)}/\lambda}(\psi_1;t^2)P_{\lambda^{(1)}}(\theta;t^2)}{\Pi_{t^2}(\theta;\psi)\Pi_{t^2}(\theta;\psi_1)}.
\end{align*}
Likewise, for the Hermitian case, we only need to prove the case $k=1$, and the larger case is similar. Let $\lambda\in\Sig_n$ be any fixed integer signature. By the definition of the product convolution and \eqref{eq: her_skew and LR coef},
\begin{align}\begin{split}
\mathbf{P}(\mu_1\boxtimes_\theta^{\her}\lambda=\lambda^{(1)})&=\sum_{\mu}\frac{P_\mu(\theta';t^2)Q_\mu(\psi_1;t^2)}{\Pi_{t^2}(\theta';\psi_1)}\frac{c_{\mu,\lambda}^{\her,\lambda^{(1)}}(t)P_{\lambda^{(1)}}(\theta;-t)}{P_\mu(\theta';t^2)P_{\lambda}(\theta;-t)}\\
&=\frac{P_{\lambda^{(1)}}(\theta;-t)\sum_\mu c_{\mu,\lambda}^{\her,\lambda^{(1)}}(t)Q_\mu(\psi_1;t^2)} {P_\lambda(\theta;-t)\Pi_{-t}(\theta;\psi_1^*)}\\
&=\frac{Q_{\lambda^{(1)}/\lambda}(\psi_1^*;-t)P_{\lambda^{(1)}}(\theta;-t)}{P_\lambda(\theta;-t)\Pi_{-t}(\theta;\psi_1^*)}
\end{split}\end{align}
Here the last row is deduced from \eqref{eq: her_skew and LR coef}. This implies
\begin{align*}
\mathbf{P}(\nu=\lambda,\mu_1\boxtimes_\theta^{\her}\nu=\lambda^{(1)})&=\frac{P_\lambda(\theta;-t)Q_\lambda(\psi;-t)}{\Pi_{-t}(\theta;\psi)}\mathbf{P}(\mu_1\boxtimes_\theta^{\alt}\lambda=\lambda^{(1)})\\
&=\frac{Q_\lambda(\psi;-t)Q_{\lambda^{(1)}/\lambda}(\psi_1^*;-t)P_{\lambda^{(1)}}(\theta;-t)}{\Pi_{-t}(\theta;\psi)\Pi_{-t}(\theta;\psi_1^*)}.
\end{align*}
Hence ends the proof.
\end{proof}

There are some antecedents to \Cref{thm: product convolution and HL process}. The operations of \Cref{defi: product convolution} are analogues of an operation $\boxtimes_{\mathbf{a}}$ defined similarly in \cite{van2020limits} in terms of the usual structure coefficients $c_{\mu,\nu}^\la(t)$, and it was shown in \cite[Proposition 2.6]{van2020limits} that these yield Hall-Littlewood processes in a similar manner. It is also worth noting that structurally identical operations for other degenerations of Macdonald polynomials appeared in classical (complex) random matrix theory in e.g. \cite{gorin2020crystallization,ahn2019fluctuations,borodin2018product,gorin2018gaussian}. 

Nonetheless, we find \Cref{thm: product convolution and HL process} quite surprising: unlike the above-mentioned results, the operations $\boxtimes_{\mathbf{a}}^{\alt}$ and $\boxtimes_{\mathbf{a}}^{\her}$ feature Hall-Littlewood polynomials with not only different variables, but different parameter $t$. Yet somehow, through the above cancellations, these operations yield usual bona fide Hall-Littlewood processes with only one parameter $t^2$ (resp. $-t$) for the alternating (resp. Hermitian) case.



\section{The Alternating Hecke Module}\label{sec:alt_hecke}

\subsection{The alternating Hecke module $H(G^{\alt},K)$}

Let $G=\GL_{2n}(F)$ be the group of all invertible $2n\times 2n$ matrices over $F$. Also, let $G^+=G\cap M_{2n\times 2n}(\mathfrak{o})$ be the subsemigroup of $G$ consisting of all matrices $x\in G$ with entries $x_{ij}\in\mathfrak{o}$, and let $K=\GL_{2n}(\mathfrak{o})=G^+\cap(G^+)^{-1}$ so that $K$ consists of all $x\in G$ with entries $x_{ij}\in\mathfrak{o}$ and $\det(x)$ a unit in $\mathfrak{o}$. The structure of the rings $L(G,K),L(G^+,K),H(G,K),H(G^+,K)$ follows from \Cref{defi: Hecke algebra}.

\begin{defi}
Let $L(G^{\alt},K)$ (resp. $L(G^{+{\alt}},K)$) denote the space of all complex-valued continuous functions of compact support on $G^{\alt}=G\cap\Alt_{2n}(F)$ (resp. $G^{+{\alt}}=G^+\cap\Alt_{2n}(\mathfrak{o})$) which are invariant with respect to $K$, i.e., such that

$$g(k^Txk)=g(x)$$

for all $x\in G$ (resp. $x\in G^+$) and $k\in K$. We may regard $L(G^{+{\alt}},K)$ as a submodule of $L(G^{\alt},K)$. 
\end{defi}

We define a multiplication of $L(G,K)$ over the module $L(G^{\alt},K)$ as follows: for $f\in L(G,K)$, $g\in L(G^{\alt},K)$,

$$(f*g)(x):=\int_Gf(y)g(y^{-1}xy^{-T})dy$$

Since $f$ and $g$ are compactly supported, the integration is over a compact set. This product satisfies $f_1*(f_2*g)=(f_1*f_2)*g$ and hence defines a module structure. (We use the same notation $*$ for the convolution of the Hecke algebra itself and the convolution of the Hecke algebra over the alternating Hecke module, but it is easy to distinguish these two based on context.) 

\quad

Each function $g\in L(G^{\alt},K)$ is constant on the orbits $\{k^Txk|k\in K\}$ in $G^{\alt}$. These orbits are compact, open, and mutually disjoint. Since $g$ has compact support, it follows that $g$ takes non-zero values on only finitely many orbits $\{k^Txk|k\in K\}$, and hence can be written as a finite linear combination of their characteristic functions. Therefore, the characteristic functions of these orbits in $G^{\alt}$ form a $\C$-basis of $L(G^{\alt},K)$. 

If we vary the definition of the module $L(G^{\alt},K)$ (resp. $L(G^{+{\alt}},K)$) by requiring the functions to take their values in $\Z$ instead of $\C$, the resulting module is the generalization of the Hecke module over $G^{\alt}$ (resp. $G^{+{\alt}}$), and we denote it by $H(G^{\alt},K)$ (resp. $H(G^{+{\alt}},K)$). Clearly we have

$$L(G^{\alt},K)=H(G^{\alt},K)\otimes_\Z\C,L(G^{+{\alt}},K)=H(G^{+{\alt}},K)\otimes_\Z\C$$

Consider an orbit $\{k^Txk|k\in K\}$, where $x\in G^{\alt}$. By multiplying $x$ by a suitable power of $\pi$ (the generator of $\mathfrak{p}$) we can bring $x$ into $G^{+{\alt}}$. Also, \Cref{prop:smith} implies that each orbit $\{k^Txk|k\in K\}$ has a unique representative of the form

$$\pi^{\alt}_\lambda=\text{diag}({\begin{pmatrix} 0 & \pi^{\lambda_1} \\ -\pi^{\lambda_1} & 0\end{pmatrix}},\ldots,{\begin{pmatrix} 0 & \pi^{\lambda_n} \\ -\pi^{\lambda_n} & 0\end{pmatrix}})$$
where $\lambda_1\ge\lambda_2\ge\ldots\ge\lambda_n$. We have $\lambda_n\ge 0$ if and only if $x\in G^{+{\alt}}$.

Let $c^{\alt}_\lambda$ denote the characteristic function of the orbit $\{k^T\pi^{\alt}_\lambda k|k\in K\}$. Then we have the $c^{\alt}_\lambda$ (resp. the $c^{\alt}_\lambda$ such that $\lambda_n\ge 0$) form a $\Z$-basis of $H(G^{\alt},K)$(resp. $H(G^{+{\alt}},K)$). The characteristic function $c_0$ of $K$ is the identity element of $H(G,K)$ and $H(G^+,K)$, and this also plays the role of the identity element when multiplying with elements in $H(G^{\alt},K)$ and $H(G^{+{\alt}},K)$. 

\begin{defi}
Let $\mu\in\Sig_{2n}$, $\nu\in\Sig_n$. Then, we define the structure coefficient $G_{\mu,\nu}^{\alt,\lambda}(\mathfrak{o})$ for the expansion of the product $c_\mu *c^{\alt}_\nu$:
\begin{equation}\label{eq:structure coefficient_alt}c_\mu *c^{\alt}_\nu=\sum_\lambda G_{\mu,\nu}^{\alt,\lambda}(\mathfrak{o})c^{\alt}_\lambda
\end{equation}
\end{defi}

\begin{prop}\label{prop: Hall alternating polynomial}
$G_{\mu,\nu}^{\alt,\lambda}(\mathfrak{o})$ has the following properties:
\begin{enumerate}[left=0pt]
\item $G_{\mu,\nu}^{\alt,\lambda}(\mathfrak{o})=0$ unless $|\lambda|=|\mu|+|\nu|$. In this case, $G_{\mu,\nu}^{\alt,\lambda}(\mathfrak{o})\in\Z_{\ge 0}$ is a non-negative integer; \label{item: alt_nonnegative integer}

\item $G_{\mu,\nu}^{\alt,\lambda}(\mathfrak{o})=G_{\mu,N[n]+\nu}^{\alt,N[n]+\lambda}(\mathfrak{o})=G_{N[2n]+\mu,\nu}^{\alt,2N[n]+\lambda}(\mathfrak{o})$ for all $N\in\Z$;

\end{enumerate}
\end{prop}
\begin{proof}
\begin{enumerate}[left=0pt]
\item We now give an interpretation of this coefficient $G_{\mu,\nu}^{\alt,\lambda}(\mathfrak{o})$. Notice that
$$G_{\mu,\nu}^{\alt,\lambda}(\mathfrak{o})=(c_\mu*c_\nu^{\alt})(\pi_\lambda^{\alt})=\int_G c_\mu(y)c_\nu^{\alt}(y^{-1}\pi_\lambda^{\alt}y^{-T})dy$$
Since $c_\mu(y)$ vanishes for $y$ outside $K\pi_\mu K$, the integration is over this orbit, which we shall write as a disjoint union of left cosets, say
\begin{equation}\label{eqref: alt_disjoint left coset}K\pi_\mu K=\bigsqcup_j y_jK\quad (y_j\in K\pi_\mu)\end{equation}
Therefore, we have
$$G_{\mu,\nu}^{\alt,\lambda}(\mathfrak{o})=\sum_j\int_{y_jK}c^{\alt}_\nu(y^{-1}\pi_\lambda^{\alt} y^{-T})dy=\sum_jc^{\alt}_\nu(y_j^{-1}\pi_\lambda^{\alt} y_j^{-T})$$
since $K$ has measure $1$. Hence $G_{\mu,\nu}^{\alt,\lambda}(\mathfrak{o})$ is equal to the number of $j$ such that $y_j^{-1}\pi^{\alt}_\lambda y_j^{-T}\in\{k\pi^{\alt}_\nu k^T|k\in K\}$, which is a non-negative integer; If $|\lambda|\ne|\mu|+|\nu|$, by checking the absolute value of the determinant we know such $y_j$ does not exist, and therefore $G_{\mu,\nu}^{\alt,\lambda}(\mathfrak{o})=0$;

\item By multiplying $\pi^N$ to the alternating matrix in the middle, we obtain $G_{\mu,\nu}^{\alt,\lambda}(\mathfrak{o})=G_{\mu,N[n]+\nu}^{\alt,N[n]+\lambda}(\mathfrak{o})$; And by multiplying $\pi^N$ to the two matrices on both sides, we obtain $G_{\mu,\nu}^{\alt,\lambda}(\mathfrak{o})=G_{N[2n]+\mu,\nu}^{\alt,2N[n]+\lambda}(\mathfrak{o})$;
\end{enumerate}
\end{proof}


\begin{rmk}
It is worth mentioning that the order of the subscripts cannot be changed, i.e., $G_{\mu,\nu}^{\alt,\lambda}(\mathfrak{o})$ is not the same as $G_{\nu,\mu}^{\alt,\lambda}(\mathfrak{o})$. 
\end{rmk}

The setting above helps us study the distribution of singular numbers of matrix products. 

\begin{prop}\label{prop: alt_transition step}
For all $\lambda\in\Sig_n$, let $$V(K^T\pi_\lambda^{\alt}K)=\int_{G^{\alt}}c_\lambda^{\alt}(x)dx$$ 
denote the volume of the orbit $\{k^T\pi_\lambda^{\alt}k|k\in K\}$, where $dx$ is the $G$-invariant measure on $G^{\alt}$ normalized by $\int_{K^TK}dx=1$. Then the probability $\mathbf{P}_{\mu,\nu}^{\alt,\lambda}:=\mathbf{P}(\SN^{\alt}(B^TAB)=\lambda)$ from \Cref{thm:Product process} has the form
$$\mathbf{P}_{\mu,\nu}^{\alt,\lambda}=\frac{G_{\mu,\nu}^{\alt,\lambda}(\mathfrak{o})V(K^T\pi_\lambda^{\alt}K)}{N_\mu V(K^T\pi_\nu^{\alt}K)}$$
where $N_\mu$ is the number of $y_j$ in \eqref{eqref: alt_disjoint left coset}, i.e. the number of disjoint left cosets in $K\pi_\mu K$.
\end{prop}

\begin{proof}
Consider the integral 
$$\mathcal{I}=\int_{G^{\alt}\times G}c_\nu^{\alt}(x)c_\mu(y)c_\lambda^{\alt}(y^Txy)dxdy.$$
On one hand, we have
\begin{equation}\label{eq: alt_x and y integral}\mathcal{I}=V(K^T\pi_\nu^{\alt}K)\int_G c_\mu(y)c_\lambda^{\alt}(y^T\pi_\nu^{\alt}y)dy=\mathbf{P}_{\mu,\nu}^{\alt,\lambda}N_\mu V(K^T\pi_\nu^{\alt}K).\end{equation}
The second equality holds because the set on which $c_\mu(y)=1$ has measure $N_\mu$ (since each coset has measure $1$), and the proportion of this set on which $c_\lambda^{\alt}(y^T\pi_\nu^{\alt}y)=1$ equals $\mathbf{P}_{\mu,\nu}^{\alt,\lambda}$. On the other hand, set $z=y^Txy$. Since the measure over $G^{\alt}$ is $G$-invariant, we have
\begin{align}\label{eq: alt_z and y integral}
\begin{split}
\mathcal{I}&=\int_{G^{\alt}\times G}c_\lambda^{\alt}(z)c_\mu(y)c_\nu^{\alt}(y^{-T}zy^{-1})dzdy\\
&=V(K^T\pi_\lambda^{\alt}K)\int_G c_\mu(y)c_\nu^{\alt}(y^{-T}\pi_\lambda^{\alt}y^{-1})dzdy=G_{\mu,\nu}^{\alt,\lambda}(\mathfrak{o})V(K^T\pi_\lambda^{\alt}K)
\end{split}
\end{align}
where the last equality is by definition of $G_{\mu,\nu}^{\alt,\lambda}(\mathfrak{o})$. The two results from \eqref{eq: alt_x and y integral} and \eqref{eq: alt_z and y integral} together give the proof.
\end{proof}

\subsection{Symmetric function interpretation}

One way of studying the Hecke ring $H(G,K)$ is the interpretation of symmetric functions. The map between the Hecke ring $H(G,K)$ and the ring of symmetric Laurent polynomials $\Z[x_1^{\pm 1},\ldots,x_n^{\pm 1}]^{S_n}$ is given by Satake isomorphism $f\mapsto\hat\omega(f)$ based on spherical functions $\omega=\omega_s\in H(G,K)$, see \cite[Chapter V.3]{Macdonald}.

Likewise, to study the $H(G,K)$-module $H(G^{\alt},K)$, we would like to connect to symmetric functions. The following result is essentially contained in \cite{Hironaka} (we will explain how to extract it from there momentarily).


\begin{thm}\label{thm: alt_Hironaka}
Define $\Z$-linear mappings
\begin{align}\begin{split}\label{eq:image of usual_alt}\psi_n:H(G,K)&\rightarrow\C[x_1^{\pm 1},\ldots,x_n^{\pm 1}]^{S_n} \\ 
 c_\mu &\mapsto q^{\langle \mu,\rho_{2n}\rangle}P_\mu(x_1q^{\frac{1}{2}},x_1q^{-\frac{1}{2}},\ldots,x_nq^{\frac{1}{2}},x_nq^{-\frac{1}{2}};q^{-1})  
\end{split}\end{align}
and
\begin{align}\begin{split}
 \label{eq:image of alt}\psi^{\alt}_n:H(G^{\alt},K)&\rightarrow\C[x_1^{\pm 1},\ldots,x_n^{\pm 1}]^{S_n} \\ 
 c_\nu^{\alt}&\mapsto q^{2\langle \nu,\rho_n\rangle}P_{\nu}(x_1,\ldots,x_n;q^{-2})  
\end{split}\end{align}
where $\langle \cdot,\cdot\rangle$ is the canonical inner product, $\rho_{2n}=\frac{1}{2}(2n-1,2n-3,\ldots,1-2n),\rho_n=\frac{1}{2}(n-1,n-3,\ldots,1-n)$. Then we have
\begin{equation}\label{eq: alt_Satake isomorphism}
\psi_n(f)\psi^{\alt}_n(g)=\psi^{\alt}_n(f*g),\quad\forall f\in H(G,K),g\in H(G^{\alt},K).
\end{equation}
In other words, the mappings $(\psi_n,\psi^{\alt}_n)$ give a module homomorphism from the $H(G,K)$-module $H(G^{\alt},K)$ to $\C[x_1^{\pm 1},\ldots,x_n^{\pm 1}]^{S_n}$, viewed as a module over itself.
\end{thm}

\begin{proof}
Given a tuple of complex numbers $z=(z_1,\ldots,z_n)\in\C^n$, \cite{Hironaka} defines maps 
\begin{equation}\label{eq: alt_usual fourier}
f\mapsto \tilde f(z)=(f*\Psi_z)/\Psi_z: H(G,K) \to \C
\end{equation}
\begin{equation}\label{eq: alt_fourier}g\mapsto\hat g(z)=\int_{G^{\alt}}g(x)\Psi_z(x^{-1})dx: H(G^{\alt},K) \to \C
\end{equation}
in \cite[Lemma 2.1]{Hironaka} and \cite[(2.5)]{Hironaka} respectively by certain integrals (The symbol $\Psi_z$ is known as the spherical function). They check in \cite[Lemma 2.1]{Hironaka} that
\begin{equation}\label{eq: alt_comm_diag_cplx}
    (f*g)^\wedge(z)=\tilde f(z)\hat g(z).
\end{equation}
Then, in \cite[(1.13)(1.14)]{Hironaka}, they compute that 
$$\tilde{c_\mu}(z)=q^{\langle\mu,\rho_{2n}\rangle}P_\mu(q^{-z_1+1/2},q^{-z_1-1/2},\ldots,q^{-z_n+1/2},q^{-z_n-1/2};q^{-1})$$
and in \cite[Theorem 3, Lemma 2.7] {Hironaka} they compute that
$$\hat c_\nu^{\alt}(z)=q^{2\langle \nu,\rho_n\rangle}P_{\nu}(q^{-z_1},\ldots,q^{-z_n};q^{-2}) $$
which (extended by linearity) serves as an alternate definition of $\tilde f(z)$ and $\hat g(z)$. Because \eqref{eq: alt_comm_diag_cplx} holds for any $z_1,\ldots,z_n \in \C$ with $\tilde f(z)$ and $\hat g(z)$ defined by \eqref{eq: alt_usual fourier} and \eqref{eq: alt_fourier} respectively, it follows that \eqref{eq: alt_Satake isomorphism} holds as long as we apply our new maps $\psi_n(f),\psi^{\alt}_n(g)$ (instead of $\tilde f,\hat g$ in the original work) and replace the complex number $q^{-z_i}$ by the variable $x_i$ for all $1\le i\le n$.
\end{proof}


\begin{cor}\label{cor: alt_polynomial in q}
$G_{\mu,\nu}^{\alt,\lambda}(\mathfrak{o})$ is a ``polynomial in $q$", i.e., there exists a polynomial $g_{\mu,\nu}^{\alt,\lambda}(t)\in\Z[t]$, independent of $\mathfrak{o}$, such that $G_{\mu,\nu}^{\alt,\lambda}(\mathfrak{o})=g_{\mu,\nu}^{\alt,\lambda}(q)$. 
\end{cor}
\begin{proof}
\eqref{eq: alt_Littlewood Richardson}, \eqref{eq:structure coefficient_alt}, \eqref{eq:image of usual_alt}, \eqref{eq:image of alt}, and \eqref{eq: alt_Satake isomorphism} together connect the coefficient $G_{\mu,\nu}^{\alt,\lambda}(\mathfrak{o})$ which corresponds to the alternating Hecke module, and $c_{\mu,\nu}^{\alt,\lambda}(q^{-1})$ which corresponds to the symmetric Laurent polynomials:
\begin{equation}\label{eq: alt_g and c}
G_{\mu,\nu}^{\alt,\lambda}(\mathfrak{o})=q^{2\langle \nu-\lambda,\rho_n\rangle+\langle\mu,\rho_{2n}\rangle+\frac{1}{2}|\mu|}c_{\mu,\nu}^{\alt,\lambda}(q^{-1})=q^{2n(\lambda)-2n(\nu)-n(\mu)+|\mu|}c_{\mu,\nu}^{\alt,\lambda}(q^{-1})
\end{equation}
The second equality holds because we always have $|\mu|+|\nu|=|\lambda|$ when the coefficient $c_{\mu,\nu}^{\alt,\lambda}(q^{-1})$ is nonzero. Hence, the ``function over $q$" statement is automatically true. 

Recall by \Cref{rmk:lr_poly} that $c_{\mu,\nu}^{\alt,\lambda}(q^{-1})$ is a polynomial in $q^{-1}$. Also, by part \eqref{item: alt_nonnegative integer} of \Cref{prop: Hall alternating polynomial}, $G_{\mu,\nu}^{\alt,\lambda}(\mathfrak{o})\in\Z$ for every $q$ as a power of prime. Therefore, the degree of the polynomial $c_{\mu,\nu}^{\alt,\lambda}(q^{-1})\in \Z[q^{-1}]$ must be less than $2n(\lambda)-2n(\nu)-n(\mu)+|\mu|$, and therefore the explicit form in \eqref{eq: alt_g and c} is a polynomial in $q$ with integer coefficients.
\end{proof}

From now on, we always write $g_{\mu,\nu}^{\alt,\lambda}(q)$ for the structure coefficient instead of $G_{\mu,\nu}^{\alt,\lambda}(\mathfrak{o})$. With the above preparation, we can start our proof of \Cref{thm:Product process}.

\begin{proof}[Proof of \Cref{thm:Product process}, alternating case] 

By multiplying a power of $\pi$ over $B$, there is no loss of generality to assume $\mu\in\Sig_{2n}^+$ has non-negative parts. By (2.9) of  \cite[Chapter V.2]{Macdonald}, we have
\begin{equation}\label{eq: alt_formula of N_mu}
N_\mu=q^{2\langle\mu,\rho_{2n}\rangle}\frac{V_{2n}(q^{-1})}{V_\mu(q^{-1})}.
\end{equation}
Also, \cite[Lemma 2.7]{Hironaka} provides the following:
\begin{equation}\label{eq: alt_volume of lambda and nu}
V(K^T\pi_\lambda^{\alt}K)=q^{4\langle\lambda,\rho_n\rangle}\frac{V_n(q^{-2})}{V_\lambda(q^{-2})},\quad V(K^T\pi_\nu^{\alt}K)=q^{4\langle\nu,\rho_n\rangle}\frac{V_n(q^{-2})}{V_\nu(q^{-2})}
\end{equation}

Apply \Cref{prop: alt_transition step} and the results from \eqref{eq: alt_g and c}, \eqref{eq: alt_formula of N_mu}, and \eqref{eq: alt_volume of lambda and nu}, we give the proof of the alternating case of \Cref{thm:Product process}:
\begin{align}\begin{split}
\mathbf{P}_{\mu,\nu}^{\alt,\lambda}&=\frac{g_{\mu,\nu}^{\alt,\lambda}(q)V(K^T\pi_\lambda^{\alt}K)}{N_\mu V(K^T\pi_\nu^{\alt}K)}\\
&=q^{2\langle \lambda-\nu,\rho_n\rangle-\langle\mu,\rho_{2n}\rangle+\frac{1}{2}|\mu|}c_{\mu,\nu}^{\alt,\lambda}(q^{-1})\frac{V_\mu(q^{-1})}{V_{2n}(q^{-1})}\frac{V_\nu(q^{-2})}{V_\lambda(q^{-2})}\\
&=q^{2n(\nu)-2n(\lambda)+n(\mu)}c_{\mu,\nu}^{\alt,\lambda}(q^{-1})\frac{V_\mu(q^{-1})}{V_{2n}(q^{-1})}\frac{V_\nu(q^{-2})}{V_\lambda(q^{-2})}\\
&=\frac{c_{\mu,\nu}^{\alt,\lambda}(q^{-1})P_\lambda(1,q^{-2},\ldots,q^{2-2n};q^{-2})}{P_\mu(1,q^{-1},\ldots,q^{1-2n};q^{-1})P_\nu(1,q^{-2},\ldots,q^{2-2n};q^{-2})}.
\end{split}\end{align}
The third row holds because we always have $|\mu|+|\nu|=|\lambda|$ when the probability is nonzero, and the last row comes from the principal specialization formula in \eqref{eq: Principal specialization formulas}.
\end{proof}

\begin{example} Let us illustrate \Cref{thm:Product process} through some small examples which can also be computed directly:
\begin{enumerate}[left=0pt]
\item Suppose $n=1$. In this case, the conditional probability is clearly $1$ when $|\lambda|=|\mu|+|\nu|$, and $0$ otherwise. This coincides with our symmetric function result in \Cref{thm:Product process}, since when we write $P_\mu(x_1,x_1t,\ldots,x_n,x_nt;t)P_\nu(x_1,\ldots,x_n;t^2)$ in the form of linear sum of $P_\lambda(x_1,\ldots,x_n;t^2)$, only the case $\lambda_1=|\mu|+|\nu|$ appears.

\quad

\item Suppose $\mu=N[2n]$. In this case, the conditional probability is one when $\lambda=2N[n]+\nu$, and zero otherwise. This coincides with our symmetric function result in \Cref{thm:Product process}, since 

$$P_{N[2n]}(x_1,x_1t,\ldots,x_n,x_nt;t)=(x_1\cdots x_n)^{2N}t^N$$

Therefore, when we write $P_\mu(x_1,x_1t,\ldots,x_n,x_nt;t)P_\nu(x_1,\ldots,x_n;t^2)$ in the form of linear sum of $P_\lambda(x_1,\ldots,x_n;t^2)$, only the case $\lambda=2N[n]+\nu$ appears.

\quad

\item Suppose $\mu=(1,0,\ldots,0),\nu=(\nu_1,\ldots,\nu_n),\lambda=(\nu_1+1,\nu_2\ldots,\nu_n)$, where $\nu_1>\nu_2$. First, let us compute the probability $\mathbf{P}_{\mu,\nu}^{\alt,\lambda}$ using the linear algebraic method. This is to ask the  probability 

$$\mathbf{P}(\SN^{\alt}(B^T\pi^{\alt}_\nu B)=\lambda)$$

where $B$ is the same as \Cref{thm:Product process}. We write the representatives of left cosets of $K\pi_\mu K$. Let $\pi$ be a generator of $\mathfrak{p}$, and $a_0=0,\ldots,a_{q-1}\in\mathfrak{o}$ such that $\{a_0+\mathfrak{p},\ldots,a_{q-1}+\mathfrak{p}\}=k=\mathfrak{o}/\mathfrak{p}$. Then, the representatives are matrices such that:
\begin{enumerate}[left=0pt]
\item They are lower triangular;

\item Among all the $2n$ diagonal entries, one of them is $\pi$, and the rest is $1$;

\item For the diagonal entry that is $\pi$, the elements  at the left of it ranges over $\{a_0,\ldots,a_{p-1}\}$;

\item For the rest of the diagonal entries that are $1$, the elements at the left of them are all zero.
\end{enumerate}
One may verify that there are indeed $1+q+q^2+\cdots+q^{2n-1}=\frac{q^{2n}-1}{q-1}$ matries of the above property, which is equal to $N_\mu=q^{2\langle \mu,\rho_{2n}\rangle}\frac{V_{2n}(q^{-1})}{V_\mu(q^{-1})}$. Among these matrices, only the following $q+1$ matrices satisfies our requirement:
\begin{enumerate}[left=0pt]
\item $\diag(1,\pi,1,\ldots,1,1)$;

\item Among the diagonal, the element $\pi$ is placed at the top.
\end{enumerate}
Therefore the conditional probability is $\frac{(q+1)(q-1)}{q^{2n}-1}=\frac{q^2-1}{q^{2n}-1}$.

\quad

On the other hand, we have 
$$P_{(1,0,\ldots,0)}(x_1,x_1t,\ldots,x_n,x_nt;t)=(1+t)(x_1+\cdots+x_n)$$
Also, if we express the product $(x_1+\cdots+x_n)P_\nu(x_1,\ldots,x_n;t^2)$ as a linear combination of Hall-Littlewood Laurent polynomials, $\lambda$ will be the largest signature under the natural ordering. Hence $c_{\mu,\nu}^{\alt,\lambda}(t)=1+t$, and
\begin{align}\begin{split}
\mathbf{P}_{\mu,\nu}^{\alt,\lambda}&=\frac{c_{\mu,\nu}^{\alt,\lambda}(q^{-1})P_\lambda(1,q^{-2},\ldots,q^{2-2n};q^{-2})}{P_\mu(1,q^{-1},\ldots,q^{1-2n};q^{-1})P_\nu(1,q^{-2},\ldots,q^{2-2n};q^{-2})}\\
&=\frac{(1+q^{-1})q^{-2n(\lambda)}V_\nu(q^{-2})}{(1+q^{-1}+\cdots+q^{1-2n})q^{-2n(\nu)}V_\lambda(q^{-2})}\\
&=\frac{1+q^{-1}}{(1+q^{-1}+\cdots+q^{1-2n})q^{2(n-1)}}=\frac{q^2-1}{q^{2n}-1}
\end{split}\end{align}
which coincides with the result of the linear algebraic method.
\end{enumerate}
\end{example}

\section{The Hermitian Hecke Module}\label{sec:her_hecke}

\subsection{The Hermitian Hecke module $H(G^{\her},K)$}

Let $G=\GL_n(F)$ be the group of all invertible $n\times n$ matrices over $F$. Also, let $G^+=G\cap M_n(\mathfrak{o})$ be the subsemigroup of $G$ consisting of all matrices $x\in G$ with entries $x_{ij}\in\mathfrak{o}$, and let $K=\GL_n(\mathfrak{o})=G^+\cap(G^+)^{-1}$ so that $K$ consists of all $x\in G$ with entries $x_{ij}\in\mathfrak{o}$ and $\det(x)$ a unit in $\mathfrak{o}$. The structure of the rings $L(G,K),L(G^+,K),H(G,K),H(G^+,K)$ follows from \Cref{defi: Hecke algebra}.

\begin{defi}
Let $L(G^{\her},K)$ (resp. $L(G^{+{\her}},K)$) denote the space of all complex-valued continuous functions of compact support on $G^{\her}=G\cap{\Her}_n(F)$ (resp. $G^{+{\her}}=G^+\cap{\Her}_n(\mathfrak{o})$) which are invariant with respect to $K$, i.e., such that

$$g(k^*xk)=g(x)$$
for all $x\in G$ (resp. $x\in G^+$) and $k\in K$. We may shall regard $L(G^{+{\her}},K)$ as a submodule of $L(G^{\her},K)$. 

\end{defi}

We define a multiplication of $L(G,K)$ over the module $L(G^{\her},K)$ as follows: for $f\in L(G,K)$, $g\in L(G^{\her},K)$,

$$(f*g)(x)=\int_Gf(y)g(y^{-1}xy^{-*})dy$$

Since $f$ and $g$ are compactly supported, the integration is over a compact set. This product is associative with respect to $L(G,K)$, i.e. $f_1*(f_2*g)=(f_1*f_2)*g$. (We use the same notation $*$ for the transition involution of a given matrix, the convolution of the Hecke algebra itself,  and the convolution of the Hecke algebra over the symplectic Hecke module, but it is easy to distinguish these three based on context.)

Each function $g\in L(G^{\her},K)$ is constant on the orbits $\{k^*xk|k\in K\}$ in $G^{\her}$. These orbits are compact, open, and mutually disjoint. Since $g$ has compact support, it follows that $g$ takes non-zero values on only finitely many orbits $\{k^*xk|k\in K\}$, and hence can be written as a finite linear combination of their characteristic functions. Therefore, the characteristic functions of these orbits in $G^{\her}$ form a $\C$-basis of $L(G^{\her},K)$. 

If we vary the definition of the module $L(G^{\her},K)$ (resp. $L(G^{+{\her}},K)$) by requiring the functions to take their values in $\Z$ instead of $\C$, the resulting module is the generalization of the Hecke ring over $G^{\her}$ (resp. $G^{+{\her}}$), and we denote it by $H(G^{\her},K)$ (resp. $H(G^{+{\her}},K)$). Clearly we have

$$L(G^{\her},K)=H(G^{\her},K)\otimes_\Z\C,L(G^{+{\her}},K)=H(G^{+{\her}},K)\otimes_\Z\C$$

Consider an orbit $\{k^*xk|k\in K\}$, where $x\in G^{\her}$. By multiplying $x$ by a suitable power of $\pi$ (the generator of $\mathfrak{p}$) we can bring $x$ into $G^{+{\her}}$. Also, \Cref{prop:smith} implies that each orbit $\{k^*xk|k\in K\}$ has a unique representative of the form

$$\pi_\lambda=\diag(\pi^{\lambda_1},\ldots,\pi^{\lambda_n})$$
where $\lambda_1\ge\lambda_2\ge\ldots\ge\lambda_n$. We have $\lambda_n\ge 0$ if and only if $x\in G^{+{\her}}$.

Let $c^{\her}_\lambda$ denote the characteristic function of the orbit $\{k^*\pi_\lambda k|k\in K\}$. Then we have the $c^{\her}_\lambda$ (resp. the $c^{\her}_\lambda$ such that $\lambda_n\ge 0$) form a $\Z$-basis of $H(G^{\her},K)$(resp. $H(G^{+{\her}},K)$). The characteristic function $c_0$ of $K$ is the identity element of $H(G,K)$ and $H(G^+,K)$, and this also plays the role of the identity element when multiplying with elements in $H(G^{\her},K)$ and $H(G^{+{\her}},K)$. 

\begin{defi}
Let $\mu,\nu\in\Sig_n$. Then, we define the structure coefficient $G_{\mu,\nu}^{\her,\lambda}(\mathfrak{o})$ for the expansion of the product $c_\mu *c^{\her}_\nu$:
\begin{equation}\label{eq:structure coefficient_her}c_\mu *c^{\her}_\nu=\sum_\lambda G_{\mu,\nu}^{\her,\lambda}(\mathfrak{o})c^{\her}_\lambda
\end{equation}
\end{defi}

\begin{prop}\label{prop: Hall Hermitian polynomial}
$G_{\mu,\nu}^{\her,\lambda}(\mathfrak{o})$ has the following properties:
\begin{enumerate}[left=0pt]
\item $G_{\mu,\nu}^{\her,\lambda}(\mathfrak{o})=0$ unless $|\lambda|=2|\mu|+|\nu|$. In this case, $G_{\mu,\nu}^{\her,\lambda}(\mathfrak{o})\in\Z_{\ge 0}$ is a non-negative integer; \label{item: her_nonnegative integer}

\item $G_{\mu,\nu}^{\her,\lambda}(\mathfrak{o})=G_{\mu,N[n]+\nu}^{\her,N[n]+\lambda}(\mathfrak{o})=G_{N[n]+\mu,\nu}^{\her,2N[n]+\lambda}(\mathfrak{o})$ for all $N\in\Z$;
\end{enumerate}
\end{prop}

\begin{proof}
\begin{enumerate}[left=0pt]
\item We now give an interpretation of this coefficient $G_{\mu,\nu}^{\her,\lambda}(\mathfrak{o})$. Notice that
$$G_{\mu,\nu}^{\her,\lambda}(\mathfrak{o})=(c_\mu*c_\nu^{\her})(\pi_\lambda)=\int_G c_\mu(y)c_\nu^{\her}(y^{-1}\pi_\lambda y^{-*})dy$$
Since $c_\mu(y)$ vanishes for $y$ outside $K\pi_\mu K$, the integration is over this orbit, which we shall write as a disjoint union of left cosets, say
\begin{equation}\label{eqref: her_disjoint left coset}K\pi_\mu K=\bigsqcup_j y_jK\quad (y_j\in K\pi_\mu)\end{equation}
Therefore, we have
$$G_{\mu,\nu}^{\her,\lambda}(\mathfrak{o})=\sum_j\int_{y_jK}c^{\her}_\nu(y^{-1}\pi_\lambda y^{-*})dy=\sum_jc^{\her}_\nu(y_j^{-1}\pi_\lambda y_j^{-*})$$
since $K$ has measure $1$. Hence $G_{\mu,\nu}^{\her,\lambda}(\mathfrak{o})$ is equal to the number of $j$ such that $y_j^{-1}\pi_\lambda y_j^{-*}\in\{k\pi_\nu k^*|k\in K\}$, which is a non-negative integer; If $|\lambda|\ne 2|\mu|+|\nu|$, by checking the absolute value of the determinant we know such $y_j$ does not exist, and therefore $G_{\mu,\nu}^{\her,\lambda}(\mathfrak{o})=0$;

\item By multiplying $\pi^N$ to the alternating matrix in the middle, we obtain $G_{\mu,\nu}^{\her,\lambda}(\mathfrak{o})=G_{\mu,N[n]+\nu}^{\her,N[n]+\lambda}(\mathfrak{o})$; And by multiplying $\pi^N$ to the two matrices on both sides, we obtain $G_{\mu,\nu}^{\her,\lambda}(\mathfrak{o})=G_{N[n]+\mu,\nu}^{\her,2N[n]+\lambda}(\mathfrak{o})$;
\end{enumerate}
\end{proof}

\begin{rmk} 
As in the alternating case, the order of the subscript cannot be changed, i.e., $G_{\mu,\nu}^{\her,\lambda}(\mathfrak{o})$ is not the same as $G_{\nu,\mu}^{\her,\lambda}(\mathfrak{o})$. 
\end{rmk}

The setting above helps us study the distribution of singular numbers of matrix products. 

\begin{prop}\label{prop: her_transition step}
For all $\lambda\in\Sig_n$, let $$V(K^*\pi_\lambda^{\her}K)=\int_{G^{\her}}c_\lambda^{\her}(x)dx$$ 
denote the volume of the orbit $\{k^*\pi_\lambda k|k\in K\}$, where $dx$ is the $G$-invariant measure on $G^{\her}$ normalized by $\int_{K^*K}dx=1$. Then the probability $\mathbf{P}_{\mu,\nu}^{\her,\lambda}:=\mathbf{P}(\SN^{\her}(B^*AB)=\lambda)$ appeared in \Cref{thm:Product process} has the form
$$\mathbf{P}_{\mu,\nu}^{\her,\lambda}=\frac{G_{\mu,\nu}^{\her,\lambda}(\mathfrak{o})V(K^*\pi_\lambda K)}{N_\mu V(K^*\pi_\nu K)}$$
where $N_\mu$ is the number of $y_j$ in \eqref{eqref: her_disjoint left coset}, i.e. the number of disjoint left cosets in $K\pi_\mu K$.
\end{prop}

\begin{proof}
Consider the integral 
$$\mathcal{I}=\int_{G^{\her}\times G}c_\nu^{\her}(x)c_\mu(y)c_\lambda^{\her}(y^*xy)dxdy.$$
On one hand, we have
\begin{equation}\label{eq: her_x and y integral}\mathcal{I}=V(K^*\pi_\nu K)\int_G c_\mu(y)c_\lambda(y^*\pi_\nu y)dy=\mathbf{P}_{\mu,\nu}^{\her,\lambda}N_\mu V(K^*\pi_\nu K).\end{equation}
The second equality holds because the set on which $c_\mu(y)=1$ has measure $N_\mu$ (since each coset has measure $1$), and the proportion of this set on which $c_\lambda^{\her}(y^T\pi_\nu^{\her}y)=1$ equals $\mathbf{P}_{\mu,\nu}^{\her,\lambda}$. On the other hand, set $z=y^*xy$. Since the measure over $G^{\her}$ is $G$-invariant, we have
\begin{align}\label{eq: her_z and y integral}
\begin{split}
\mathcal{I}&=\int_{G^{\her}\times G}c_\lambda^{\her}(z)c_\mu(y)c_\nu^{\her}(y^{-*}zy^{-1})dzdy\\
&=V(K^*\pi_\lambda K)\int_G c_\mu(y)c_\nu^{\her}(y^{-*}\pi_\lambda y^{-1})dzdy=G_{\mu,\nu}^{\her,\lambda}(\mathfrak{o})V(K^*\pi_\lambda K)
\end{split}
\end{align}
where the last equality is by definition of $G_{\mu,\nu}^{\her,\lambda}(\mathfrak{o})$. The two results from \eqref{eq: her_x and y integral} and \eqref{eq: her_z and y integral} together give the proof.
\end{proof}

\subsection{Symmetric function interpretation}

As an analogy of the alternating case, to study the $H(G,K)$-module $H(G^{\her},K)$, we would like to connect to symmetric functions. The following result is essentially contained in \cite{Hironaka_Hermitian} and \cite{Hironaka_Hermitian_and_Symmetric_I} (we will explain how to extract it from there momentarily).

\begin{thm}\label{thm: her_Hironaka}
Define $\Z$-linear mappings
\begin{align}\begin{split}\label{eq:image of usual_her}\psi_n:H(G,K)&\rightarrow\C[x_1^{\pm 1},\ldots,x_n^{\pm 1}]^{S_n}\\
c_\mu&\mapsto q^{2\langle \mu,\rho_n\rangle}P_\mu(x_1^2,\ldots,x_n^2;q^{-2})\end{split}\end{align}
and
\begin{align}\begin{split}\label{eq:image of her}
\psi^{\her}_n:H(G^{\her},K)&\rightarrow\C[x_1^{\pm 1},\ldots,x_n^{\pm 1}]^{S_n}\\
c_\nu^{\her}&\mapsto (-1)^{n(\nu)+|\nu|}q^{\langle \nu,\rho_n\rangle}P_\nu(x_1,\ldots,x_n;-q^{-1})
\end{split}\end{align}
where $\langle \cdot,\cdot\rangle$ is the canonical inner product, $\rho_n=\frac{1}{2}(n-1,n-3,\ldots,1-n)$. Then we have
\begin{equation}\label{eq: her_Satake isomorphism}
\psi_n(f)\psi^{\her}_n(g)=\psi_n^{\her}(f*g),\quad\forall f\in H(G,K),g\in H(G^{\her},K).
\end{equation}
Thus the mappings $(\psi_n,\psi^{\her}_n)$ give a module homomorphism from the $H(G,K)$-module $H(G^{\her},K)$ to $\C[x_1^{\pm 1},\ldots,x_n^{\pm 1}]^{S_n}$, viewed as a module over itself.
\end{thm}

\begin{proof}
Given a tuple of complex numbers $z=(z_1,\ldots,z_n)\in\C^n$, \cite{Hironaka_Hermitian_and_Symmetric_I} defines maps 
\begin{equation}\label{eq: her_usual fourier}
f\mapsto \tilde f(z)=(f*\zeta(\cdot;z))/\zeta(\cdot;z)): H(G,K) \to \C
\end{equation}
\begin{equation}\label{eq: her_fourier}g\mapsto\hat g(z)=\int_{G^{\her}}g(x)\zeta(x^{-1};z)dx: H(G^{\her},K) \to \C
\end{equation}
in \cite[Lemma 1.9]{Hironaka_Hermitian_and_Symmetric_I} and \cite[(1.10)]{Hironaka_Hermitian_and_Symmetric_I} respectively by certain integrals (The symbol $\zeta(\cdot;z))$ is known as the spherical function). She checks in \cite[Lemma 1.11]{Hironaka_Hermitian_and_Symmetric_I} that
\begin{equation}\label{eq: her_comm_diag_cplx}
    (f*g)^\wedge(z)=\tilde f(z)\hat g(z).
\end{equation}
Then, in \cite[page 209]{Hironaka_Hermitian_and_Symmetric_I}, she computes that\footnote{In the original paper \cite{Hironaka_Hermitian_and_Symmetric_I}, Hironaka actually wrote $\tilde{c_\mu}(z)=q^{2\langle\mu,\rho_{n}\rangle}P_\mu(q^{-4z_1},\ldots,q^{-4z_n};q^{-2})$. However, the definition of $z$ in \cite[(2.2)]{Hironaka_Hermitian} and \cite[(1.6)]{Hironaka_Hermitian_and_Symmetric_I} are twice as different. Therefore, we revise the explicit form so that the definition of $z$ in our proof becomes uniform.}
$$\tilde{c_\mu}(z)=q^{2\langle\mu,\rho_{n}\rangle}P_\mu(q^{-2z_1},\ldots,q^{-2z_n};q^{-2})$$
and in \cite[(2.6),(2.7)] {Hironaka_Hermitian} she computes that
$$\hat c_\nu^{\her}(z)=(-1)^{n(\nu)+|\nu|}q^{\langle \nu,\rho_n\rangle}P_{\nu}(q^{-z_1},\ldots,q^{-z_n};-q^{-1}) $$
which (extended by linearity) serves as an alternate definition of $\tilde f(z)$ and $\hat g(z)$. Because \eqref{eq: her_comm_diag_cplx} holds for any $z_1,\ldots,z_n$ with $\tilde f(z)$ and $\hat g(z)$ defined by \eqref{eq: her_usual fourier} and \eqref{eq: her_fourier} respectively, it follows that \eqref{eq: her_Satake isomorphism} holds as long as we apply our new maps $\psi_n(f),\psi^{\her}_n(g)$ (instead of $\tilde f,\hat g$ in the original work) and replace the complex number $q^{-z_i}$ by the variable $x_i$ for all $1\le i\le n$.
\end{proof}

\begin{cor}\label{cor: her_polynomial in q}
$G_{\mu,\nu}^{\her,\lambda}(\mathfrak{o})$ is a ``polynomial in $q$", i.e., there exists a polynomial $g_{\mu,\nu}^{\her,\lambda}(t)\in\Z[t]$, independent of $\mathfrak{o}$, such that $G_{\mu,\nu}^{\her,\lambda}(\mathfrak{o})=g_{\mu,\nu}^{\her,\lambda}(q)$. 
\end{cor}
\begin{proof}
\eqref{eq: her_Littlewood Richardson}, \eqref{eq:structure coefficient_her}, \eqref{eq:image of usual_her}, \eqref{eq:image of her}, and \eqref{eq: her_Satake isomorphism} together connect the coefficient $g_{\mu,\nu}^{\her,\lambda}(q)$ which corresponds to the Hermitian Hecke module, and $c_{\mu,\nu}^{\her,\lambda}(q^{-1})$ which corresponds to the symmetric Laurent polynomials: (we eliminate $(-1)^{|\lambda|-|\nu|}$ since when the coefficient $c_{\mu,\nu}^{\her,\lambda}(t)$ is nontrivial, $|\lambda|$ and $|\nu|$ must have the same parity): 
\begin{equation}\label{eq: her_g and c}
G_{\mu,\nu}^{\her,\lambda}(\mathfrak{o})=(-1)^{n(\lambda
)-n(\nu)}q^{\langle\nu-\lambda+2\mu,\rho_n\rangle}c_{\mu,\nu}^{\her,\lambda}(q^{-1})=(-1)^{n(\lambda
)-n(\nu)}q^{n(\lambda)-n(\nu)-2(\mu)}c_{\mu,\nu}^{\her,\lambda}(q^{-1})
\end{equation} 
The second equality holds because we always have $2|\mu|+|\nu|=|\lambda|$ when the coefficient $c_{\mu,\nu}^{\her,\lambda}(q^{-1})$ is nonzero. Hence, the ``function over $q$" statement is automatically true. Also, by part \eqref{item: her_nonnegative integer} of \Cref{prop: Hall Hermitian polynomial}, $G_{\mu,\nu}^{\her,\lambda}(\mathfrak{o})\in\Z$ for every $q$ as a power of prime. Therefore, the degree of the polynomial $c_{\mu,\nu}^{\her,\lambda}(q^{-1})\in \Z[q^{-1}]$ must be less than $n(\lambda)-n(\nu)-2(\mu)$, and therefore the explicit form in \eqref{eq: her_g and c} is a polynomial in $q$ with integer coefficients.
\end{proof}

From now on, we always write $g_{\mu,\nu}^{\her,\lambda}(q)$ for the structure coefficient instead of $G_{\mu,\nu}^{\her,\lambda}(\mathfrak{o})$. With the above preparation, we can start our proof of \Cref{thm:Product process}.

\begin{proof}[Proof of \Cref{thm:Product process}, Hermitian case] 
By multiplying a power of $\pi$ over $B$, there is no loss of generality to assume $\mu\in\Sig_n^+$ has non-negative parts. By (2.9) of \cite[Chapter V.2]{Macdonald}, we have
\begin{equation}\label{eq: her_formula of N_mu}
N_\mu=q^{4\langle\mu,\rho_n\rangle}\frac{V_n(q^{-2})}{V_\mu(q^{-2})}
\end{equation}
Also, \cite[(2.7)]{Hironaka_Hermitian} provides the following:
\begin{equation}\label{eq: her_volume of lambda and nu}
V(K^*\pi_\lambda K)=q^{2\langle\lambda,\rho_n\rangle}\frac{V_n(-q^{-1})}{V_\lambda(-q^{-1})},\quad V(K^*\pi_\nu K)=q^{2\langle\nu,\rho_n\rangle}\frac{V_n(-q^{-1})}{V_\nu(-q^{-1})}
\end{equation}
Apply \Cref{prop: her_transition step} and the results from \eqref{eq: her_g and c}, \eqref{eq: her_formula of N_mu}, and \eqref{eq: her_volume of lambda and nu}, we give the proof of the Hermitian case of \Cref{thm:Product process}:
\begin{align}\begin{split}
\mathbf{P}_{\mu,\nu}^{\her,\lambda}&=\frac{g_{\mu,\nu}^{\her,\lambda}(q)V(K^*\pi_\lambda K)}{N_\mu V(K^*\pi_\nu K)}\\
&=(-1)^{n(\lambda
)-n(\nu)}q^{\langle\lambda-\nu-2\mu,\rho_n\rangle}c_{\mu,\nu}^{\her,\lambda}(q^{-1})\frac{V_\mu(q^{-2})}{V_{n}(q^{-2})}\frac{V_\nu(-q^{-1})}{V_\lambda(-q^{-1})}\\
&=(-1)^{n(\lambda
)-n(\nu)}q^{n(\nu)-n(\lambda)+2n(\mu)}c_{\mu,\nu}^{\her,\lambda}(q^{-1})\frac{V_\mu(q^{-2})}{V_{n}(q^{-2})}\frac{V_\nu(-q^{-1})}{V_\lambda(-q^{-1})}\\
&=\frac{c_{\mu,\nu}^{\her,\lambda}(q^{-1})P_\lambda(1,-q^{-1},\ldots,(-q)^{1-n};-q^{-1})}{P_\mu(1,q^{-2},\ldots,q^{2-2n};q^{-2})P_\nu(1,-q^{-1},\ldots,(-q)^{1-n};-q^{-1})}
\end{split}\end{align}
The third row holds because we always have $2|\mu|+|\nu|=|\lambda|$ when the probability is nonzero, and the last row comes from the principal specialization formula in \eqref{eq: Principal specialization formulas}.
\end{proof}

\begin{example} 
Let us illustrate \Cref{thm:Product process} through some small examples which can be computed by other means:
\begin{enumerate}[left=0pt]
\item $\mu=N[n]$. In this case, the conditional probability is clearly $1$ when $\lambda=2N[n]+\nu$, and $0$ otherwise. This coincides with our symmetric function result in \Cref{thm:Product process}, since 

$$P_{N[n]}(x_1^2,x_2^2,\ldots,x_n^2;t^2)=(x_1\cdots x_n)^{2M}$$
Therefore, when we write $P_\mu(x_1^2,x_2^2,\ldots,x_n^2;t^2)P_\nu(x_1,\ldots,x_n;-t)$ in the form of linear sum of $P_\lambda(x_1,\ldots,x_n;-t)$, only the case $\lambda=2N[n]+\nu$ appears.

\item $\mu=(1,0,\ldots,0),\nu=(\nu_1,\ldots,\nu_n),\lambda=(\nu_1+2,\ldots,\nu_n)$, where $\nu_1>\nu_2$. First, let us compute the probability $\mathbf{P}_{\mu,\nu}^{\her,\lambda}$ by linear algebraic method. This is to ask the probability 

$$\mathbf{P}(\SN^{\her}(B^*\pi_\nu B)=\lambda)$$
where $B$ is the same as \Cref{thm:Product process}.
We write the representatives of left cosets of $K\pi_\mu K$. Let $\pi$ be a generator of $\mathfrak{p}$, and $a_0=0,\ldots,a_{q-1}\in\mathfrak{o}$ such that $\{a_0+\mathfrak{p},\ldots,a_{q-1}+\mathfrak{p}\}=k=\mathfrak{o}/\mathfrak{p}$. Then, the representatives are matrices such that:
\begin{enumerate}[left=0pt]
\item They are lower triangular;

\item Among all the $n$ diagonal entries, one of them is $\pi$, and the rest is $1$;

\item For the diagonal entry that is $\pi$, the elements  at the left of it ranges over $\{a_0,\ldots,a_{p-1}\}$;

\item For the rest of the diagonal entries that are $1$, the elements at the left of them are all zero.
\end{enumerate}
One may verify that there are indeed $1+q^2+q^4+\cdots+q^{2n-2}=\frac{q^{2n}-1}{q^2-1}$ matries of the above property, which is equal to $N_\mu=q^{4\langle \mu,\rho\rangle}\frac{V_n(q^{-2})}{V_\mu(q^{-2})}$. Among these matrices, only $\diag(\pi,1,\ldots,1)$ satisfies our requirement. Hence, the conditional probability is $\frac{q^2-1}{q^{2n}-1}$. On the other hand, we have

$$P_{(1,0,\ldots,0)}(x_1^2,\ldots,x_n^2;t^2)=x_1^2+\cdots+x_n^2$$
Also, if we express the product $(x_1^2+\cdots+x_n^2)P_\nu(x_1,\ldots,x_n;-t)$ into the sum of Hall-Littlewood Laurent polynomials, $\lambda$ will be the largest signature under the natural ordering. Hence $c_{\mu,\nu}^{\her,\lambda}(t)=1$, and

\begin{align}\begin{split}
\mathbf{P}_{\mu,\nu}^{\her,\lambda}&=\frac{c_{\mu,\nu}^{\her,\lambda}(q^{-1})P_\lambda(1,-q^{-1},\ldots,(-q)^{1-n};q^{-1})}{P_\mu(1,q^{-2},\ldots,q^{2-2n};q^{-2})P_\nu(1,-q^{-1},\ldots,(-q)^{1-n});q^{-1})}\\
&=\frac{q^{-n(\lambda)}V_\nu(-q^{-1})}{(1+q^{-2}+\cdots+q^{2-2n})q^{-n(\nu)}V_\lambda(-q^{-1})}\\
&=\frac{1}{(1+q^{-2}+\cdots+q^{2-2n})q^{2(n-1)}}=\frac{q^2-1}{q^{2n}-1}
\end{split}\end{align}
which coincides with the linear algebraic method.
\end{enumerate}
\end{example}

\section{The Hall-Littlewood Measure, Product Process, and Corners Process}\label{sec:HL}

In this section, we combine \Cref{thm:Product process} with Hall-Littlewood combinatorics to prove \Cref{thm:Corner process} and other results in the introduction section. The joint distributions of corners will be given, which enables us to study the corner of a Haar distributed matrix with singular numbers all zero. The proofs we give are somewhat different for the alternating and Hermitian cases, but the only real reason for this is that for the alternating case one can take a shortcut and reduce to results of \cite{van2020limits}.

To prepare for the proof of \Cref{thm:Corner process} we give several technical lemmas, which we prove later in this section. The first, for the alternating case, relates the product operation of \Cref{thm:Product process} to products of matrices with no symmetry restrictions, which were understood in \cite{van2020limits}.

\begin{lemma}\label{thm:alt_to_un}
Let $\lambda \in \Sig_n, \nu \in \Sig_{n-1}$ be fixed, let $r\in\mathbb{Z}_{\ge 0}$ be sufficiently large, and let $\kappa=(r+|\lambda|-|\nu|,\nu_1,\ldots,\nu_{n-1})$. Let $A,B \in \Mat_n(\tilde{F})$ be two matrices with i.i.d. additive Haar entries chosen from a degree $2$ unramified extension $\tilde{F}$ of $F$. Then
\begin{equation}
\mathbf{P}_{(r)\lambda}^{\alt,\kappa} = \mathbf{P}(\SN(AB)=\kappa\mid \SN(A)=\lambda, \SN(B)=(r)),
\end{equation}
where $\mathbf{P}_{(r)\lambda}^{\alt,\kappa}$ is as defined in \Cref{thm:Product process}.
\end{lemma}

For the Hermitian case we argue slightly differently, giving an explicit form of the Pieri case of the Littlewood-Richardson-type coefficients defined earlier.

\begin{lemma}\label{thm:explicit_her_lr_coefs}
    For any $\nu \in \Sig_n, \mu \in \Sig_{n-1}$ and $D$ large enough so that all $n$-tuples below are valid signatures, 
    \begin{equation}
        c_{(D,0[n-1]),\nu}^{\her,(2D+|\nu|-|\mu|,\mu_1,\ldots,\mu_{n-1})}(t) = \sum_{\kappa \in \Sig_{n-1}} P_{\nu/\kappa}(1;-t)Q_{\mu/\kappa}(-1;-t),
    \end{equation}
    where $c_{\mu,\nu}^{\her,\la}(t)$ is as defined in \eqref{eq: her_Littlewood Richardson}. 
\end{lemma}

\begin{proof}[Proof of \Cref{thm:Corner process}]
We begin with the alternating case, \eqref{item:corner_alternating} of \Cref{thm:Corner process}, and further begin with the odd case \eqref{eq:odd_alt_corner}. The probability that we wish to compute is equal to

\begin{align}\begin{split}
\mathbf{P}_{2n<2n+1}^{\alt}(\lambda\mid\nu)&=\mathbf{P}(\SN^{\alt}(\diag(1,\ldots,1,0)A\diag(1,\ldots,1,0))=\lambda)\\
&=\mathbf{P}(\SN^{\alt}(\diag(1,\ldots,1,0)U^T\pi_{\nu}^{\alt}U\diag(1,\ldots,1,0))=\lambda)\\
&=\mathbf{P}(\SN^{\alt}(B^T\pi_{\nu}^{\alt}B)=\lambda)
\end{split}\end{align}
where $U\in\GL_{2n+1}(\mathfrak{o})$ is Haar distributed, $B$ is the $2n\times2n$ left upper corner of $U$, and all diagonal matrices are ${(2n+1)\times(2n+1)}$. Also, the symbol $\pi_{\nu}^{\alt}$ refers to the $(2n+1)\times(2n+1)$ matrix in the second row, and refers to the $2n\times 2n$ matrix in the third row. Then by \cite[Theorem 1.3]{van2020limits}, $\SN(B)$ is distributed with respect to the Hall-Littlewood measure:
$$\mathbf{P}(\SN(B)=(m))=\frac{Q_{(m)}(1;t)P_{(m)}(t,t^2,\ldots,t^{2n};t)}{\Pi_t(1;t,t^2,\ldots,t^{2n})}=\frac{Q_{(m)}(t;t)P_{(m)}(1,t,\ldots,t^{2n-1};t)}{\Pi_t(t;1,t,\ldots,t^{2n-1})}$$
Then we apply the result from the product convolution \eqref{eq: alt_conditional distribution of HL process} in the Preliminaries section: set $\theta=(1,t^2,\ldots,t^{2n-2}),\theta'=(1,t,\ldots,t^{2n-1}),\psi_1=(t)$. Then, we immediately get 

$$\mathbf{P}_{2n<2n+1}^{\alt}(\lambda\mid\nu)=\frac{Q_{\lambda/\nu}(\psi_1;t^2)P_{\lambda}(\theta;t^2)}{P_\nu(\theta;t^2)\Pi_{t^2}(\theta;\psi_1)}=\frac{Q_{\lambda/\nu}(t;t^2)P_\lambda(1,t^2,\ldots,t^{2n-2};t^2)}{P_\nu(1,t^2,\ldots,t^{2n-2};t^2)\Pi_{t^2}(t;1,t^2,\ldots,t^{2n-2})},$$
proving \eqref{eq:odd_alt_corner}.

For \eqref{eq:even_alt_corner}, as an analogy of the proof for the usual corner process in \cite[Theorem 1.3]{van2020limits}, we regard $\mathbf{P}_{2n-1<2n}^{\alt}(\nu\mid\lambda)$ as the limit of $\mathbf{P}_{(r)\lambda}^{\alt,\kappa}$ when $r$ goes to infinity. By \Cref{thm:alt_to_un} and usual corner process for matrices over $\tilde F$, we have
\begin{align}\begin{split}
\mathbf{P}_{2n-1<2n}^{\alt}(\nu\mid\lambda)&=\lim_{r\rightarrow\infty}\mathbf{P}_{(r)\lambda}^{\alt,\kappa}=\lim_{r\rightarrow\infty}\mathbf{P}_{(r)\lambda}^{\text{un},\kappa}\\
&=\mathbf{P}(\SN(\tilde{A}_{n-1})=\nu)=\frac{P_{\lambda/\nu}(1;t^2)P_\nu(t^2,t^4,\ldots,t^{2n-2};t^2)}{P_\lambda(1,t^2,\ldots,t^{2n-2};t^2)}\end{split}\end{align}
where $\tilde{A}_{n-1}$ is the $(n-1)\times n$ corner of $\tilde A\in\Mat_{n}(\tilde F)$, which is the random matrix with fixed singular numbers $\lambda$ and distribution invariant under $\GL_n(\tilde{\mathfrak{o}})$ on both sides. 

We now proceed to the Hermitian case, part \eqref{item:hermitian_corner} of \Cref{thm:Corner process}. Note that for any $r > \max(\nu_1,\la_1)$, by reducing modulo $\pi^r$ we have
\begin{align}\label{eq:fix_r_herm}
    \begin{split}
            \mathbf{P}(\SN(A_{n-1}) = \nu) &= \mathbf{P}(\SN(\diag(1[n-1],\pi^r) A \diag(1[n-1],\pi^r)) = (2r+|\la|-|\nu|,\nu)) \\ 
            &= \bP_{(r,0[n-1]),\la}^{\her,(2r+|\la|-|\nu|,\nu)},
    \end{split}
\end{align}
where the right hand side is as defined in \Cref{thm:Product process} and $A_{n-1}$ is as in our theorem statement. Hence it suffices to compute the right hand side of \eqref{eq:fix_r_herm} for sufficiently large $r$. Combining \Cref{thm:Product process} and the explicit formula \Cref{thm:explicit_her_lr_coefs}, we have
\begin{multline}
    \bP_{(r,0[n-1]),\la}^{\her,(2r+|\la|-|\nu|,\nu)} = \frac{P_{(2r+|\la|-|\nu|,\nu)}(1,-q^{-1},\ldots,(-q^{-1})^{n-1};-q^{-1})}{P_{(r,0[n-1])}(1,q^{-2},\ldots,(q^{-2})^{n-1};q^{-2}) P_\la(1,-q^{-1},\ldots,(-q^{-1})^{n-1};-q^{-1})} \\ 
    \times \sum_{\kappa \in \Sig_{n-1}} P_{\la/\kappa}(1;-q^{-1})Q_{\nu/\kappa}(-1;-q^{-1}).
\end{multline}
Moving terms, the above is
\begin{equation}\label{eq:wts=1}
    \frac{P_{(2r+|\la|-|\nu|,\nu)}(1,-q^{-1},\ldots,(-q^{-1})^{n-1};-q^{-1})\Pi_{-q^{-1}}(-1;-q^{-1},\ldots,(-q^{-1})^{n-1})}{P_\nu(-q^{-1},\ldots,(-q^{-1})^{n-1};-q^{-1})P_{(r,0[n-1])}(1,q^{-2},\ldots,(q^{-2})^{n-1};q^{-2})} \cdot \text{RHS\eqref{eq:herm_corner}},
\end{equation}
so we must show the prefactor in \eqref{eq:wts=1} is $1$. This is a straightforward computation from the definition of the Cauchy kernel \eqref{eq: Cauchy kernel} and principal specialization formulas \eqref{eq: Principal specialization formulas}. 
\end{proof}

We now prove the lemmas used above.

\begin{proof}[Proof of \Cref{thm:alt_to_un}]
A simple computation with the branching rule \eqref{eq:Pbranch_with_GT} shows that for $m>l\ge 0,P_{(m)/(l)}(x,xt;t)=P_{(m)/(l)}(x;t^2)(1+t)^{\textbf{1}_{l=0}}$. Thus by the branching rule of the Hall-Littlewood symmetric function, we have
\begin{align}\begin{split}
P_{(r)}(x_1,x_1t,\ldots,x_n,x_nt;t)&=\sum_{0\le m_1\le\ldots\le m_n=r}P_{(m_1)}(x_1,x_1t;t)\cdots P_{(m_n)/(m_{n-1})}(x_n,x_nt;t)\\
&=(1+t)\sum_{0\le m_1\le\ldots\le m_n=r}P_{(m_1)}(x_1;t^2)\cdots P_{(m_n)/(m_{n-1})}(x_n;t^2)\\
&=(1+t)P_{(r)}(x_1,\ldots,x_n;t^2)
\end{split}\end{align}
Hence we have $c_{(r)\lambda}^{\alt,\kappa}(t)=(1+t)c_{(r)\lambda}^\kappa(t^2)$, and
\begin{align}\begin{split}
\mathbf{P}_{(r)\lambda}^{\alt,\kappa}&=\frac{c_{(r)\lambda}^{\alt,\kappa}(q^{-1})P_\lambda(1,q^{-2},\ldots,q^{2-2n};q^{-2})}{P_\mu(1,q^{-1},\ldots,q^{1-2n};q^{-1})P_\nu(1,q^{-2},\ldots,q^{2-2n};q^{-2})}\\
&=\frac{c_{(r)\lambda}^\kappa(q^{-2})P_\kappa(1,q^{-2},\ldots,q^{2-2n};q^{-2})}{P_{(r)}(1,q^{-2},\ldots,q^{2-2n};q^{-2})P_\lambda(1,q^{-2},\ldots,q^{2-2n};q^{-2})}\\
&=\mathbf{P}_{(r)\lambda}^{\text{un},\kappa},
\end{split}\end{align}
which completes the proof.
\end{proof}

The following basic combinatorial lemma is necessary to prove \Cref{thm:explicit_her_lr_coefs}.
\begin{lemma}
\label{thm:factor_hl}
Let $\tau_1,\ldots,\tau_{\binom{n}{k}} \in S_n$ be a set of coset representatives for $S_n/(S_k \times S_{n-k})$, let $\mu \in \Sig_k, \nu \in \Sig_{n-k}$, and define
\begin{equation}
\la(D) = (\mu_1+D,\ldots,\mu_k+D,\nu_1,\ldots,\nu_n)
\end{equation}
for all $D$ large enough such that this is a valid signature. Then 
\begin{multline}\label{eq:factor_hl}
P_{\la(D)}(x_1,\ldots,x_n;t) = \\ 
\sum_{j=1}^{\binom{n}{k}} \tau_j\left( \Pi_t(x_1^{-1},\ldots,x_k^{-1};x_{k+1},\ldots,x_n) x_1^D \cdots x_k^D P_\mu(x_1,\ldots,x_k;t) P_\nu(x_{k+1},\ldots,x_n;t)\right).
\end{multline}
\end{lemma}
\begin{proof}
Direct computation:
\begin{align}
\begin{split}
&P_{\la(D)}(x_1,\ldots,x_n;t) = \frac{(1-t)^n}{\Pi_i (t;t)_{m_i(\la(D))}}\sum_{\sigma \in S_n} \sigma\left((x_1 \cdots x_k)^D \prod_{i=1}^k x_i^{\mu_i} \prod_{i=k+1}^n x_i^{\nu_{i-k}} \prod_{1 \leq i < j \leq n} \frac{x_i - t x_j}{x_i - x_j}\right) \\ 
&= \sum_{j=1}^{\binom{n}{k}} \tau_j\left( \Pi_t(x_1^{-1},\ldots,x_k^{-1};x_{k+1},\ldots,x_n) x_1^D \cdots x_k^D \frac{(1-t)^k}{\Pi_i (t;t)_{m_i(\mu)}}   \right. \\ 
&\left.\times \sum_{\sigma \in S_k}\sigma\left(\prod_{i=1}^k x_i^{\mu_i}  \prod_{1 \leq i < j \leq n} \frac{x_i - t x_j}{x_i - x_j} \right)\frac{(1-t)^{n-k}}{\Pi_i (t;t)_{m_i(\nu)}} \sum_{\sigma' \in S_{n-k}} \sigma'\left( \prod_{i=k+1}^n x_i^{\nu_{i-k}} \prod_{k+1 \leq i < j \leq n} \frac{x_i - t x_j}{x_i - x_j} \right)\right) \\ 
&=\text{RHS\eqref{eq:factor_hl}}.
\end{split}
\end{align}
\end{proof}

\begin{rmk}
    \Cref{thm:factor_hl} is inspired by, and is morally a special case of, a result \cite[Theorem 1.5]{van2020limits} which was proven in the generality of Macdonald polynomials. We had to reprove it here because we needed an exact equality rather than the slightly weaker limiting statement which was given in that result, making it not literally a special case.
\end{rmk}

\begin{proof}[Proof of \Cref{thm:explicit_her_lr_coefs}]
We have
    \begin{multline}\label{eq:lr_herm_pieri_spec_case}
        P_{(D,0[n-1])}(x_1^2,\ldots,x_1^2;t^2) P_\nu(x_1,\ldots,x_n;-t) = \\
        \sum_{\la} c_{(D,0[n-1]),\nu}^{\her,(2D+|\nu|-|\la|,\la)}(t) P_{(2D+|\nu|-|\la|,\la)}(x_1,\ldots,x_n;-t)
    \end{multline}
    by \eqref{eq: her_Littlewood Richardson} and homogeneity of Hall-Littlewood polynomials, where the sum is over all $\la  \in \Sig_{n-1}$ for which $(2D+|\nu|-|\la|,\la) \in \Sig_n$. Since Hall-Littlewood polynomials form a basis, it therefore suffices to show
    \begin{multline}\label{eq:lr_suffices}
        P_{(D,0[n-1])}(x_1^2,\ldots,x_1^2;t^2) P_\nu(x_1,\ldots,x_n;-t) \\  = \sum_{\la} \left( \sum_{\kappa \in \Sig_{n-1}} P_{\nu/\kappa}(1;-t)Q_{\la/\kappa}(-1;-t)\right) P_{(2D+|\nu|-|\la|,\la)}(x_1,\ldots,x_n;-t)
    \end{multline} 
    Applying \Cref{thm:factor_hl} to the $D$-dependent polynomial on both sides and letting $\tau_1,\ldots,\tau_n$ be coset representatives of $S_n/(S_1 \times S_{n-1})$, \eqref{eq:lr_suffices} is equivalent to 
    \begin{multline}\label{eq:second_lr_suffices}
            \sum_{j=1}^n \tau_j(\Pi_{t^2}(x_1^{-2};x_2^{2},\ldots,x_n^{2}))P_\nu(x_1,\ldots,x_n;-t) = \\  \sum_{j=1}^n \tau_j\left( \sum_{\la,\kappa \in \Sig_{n-1}} P_{\nu/\kappa}(1;-t)Q_{\la/\kappa}(-1;-t) \Pi_{-t}(x_1^{-1};x_2,\ldots,x_n) x_1^{|\nu|-|\la|}P_\la(x_2,\ldots,x_n;-t)\right).
    \end{multline}
    Since $P_\nu(x_1,\ldots,x_n;-t)$ is symmetric, we may bring it inside $\tau_j$, so to show \eqref{eq:second_lr_suffices} it suffices to show the term-by-term equality 
    \begin{multline}
                \label{eq:cauchy_swap}
        \Pi_{t^2}(x_1^{-2};x_2^{2},\ldots,x_n^{2})P_\nu(x_1,\ldots,x_n;-t) = \\  \sum_{\la \in \Sig_{n-1}} \left(\sum_{\kappa \in \Sig_{n-1}} P_{\nu/\kappa}(1;-t)Q_{\la/\kappa}(-1;-t)\right) \Pi_{-t}(x_1^{-1};x_2,\ldots,x_n) x_1^{|\nu|-|\la|}P_\la(x_2,\ldots,x_n;-t). 
    \end{multline}
    Writing out the Cauchy kernel explicitly and completing the square we see 
    \begin{equation}
        \Pi_{t^2}(x_1^{-2};x_2^{2},\ldots,x_n^{2}) = \Pi_{-t}(x_1^{-1};x_2,\ldots,x_n)\Pi_{-t}(-x_1^{-1};x_2,\ldots,x_n).
    \end{equation}
    By the branching rule,
    \begin{equation}
        P_\nu(x_1,\ldots,x_n;-t) = \sum_{\kappa \in \Sig_{n-1}} P_{\nu/\kappa}(x_1;-t)P_\kappa(x_2,\ldots,x_n;-t).
    \end{equation}
    Hence by the skew Cauchy identity,
    \begin{align}
        \begin{split}
            \text{LHS\eqref{eq:cauchy_swap}} &= \Pi_{-t}(x_1^{-1};x_2,\ldots,x_n) \sum_\kappa P_{\nu/\kappa}(x_1;-t) \left( \Pi_{-t}(-x_1^{-1};x_2,\ldots,x_n)P_\kappa(x_2,\ldots,x_n;-t)\right) \\ 
            &= \Pi_{-t}(x_1^{-1};x_2,\ldots,x_n) \sum_\kappa P_{\nu/\kappa}(x_1;-t) \sum_{\la \in \Sig_{n-1}} Q_{\la/\kappa}(-x_1^{-1};-t) P_\la(x_2,\ldots,x_n;-t) \\ 
            &= \text{RHS\eqref{eq:cauchy_swap}},
        \end{split}
    \end{align}
    where $\kappa$ is sum over $\Sig_{n-1}$ as above. Therefore we complete the proof.\end{proof}

From \Cref{thm:Corner process} it is easy to derive explicit expressions as Hall-Littlewood measures for the distributions of singular numbers of matrices with i.i.d. Haar entries and corners of invertible matrices. The former was already stated as \Cref{thm:The Hall-Littlewood measure of i.i.d entries}, and the latter we give now. It is a direct analogy to \cite[Theorem 1.3(1)]{van2020limits}, which showed that the singular numbers of a corner of a Haar-distributed element of $\GL_N(\Z_p)$ have Hall-Littlewood distribution. 

\begin{cor}\label{thm:Corner of the invertible matrix}
Let $t=1/q$, and $m>n\ge 1$ be integers.
\begin{enumerate}[left=0pt]
\item(Alternating case)\begin{enumerate}[left=0pt]
\item Let $A$ be the top left $2n\times2n$ submatrix of a Haar-distributed element of $\Alt_{2m}(\mathfrak{o})\cap \GL_{2m}(\mathfrak{o})$. Then $\SN^{\alt}(A)\in\Sig_n^+$ has distribution given by the Hall-Littlewood measure $\mathbf{P}_{2n<2m}^{\alt}(\cdot\mid0)$ defined by

\begin{equation}
\mathbf{P}_{2n<2m}^{\alt}(\lambda\mid0)=\frac{P_\lambda(1,t^2,\ldots,t^{2n-2};t^2)Q_\lambda(t,t^3,\ldots,t^{2m-2n-1};t^2)}{\Pi_{t^2}(1,t^2,\ldots,t^{2n-2};t,t^3,\ldots,t^{2m-2n-1})}.
\end{equation}

\item Let $A$ be the top left $(2n+1)\times(2n+1)$ submatrix of a Haar-distributed element of $\Alt_{2m}(\mathfrak{o})\cap \GL_{2m}(\mathfrak{o})$. Then $\SN^{\alt}(A)\in\Sig_n^+$ has distribution given by the Hall-Littlewood measure $\mathbf{P}_{2n+1< 2m}^{\alt}(\cdot\mid0)$ defined by

\begin{equation}
\mathbf{P}_{2n+1< 2m}^{\alt}(\lambda\mid0)=\frac{P_\lambda(1,t^2,\ldots,t^{2n-2};t^2)Q_\lambda(t^3,t^5,\ldots,t^{2m-2n-1};t^2)}{\Pi_{t^2}(1,t^2,\ldots,t^{2n-2};t^3,t^5,\ldots,t^{2m-2n-1})}.
\end{equation}
\end{enumerate}

\item(Hermitian case)Let $A$ be the top left $n\times n$ submatrix of a Haar-distributed element of $\Her_m(\mathfrak{o})\cap \GL_m(\mathfrak{o})$.  Then $\SN^{\her}(A)
\in\Sig_n^+$ has distribution given by the Hall-Littlewood measure $\mathbf{P}_{n<m}^{\her}(\cdot\mid0)$ defined by

\begin{equation}
\mathbf{P}_{n<m}^{\her}(\lambda\mid0)=\frac{P_\lambda(1,-t,t^2,\ldots,(-t)^{n-1};-t)Q_\lambda(t,-t^2,\ldots,-(-t)^{m-n};-t)}{\Pi_{-t}(1,-t,t^2,\ldots,(-t)^{n-1};t,-t^2,\ldots,-(-t)^{m-n})}.
\end{equation}
\end{enumerate}
\end{cor}

\begin{proof}
We begin with part (a), the even alternating case. Set

$$\lambda_{2n}^{2m-2}=(\lambda^{(2n)},\lambda^{(2n+2)},\ldots,\lambda^{(2m-2)}),\nu_{2n+1}^{2m-3}=(\nu^{(2n+1)},\nu^{(2n+3)},\ldots,\nu^{(2m-3)})$$
where $\lambda^{(2n)},\lambda^{(2n+2)},\ldots,\lambda^{(2m-2)}$ denote the singular numbers of the $(2n)\times(2n),(2n+2)\times (2n+2),\ldots,(2m-2)\times(2m-2)$ left upper corner, and $\nu^{(2n+1)},\nu^{(2n+3)},\ldots,\nu^{(2m-3)}$ denote the singular numbers of the $(2n+1)\times (2n+1),\ldots,(2m-3)\times (2m-3)$ left upper corner. Then by \Cref{thm:Corner process}, the $(2m-1)\times(2m-1)$ left upper corner is invertible, and the joint variable $(\lambda_{2n}^{2m-2},\nu_{2n+1}^{2m-3})$ is distributed with respect to the Hall-Littlewood process, which has the form

$$\mathbf{P}^{\alt}(\lambda_{2n}^{2m-2},\nu_{2n+1}^{2m-3})=\mathbf{P}_{2m-2<2m-1}^{\alt}(\lambda^{(2m-2)}\mid 0)\cdots\mathbf{P}_{2n<2n+1}^{\alt}(\lambda^{(2n)}\mid\nu^{(2n+1)})$$
which by \Cref{thm:Corner process}, is furthermore equal to

$$Q_{\lambda^{(2m-2)}}(t;t^2)P_{\lambda^{(2m-2)}/\nu^{(2m-3)}}(1;t^2)Q_{\lambda^{(2m-4)}/\nu^{(2m-3)}}(t^{-1};t^2)\cdots P_{\lambda^{(2n+2)}/\nu^{(2n+1)}}(t^{2m-2n-4};t^2)\cdot$$
$$Q_{\lambda^{(2n)}/\nu^{(2n+1)}}(t^{2n-2m+3};t^2)P_{\lambda^{(2n)}}(t^{2m-2n-2},t^{2m-2n},\ldots,t^{2m-4};t^2)/\Pi
$$
where $\Pi$ has the form 

$$\Pi=\Pi_{t^2}(1;t)\Pi_{t^2}(1;t,t^3)\cdots\Pi_{t^2}(1;t,\ldots,t^{2m-2n-3})\Pi_{t^2}(1,\ldots,t^{2n-2};t,t^3,\ldots,t^{2m-2n-1})$$

Therefore, we have $\mathbf{P}_{2n<2m}^{\alt}(\lambda\mid0)=\sum_{\lambda^{(2n)}=\lambda}\mathbf{P}^{\alt}(\lambda_{2n}^{2m},\nu_{2n+1}^{2m-1})$, which is equal by repeated applications of the Cauchy identity to

\begin{align}
&\sum_{\lambda^{(2m-2)},\ldots,\nu^{(2n+1)}}Q_{\lambda^{(2m-2)}}(t;t^2)\cdots Q_{\lambda/\nu^{(2n+1)}}(t^{2n-2m+3};t^2)P_{\lambda}(t^{2m-2n-2},\ldots,t^{2m-4};t^2)/\Pi\nonumber\\
&=\sum_{\nu^{(2m-3)},\ldots,\nu^{(2n+1)}}\Pi_{t^2}(1;t)Q_{\nu^{(2m-3)}}(t;t^2)\cdots P_{\lambda}(t^{2m-2n-2},\ldots,t^{2m-4};t^2)/\Pi\nonumber\\
&=\ldots\nonumber\\
&=\Pi_{t^2}(1;t)\cdots\Pi_{t^2}(1;t,\ldots,t^{2m-2n-3})Q_\lambda(t^{2n-2m+3},\ldots,t;t^2)P_\lambda(t^{2m-2n-2},\ldots,t^{2m-4};t^2)/\Pi\nonumber\\
&=P_\lambda(1,t^2,\ldots,t^{2n-2};t^2)Q_\lambda(t,t^3,\ldots,t^{2m-2n-1};t^2)/\Pi_{t^2}(1,t^2,\ldots,t^{2n-2};t,t^3,\ldots,t^{2m-2n-1})\nonumber
\end{align}
which ends the proof.

We now consider the odd alternating case. Motivated by the above, we have $\mathbf{P}_{2n+1<2m}^{\alt}(\lambda\mid0)=\sum_{\nu^{(2n+1)}=\lambda}\mathbf{P}^{\alt}(\lambda_{2n}^{2m},\nu_{2n+1}^{2m-1})$, which is equal (again by repeated applications of the Cauchy identity) to

\begin{align}
&\sum_{\lambda^{(2m-2)},\ldots,\lambda^{(2n+2)},\lambda^{(2n)}}Q_{\lambda^{(2m-2)}}(t;t^2)\cdots Q_{\lambda^{(2n)}/\lambda}(t^{2n-2m+3};t^2)P_{\lambda^{(2n)}}(t^{2m-2n-2},\ldots,t^{2m-4};t^2)/\Pi\nonumber\\
&=\ldots\nonumber\\
&=\Pi_{t^2}(1;t)\cdots\Pi_{t^2}(1;t,\ldots,t^{2m-2n-3})\times\nonumber\\
&\sum_{\lambda^{(2n)}}Q_\lambda(t^{2n-2m+5},\ldots,t;t^2)Q_{\lambda^{(2n)}/\lambda}(t^{2n-2m+3};t^2)P_{\lambda^{(2n)}}(t^{2m-2n-2},\ldots,t^{2m-4};t^2)/\Pi\nonumber\\
&=P_\lambda(1,t^2,\ldots,t^{2n-2};t^2)Q_\lambda(t^3,t^5,\ldots,t^{2m-2n-1};t^2)/\Pi_{t^2}(1,t^2,\ldots,t^{2n-2};t^3,t^5,\ldots,t^{2m-2n-1})\nonumber
\end{align}
which ends the proof.

We now show the Hermitian case. Set
$$\lambda_n^{m-1}=(\lambda^{(n)},\lambda^{(n+1)},\ldots,\lambda^{(m-1)})$$
where $\lambda^{(n)},\lambda^{(n+1)},\ldots,\lambda^{(m-1)}$ denote the singular numbers of the $n\times n,(n+1)\times (n+1),\ldots,(m-1)\times(m-1)$ left upper corner. Then by \Cref{thm:Corner process}, the variable $\lambda_n^{m-1}$ is distributed with respect to the process with the form

$$\mathbf{P}^{\her}(\lambda_n^{m-1})=\mathbf{P}_{m-1<m}^{\her}(\lambda^{(m-1)}\mid 0)\cdots\mathbf{P}_{n<n+1}^{\her}(\lambda^{(n)}\mid\lambda^{(n+1)})$$
which is furthermore equal to

$$\sum_{\kappa^{(m-1)},\ldots,\kappa^{(n+1)}}Q_{\lambda^{(m-1)}}(1;-t^{-1})P_{\lambda^{(m-1)}/\kappa^{(m-1)}}(t;-t)Q_{\lambda^{(m-2)}/\kappa^{(m-1)}}(-t^{-1};-t)\cdots$$

$$P_{\lambda^{(n+1)}/\kappa^{(n+1)}}(-(-t)^{m-n-1};-t)Q_{\lambda^{(n)}/\kappa^{(n+1)}}((-t)^{n-m+1};-t)P_{\lambda^{(n)}}(-(-t)^{m-n},\ldots,-(-t)^{m-1})/\Pi$$
where $\Pi$ has the form
$$\Pi=\Pi_{-t}(1;-t)\cdots\Pi_{-t}(1;-t,\ldots-(-t)^{m-n-1})\Pi_{-t}(1,\ldots,(-t)^{n-1};t,\ldots,-(-t)^{m-n})$$
Therefore, we have $\mathbf{P}_{n<m}^{\her}(\lambda\mid0)=\sum_{\lambda^{(n)}=\lambda}\mathbf{P}^{\her}(\lambda_n^{m-1})$, which is equal (again by repeated applications of the Cauchy identity) to

\begin{align*}
&\sum_{\lambda^{(m-1)},\ldots,\kappa^{(n+1)}}Q_{\lambda^{(m-1)}}(1;-t)\cdots Q_{\lambda/\kappa^{(n+1)}}((-t)^{n-m-1};-t)P_{\lambda}(-(-t)^{m-n},\ldots,-(-t)^{m-1};-t)/\Pi\nonumber\\
&=\sum_{\kappa^{(m-1)},\ldots,\kappa^{(n+1)}}\Pi_{-t}(1;-t)Q_{\kappa^{(m-1)}}(1;-t)\cdots P_\lambda(-(-t)^{m-n},\ldots,-(-t)^{m-1};-t)/\Pi\nonumber\\
&=\ldots\nonumber\\
&=\Pi_{-t}(1;-t)\cdots\Pi_{-t}(1;-t,\ldots,(-t)^{m-n-1})\times\\
&Q_\lambda((-t)^{n-m+1},\ldots,1;-t)P_\lambda(-(-t)^{m-n},\ldots,-(-t)^{m-1};-t)/\Pi\nonumber\\
&=P_\lambda(1,\ldots,(-t)^{n-1};-t)Q_\lambda(t,\ldots,-(-t)^{m-n};-t)/\Pi_{-t}(1,\ldots,(-t)^{n-1};t,-t^2,\ldots,-(-t)^{m-n})\nonumber
\end{align*}
which ends the proof.
\end{proof}

Using \Cref{thm:Corner of the invertible matrix}, we may recover formulas for the distribution of singular numbers of alternating and Hermitian matrices with i.i.d. Haar entries, recovering results of \cite{bhargava2013modeling}, \cite{fulman2018random} and \cite{lee2022universality} as mentioned in the Introduction. In the Hermitian case, the $n \to \infty$ limit of this measure was also found in random matrix theory over finite fields \cite[Section 4.4]{fulman_thesis}.

\begin{cor}\label{thm:The Hall-Littlewood measure of i.i.d entries}

Let $t=1/q$, and $n\ge 1$ be an integer.

\begin{enumerate}[left=0pt]
\item(Alternating case)\begin{enumerate}[left=0pt]
\item Let $A\in\Alt_{2n}(\mathfrak{o})$ be random with i.i.d entries above the diagonal, distributed according to the additive Haar measure on $\mathfrak{o}$. Then, the singular numbers of $A$ are distributed with respect to the Hall-Littlewood measure $\mathbf{P}^{\alt}_{2n}$ defined by

\begin{align}
\begin{split}
\mathbf{P}^{\alt}_{2n}(\lambda)&=\frac{P_\lambda(1,t^2,\ldots,t^{2n-2};t^2)Q_\lambda(t,t^3,\ldots;t^2)}{\prod_{t^2}(1,t^2,\ldots,t^{2n-2};t,t^3,\ldots)}\\
&=t^{4n(\lambda)+|\lambda|}\frac{(1-t)(1-t^2)\cdots(1-t^{2n})}{\Pi_{i\ge 0}(t^2;t^2)_{m_i(\lambda)}}.
\end{split}
\end{align}

\item Let $A\in\Alt_{2n+1}(\mathfrak{o})$ be random with i.i.d entries above the diagonal, distributed according to the additive Haar measure on $\mathfrak{o}$. Then, the singular numbers of $A$ are distributed with respect to the Hall-Littlewood measure $\mathbf{P}^{\alt}_{2n+1}$ defined by

\begin{align}
\begin{split}
\mathbf{P}^{\alt}_{2n+1}(\lambda)&=\frac{P_\lambda(1,t^2,\ldots,t^{2n-2};t^2)Q_\lambda(t^3,t^5,\ldots;t^2)}{\Pi_{t^2}(1,t^2,\ldots,t^{2n-2};t^3,t^5,\ldots)}\\
&=t^{4n(\lambda)+3|\lambda|}\frac{(1-t^2)(1-t^3)\cdots(1-t^{2n+1})}{\prod_{i\ge 0}(t^2;t^2)_{m_i(\lambda)}}.
\end{split}
\end{align}
\end{enumerate}

\item (Hermitian case) Let $A\in\Her_n(\mathfrak{o})$ be random with i.i.d entries above the diagonal, distributed according to the additive Haar measure on $\mathfrak{o}$, and i.i.d entries over the diagonal, distributed according to the additive Haar measure on $\mathfrak{o}\cap F_0$ (the entries above and over the diagonal are also independent). Then, the singular numbers of $A$ are distributed with respect to the Hall-Littlewood measure $\mathbf{P}^{\her}_n$ defined by

\begin{align}
\begin{split}
\mathbf{P}^{\her}_n(\lambda)&=\frac{P_\lambda(1,-t,t^2,\ldots,(-t)^{n-1};-t)Q_\lambda(t,-t^2,\ldots;-t)}{\Pi_{-t}(1,-t,t^2,\ldots,(-t)^{n-1};t,-t^2,t^3,\ldots)}\\
&=t^{2n(\lambda)+|\lambda|}\frac{(1-t^2)(1-t^4)\cdots(1-t^{2n})}{\prod_{i\ge 0}(-t;-t)_{m_i(\lambda)}}.
\end{split}
\end{align}
\end{enumerate}
\end{cor}

\begin{proof}[Proof of \Cref{thm:The Hall-Littlewood measure of i.i.d entries}] We prove the Hermitian case here, and the alternating case is similar. Let $m>n$ be sufficiently large, and $A_m\in\Her_m(\mathfrak{o})$ denote the random matrix with i.i.d entries above the diagonal, distributed according to the additive Haar measure on $\mathfrak{o}$, and i.i.d entries on the diagonal, distributed according to the additive Haar measure on $\mathfrak{o}\cap F_0$. Letting $A$ become the $n\times n$ upper left corner of $A_m$, we view $\mathbf{P}_{n<m}^{\her}(\lambda\mid0)$ as the conditional probability 
$$\mathbf{P}_{n<m}^{\her}(\lambda\mid0)=\mathbf{P}(\SN^{\her}(A)=\lambda\mid \SN^{\her}(A_m)=0).$$
Let $\mathbf{P}_{m>n}^{\her}(0\mid\lambda)$ denote the conditional probability
\begin{align}\begin{split}
\mathbf{P}_{m>n}^{\her}(0\mid\lambda)&:=\mathbf{P}( \SN^{\her}(A_m)=0\mid \SN^{\her}(A)=\lambda)\\
&=\mathbf{P}(\SN^{\her}(A_m)=0\mid A=\diag(1,\ldots,1,\pi^{\lambda_1},\ldots,\pi^{\lambda_l})).
\end{split}\end{align}
Then we have
\begin{equation}\label{eq:Haar matrix}
\mathbf{P}_n^{\her}(\lambda)\mathbf{P}_{m>n}^{\her}(0\mid\lambda)=\mathbf{P}(\SN^{\her}(A)=\lambda, \SN^{\her}(A_m)=0)=\mathbf{P}_m^{\her}(0)\mathbf{P}_{n<m}^{\her}(\lambda\mid0)
\end{equation}
Now let $l=l(\lambda)$. Then we must have $l+n\le m$, otherwise there is no way the matrix $A_m$ can be invertible. Conditioned on $A=\diag(1,\ldots,1,\pi^{\lambda_1},\ldots,\pi^{\lambda_l})$, denote the random matrix
$$A_m=\begin{pmatrix}I_{n-l} & 0 & B\\ 0 & \diag(\pi^{\lambda_1},\ldots,\pi^{\lambda_l}) & C\\ 
B^* & C^* & D\end{pmatrix}$$
where $B,C\in M_{l,(m-n)}(\mathfrak{o})$ have unrestricted Haar distributed entries, and $D\in\Her_{m-n}(\mathfrak{o})$ has i.i.d. Haar distributed entries subject to the Hermitian condition. Hence we have
\begin{align}\label{eq:cbd}
\begin{split}
\mathbf{P}_{m>n}^{\her}(0\mid\lambda)&=\mathbf{P}(\SN^{\her}\begin{pmatrix}I_{n-l} & 0 & B\\ 0 & \diag(\pi^{\lambda_1},\ldots,\pi^{\lambda_l}) & C\\ 
B^* & C^* & D\end{pmatrix}=0)\\
&=\mathbf{P}(\begin{pmatrix}I_{n-l} & 0 & B'\\ 0 & 0 & C'\\ B'^* & C'^* & D'\end{pmatrix}\in\GL_m(k)\text{ is invertible})\\
&=\mathbf{P}(\begin{pmatrix}I_{n-l} & 0 & 0\\ 0 & 0 & C'\\ 0 & C'^* & D'-B'^*B'\end{pmatrix}\in\GL_m(k)\text{ is invertible})\\
&=\mathbf{P}(\begin{pmatrix}0 & C'\\ C'^* & D'-B'^*B'\end{pmatrix}\in\GL_{m-n+l}(k)\text{ is invertible})
\end{split}
\end{align}
where the prime symbol refers to reducing modulo $\mathfrak{p}$, the first equality is essentially by definition, the second is by reducing modulo $\mathfrak{p}$, and the third is by using the $(n-l)\times(n-l)$ identity matrix on the left upper side to eliminate the first $(n-l)$ rows and columns by left- and right-multiplying by the appropriate element of $\GL_m(k)$. Since the Haar measure is additive and multiplicative (under general linear group over $\mathfrak{o}$ on both sides) invariant, $C'\in M_{l,(m-n)}(k)$, and $D', D'-B'^*B'\in\Her_{m-n}(k)$ are all distributed uniformly. When the matrix 
$$\begin{pmatrix}0 & C'\\ C'^* & D'-B'^*B'\end{pmatrix}$$ 
is invertible, $C'$ has to be full rank, so we can perform a linear transformation to replace $C'$ with $\diag_{l\times(m-n)}(1,\ldots,1)$ and eliminate $l$ rows and columns of the matrix $D'-B'^*B'$ to yield a new $(m-n-l) \times (m-n-l)$ matrix $D''$. Because $D'-B'^*B'$ is uniform and is independent of $C'$, $D''$ is a uniformly random $(m-n-l) \times (m-n-l)$ Hermitian matrix, and the desired event of \eqref{eq:cbd} now holds if and only if $D''$ is invertible. Hence

\begin{align}\begin{split}
\mathbf{P}_{m>n}^{\her}(0\mid\lambda)&=\mathbf{P}(\rank(C')=l)\mathbf{P}_{m-n-l}^{\her}(0)\\
&=(1-q^{2n-2m})\cdots(1-q^{2n-2m-2l+2})\frac{|\Her_{m-n-l}(k)\cap\GL_{m-n-l}(k)|}{|\Her_{m-n-l}(k)|},
\end{split}\end{align}
where we use \Cref{thm:invertible probability} for the second equality. To sum up, by applying the explicit form of the Cauchy kernel \eqref{eq: Cauchy kernel} and principal specialization formulas \eqref{eq: Principal specialization formulas}, we have

\begin{align}\begin{split}
\mathbf{P}_n^{\her}(\lambda)&=\frac{\mathbf{P}_m^{\her}(0)\mathbf{P}_{n<m}^{\her}(\lambda\mid0)}{\mathbf{P}_{m>n}^{\her}(0\mid\lambda)}\\
&=\frac{|\Her_m(k)\cap \GL_m(k)||\Her_{m-n-l}(k)|}{|\Her_m(k)||\Her_{m-n-l}(k)\cap\GL_{m-n-l}(k)|}\frac{\mathbf{P}_{n<m}^{\her}(\lambda\mid0)}{(1-q^{2n-2m})\cdots(1-q^{2n-2m-2l+2})}\\
&=\frac{(1-(-q)^{n-m})\cdots(1-(-q)^{n-m-l+1})P_\lambda(1,\ldots,(-t)^{n-1};-t)Q_\lambda(t,\ldots,-(-t)^{m-n};-t)}{(1-q^{2n-2m})\cdots(1-q^{2n-2m-2l+2})\Pi_{-t}(1,-t,t^2,\ldots,(-t)^{n-1};t,-t^2,\ldots,-(-t)^{m-n})}\\
&=\frac{P_\lambda(1,-t,t^2,\ldots,(-t)^{n-1};-t)Q_\lambda(t,-t^2,\ldots;-t)}{\Pi_{-t}(1,-t,t^2,\ldots,(-t)^{n-1};t,-t^2,t^3,\ldots)}\\
&=t^{2n(\lambda)+|\lambda|}\frac{(1-t^2)(1-t^4)\cdots(1-t^{2n})}{\prod_{i\ge 0}(-t;-t)_{m_i(\lambda)}}.
\end{split}\end{align}
The first row is from \eqref{eq:Haar matrix}, the second is simply rephrasing probabilities in terms of the sizes of finite sets, the third row is from \Cref{thm:invertible probability} and \Cref{thm:Corner of the invertible matrix}, the fourth is from the explicit formula for the Cauchy kernel \eqref{eq: Cauchy kernel}, and the last follows from the principal specialization formulas in \Cref{prop:hl_principal_formulas}. This ends the proof.
\end{proof}

\begin{rmk}
As an analogy of the usual case, see \cite[Corollary 1.4]{van2020limits}, the distributions appearing in \Cref{thm:The Hall-Littlewood measure of i.i.d entries} are the limits of those in \Cref{thm:Corner of the invertible matrix}, i.e. 

$$\mathbf{P}_{2n}^{\alt}(\lambda)=\lim_{m\rightarrow\infty}\mathbf{P}_{2n<2m}^{\alt}(\lambda\mid0),\quad\mathbf{P}_{2n+1}^{\alt}(\lambda)=\lim_{m\rightarrow\infty}\mathbf{P}_{2n+1<2m}^{\alt}(\lambda\mid0)$$

$$\mathbf{P}_n^{\her}(\lambda)=\lim_{m\rightarrow\infty}\mathbf{P}_{n<m}^{\her}(\lambda\mid0).$$
This is intuitively true because when $m$ is very large, the  corner has very little effect on the singular numbers of the whole matrix; one may use this to give an alternate proof of \Cref{thm:The Hall-Littlewood measure of i.i.d entries} similarly to the proof of \cite[Corollary 1.4]{van2020limits}.
\end{rmk}

With the preparation above, we now turn back to the proof of \Cref{thm:corners_nomarkov_intro}. 

\begin{proof}[Proof of \Cref{thm:corners_nomarkov_intro}]
\Cref{thm:The Hall-Littlewood measure of i.i.d entries} gives the distribution of the singular numbers of the $2n\times 2n$ alternating (resp. $n\times n$ Hermitian) Haar matrix. Then, iterating the form of the corner process in \Cref{thm:Corner process} towards the $2\times 2$ alternating (resp. $1\times 1$ Hermitian) corner brings us the joint distribution. Explicitly, for the alternating case, we have
\begin{equation*}
    \text{LHS of \eqref{eq:alt_corner_proc_intro}}
    =\mathbf{P}_{2n}^{\alt}(\la^{(2n)})\mathbf{P}_{2n-1<2n}^{\alt}(\nu^{(2n-1)}\mid\la^{(2n)})\cdots\mathbf{P}_{2<3}^{\alt}(\la^{(2)}\mid\nu^{(3)})=\text{RHS of \eqref{eq:alt_corner_proc_intro}}
\end{equation*}
Likewise, for the Hermitian case, we have
\begin{equation*}
    \text{LHS of \eqref{eq:herm_corner_proc_intro}}
    =\mathbf{P}_n^{\her}(\la^{(n)})\mathbf{P}_{n-1<n}^{\her}(\la^{(n-1)}\mid\la^{(n)})\cdots\mathbf{P}_{1<2}^{\her}(\la^{(1)}\mid\la^{(2)})=\text{RHS of \eqref{eq:herm_corner_proc_intro}}
\end{equation*}
\end{proof}

As a separate corollary of \Cref{thm:Product process}, we obtain formulas for the joint distribution of 
$$\SN^{\alt}(A),\SN^{\alt}(B_1^TAB_1),\ldots,\SN^{\alt}(B_k^T\cdots B_1^TAB_1\cdots B_k)$$
for the alternating case, where $A\in\Alt_{2n}(\mathfrak{o}),B_1,\ldots,B_k\in M_{2n\times 2n}(\mathfrak{o})$ is random with i.i.d. additive Haar entries subject to the alternating restriction. We also obtain the joint distribution of 
$$\SN^{\her}(A),\SN^{\her}(B_1^*AB_1),\ldots,\SN^{\her}(B_k^*\cdots B_1^*AB_1\cdots B_k)$$
for the Hermitian case, where $A\in\Her_n(\mathfrak{o}),B_1,\ldots,B_k\in M_{n\times n}(\mathfrak{o})$ is random with i.i.d. additive Haar entries subject to the Hermitian restriction. These processes are the correct analogues of the product process considered in \cite{van2020limits}, which was shown to yield a Hall-Littlewood process in \cite[Corollary 3.4]{van2020limits}. Similarly, we find that the alternating and Hermitian settings also yield Hall-Littlewood processes:

\begin{cor}\label{thm:joint distribution}
Let $t=1/q$, and $n\ge 1$ be an integer.
\begin{enumerate}[left=0pt]
\item (Alternating case) Let $A\in\Alt_{2n}(\mathfrak{o})$ be the same as in \Cref{thm:The Hall-Littlewood measure of i.i.d entries}, and $B_1,\ldots,B_k\in M_{2n\times 2n}(\mathfrak{o})$ be i.i.d random matrices with i.i.d entries distributed by the additive Haar measure on $\mathfrak{o}$. Then, the joint distribution  
$$\mathbf{P}(\SN^{\alt}(A)=\lambda,\SN^{\alt}(B_1^tAB_1)=\lambda^{(1)},\ldots,\SN^{\alt}(B_k^t\cdots B_1^tAB_1\cdots B_k)=\lambda^{(k)})$$
is given by

\begin{equation}
\frac{Q_{\lambda}(t,t^3,\ldots;t^2)Q_{\lambda^{(1)}/\lambda}(t,t^2,\ldots;t^2)\cdots Q_{\lambda^{(k)}/\lambda^{(k-1)}}(t,t^2,\ldots;t^2)P_{\lambda^{(k)}}(1,t^2,\ldots,t^{2n-2};t^2)}{\Pi_{t^2}(1,t^2,\ldots,t^{2n-2};t,t^3,\ldots)\Pi_{t^2}(1,t^2,\ldots,t^{2n-2};t,t^2,\ldots)^k}
\end{equation}

\item (Hermitian case) Let $A\in\Her_n(\mathfrak{o})$ be the same as in \Cref{thm:The Hall-Littlewood measure of i.i.d entries}, and $B_1,\ldots,B_k\in M_{n\times n}(\mathfrak{o})$ be i.i.d random matrices with i.i.d entries distributed by the additive Haar measure on $\mathfrak{o}$. Then, the joint distribution  
$$\mathbf{P}(\SN^{\her}(A)=\lambda,\SN^{\her}(B_1^*AB_1)=\lambda^{(1)},\ldots,\SN^{\her}(B_k^*\cdots B_1^*AB_1\cdots B_k)=\lambda^{(k)})$$
is given by

\begin{equation}
\frac{Q_{\lambda}(t,-t^2,\ldots;-t)Q_{\lambda^{(1)}/\lambda}(t,-t,\ldots;-t)\cdots Q_{\lambda^{(k)}/\lambda^{(k-1)}}(t,-t,\ldots;-t)P_{\lambda^{(k)}}(1,\ldots,(-t)^{n-1};-t)}{\Pi_{-t}(1,-t,\ldots,(-t)^{n-1};t,-t^2,\ldots)\Pi_{-t}(1,-t,\ldots,(-t)^{n-1};t,-t,\ldots)^k}
\end{equation}
Here $t,-t,\ldots$ is the abbreviation of $t,-t,t^2,-t^2,t^3,-t^3,\ldots$, and $t,-t^2,\ldots$ is the abbreviation of $t,-t^2,t^3,-t^4,t^5,\ldots$.
\end{enumerate}
\end{cor}

\begin{proof}[Proof of \Cref{thm:joint distribution}]
Due the formula in \Cref{thm:Product process}, the product of random matrices can be regarded as product convolution defined in \Cref{defi: product convolution} with specializations $\theta=(1,t^2,\ldots,t^{2n-2}),\theta'=(1,t,\ldots,t^{2n-1})$ (resp. $\theta=(1,-t,\ldots,(-t)^{n-1}),\theta'=(1,t^2,\ldots,t^{2n-2})$) for the alternating (resp. Hermitian) case. Therefore, by \Cref{thm: product convolution and HL process}, the joint distribution of matrix product has the form of the Hall-Littlewood process.  

We first show the alternating case. By \cite[Corollary 1.4]{van2020limits}, we know for all $i\ge 1$, $\SN(B_i)$ is distributed with respect to the Hall-Littlewood measure of parameter $t$ with specializations $\theta'$ and $\psi_i=(t,t^2,t^3,\ldots)$.  Then we deduce from \eqref{eq: alt_joint distribution of product convolution} that the probability 

$$\mathbf{P}(\SN^{\alt}(B_\tau^T\cdots B_1^TAB_1\cdots B_\tau)=\lambda^{(\tau)},0\le\tau\le k)$$ 
is equal to

\begin{equation}\frac{Q_{\lambda}(t,t^3,\ldots;t^2)Q_{\lambda^{(1)}/\lambda}(t,t^2,\ldots;t^2)\cdots Q_{\lambda^{(k)}/\lambda^{(k-1)}}(t,t^2,\ldots;t^2)P_{\lambda^{(k)}}(1,t^2,\ldots,t^{2n-2};t^2)}{\Pi_{t^2}(1,t^2,\ldots,t^{2n-2};t,t^3,\ldots)\Pi_{t^2}(1,t^2,\ldots,t^{2n-2};t,t^2,\ldots)\cdots\Pi_{t^2}(1,t^2,\ldots,t^{2n-2};t,t^2,\ldots)}
\end{equation}

In particular, $\SN^{\alt}(B_k^T\cdots B_1^TAB_1\cdots B_k)$ is distributed with respect to the Hall-Littlewood measure with indeterminate $t^2$: the probability $\mathbf{P}(\SN^{\alt}(B_k^T\cdots B_1^TAB_1\cdots B_k)=\lambda)$ is equal to

\begin{equation}
\frac{P_\lambda(1,t^2,\ldots,t^{2n-2};t^2)Q_\lambda(t,t^3,t^5,\ldots,t,t^2,\ldots,t,t^2,\ldots,\ldots,t,t^2,\ldots;t^2)}{\Pi_{t^2}(1,t^2,\ldots,t^{2n-2};t,t^3,t^5,\ldots,t,t^2,\ldots,t,t^2,\ldots,\ldots,t,t^2,\ldots;t^2)}
\end{equation}
where there are $k$ copies of $t,t^2,\ldots$ in the specialization.

Likewise, for the Hermitian case, for all $i\ge 1$, $\SN(B_i)$ is distributed with respect to the Hall-Littlewood measure of parameter $t^2$ with specializations $\theta'$ and $\psi_i=(t^2,t^4,\ldots)$. Set $\psi_i^*=(t,-t,t^2,-t^2,t^3,-t^3,\ldots)$. Then we deduce from \eqref{eq: her_joint distribution of product convolution} that the probability $\mathbf{P}(\SN^{\her}(B_\tau^*\cdots B_1^*AB_1\cdots B_\tau)=\lambda^{(\tau)},0\le\tau\le k)$ is equal to

\begin{equation}\frac{Q_{\lambda}(t,-t^2,\ldots;-t)Q_{\lambda^{(1)}/\lambda}(t,-t,\ldots;-t)\cdots Q_{\lambda^{(k)}/\lambda^{(k-1)}}(t,-t,\ldots;-t)P_{\lambda^{(k)}}(1,\ldots,(-t)^{n-1};-t)}{\Pi_{-t}(1,-t,\ldots,(-t)^{n-1};t,-t^2,\ldots)\Pi_{-t}(1,-t,\ldots,(-t)^{n-1};t,-t,\ldots)^k}
\end{equation}

In particular, $\SN^{\her}(B_k^*\cdots B_1^*AB_1\cdots B_k)$ is distributed with respect to the Hall-Littlewood measure with indeterminate $-t$: the probability $\mathbf{P}(\SN^{\her}(B_k^*\cdots B_1^*AB_1\cdots B_k)=\lambda)$ is equal to

\begin{equation}
\frac{P_\lambda(1,-t,\ldots,(-1)^{n-1};-t)Q_\lambda(t,-t^2,t^3,\ldots,t,-t,\ldots,t,-t,\ldots,\ldots,t,-t,\ldots;-t)}{\Pi_{-t}(1,-t,\ldots,(-t)^{n-1};t,-t^2,t^3,\ldots,t,-t,\ldots,t,-t,\ldots,\ldots,t,-t,\ldots;-t)}
\end{equation}
where there are $k$ copies of $t,-t,t^2,-t^2,\ldots$ in the specialization.
\end{proof}

\begin{rmk}
    \label{rmk:ennola}
    The Ennola duality mentioned in \Cref{rmk:formal_macdonald_process} is the statement that (1) the conjugacy classes of the finite groups $\GL_n(\F_q)$ and $U_n(\F_{q^2})$ are canonically in bijection, as are the characters, and (2) each entry in one of their character tables (i.e. the values of characters on each conjugacy class) is a function $f^{\GL}(q)$ (resp. $f^{\U}(q)$) of $q$ such that $f^{\GL}(q) = f^{\U}(-q)$. 

    \Cref{thm:Corner process}, combined with previous work \cite{van2020limits}, gives a very similar statement concerning the corners processes associated with $\GL_n(F_0)$ and $\Her_n(F)$, where $F_0$ is a non-archimedean local field with residue field $\F_q$ and $F$ is one with residue field $\F_{q^2}$. Namely, \cite{van2020limits} yields\footnote{To obtain \eqref{eq:compare_ennola}, combine the fact that removing the bottom row of $A$ changes the singular numbers as in \cite[(1.4)]{van2020limits}, and subsequently removing the rightmost column to obtain $A_{n-1}$ changes the singular numbers as in \cite[(1.3)]{van2020limits}. One should take $n=N$ and $k=d=1$ in these equations, and replace $\mu$ in (1.4) by $\kappa$ and $\la$ in (1.3) by $\mu$ to match the notation of our \eqref{eq:compare_ennola}. One must also extend \cite[Theorem 1.3]{van2020limits} from $\Q_p$ to general $F$, but this extension holds easily in light of \cite[Remark 4]{van2020limits}.} that for a matrix $A \in \Mat_n(F_0)$ with singular numbers $\la$ and distribution invariant under $\GL_n(\mf{o}\cap F_0)$ on the right and left, and $A_{n-1}$ its top-left $(n-1) \times (n-1)$ submatrix,
    \begin{equation}\label{eq:compare_ennola}
        \mathbf{P}(\SN(A_{n-1}) = \nu) = \sum_{\kappa \in \Sig_{n-1}} \frac{P_{\la/\kappa}(1;t) Q_{\nu/\kappa}(1;t) P_\nu(t,\ldots,t^{n-1};t)}{P_\la(1,\ldots,t^{n-1};t)\Pi_t(1;t,\ldots,t^{n-1})},
    \end{equation}
    where $t=1/q$. Now, \eqref{eq:compare_ennola} is exactly the same as \eqref{eq:herm_corner} from \Cref{thm:Corner process}, except that $t=-1/q$ there and the sign of the argument of $Q_{\nu/\kappa}$ is flipped. We hope that a better understanding of Ennola duality can elucidate why this relation between \eqref{eq:herm_corner} and \eqref{eq:compare_ennola} holds, and how it may be generalized. It is also interesting that the partition $\kappa$ in \eqref{eq:compare_ennola} has a very concrete meaning---the singular numbers of $A$ after the bottom row is removed---but we are not aware of any meaning for the partition $\kappa$ in \eqref{eq:herm_corner}. As discussed in \Cref{rmk:formal_macdonald_process}, $\kappa$ cannot be an intermediate step of Markov dynamics as it is in \eqref{eq:compare_ennola}, because the transition probabilities would be negative.
\end{rmk}

\subsection{Explicit examples.} Since principally specialized Hall-Littlewood polynomials have very explicit formulas in \eqref{eq: Principal specialization formulas}, we would be remiss not to explain what the probabilities in \Cref{thm:corners_nomarkov_intro} and related results reduce to. In the alternating case, we remark that there is a relatively simple explicit probabilistic description/sampling algorithm for the Markov dynamics in the odd-to-even step of the alternating corners process (\Cref{thm:Corner process} (a)), see \cite[Proposition 5.3]{van2020limits}. An explicit description of the Markov dynamics in the even-to-odd step (\Cref{thm:Corner process} (b)) can be obtained similarly to that one. 

\begin{example}
Consider the $2\times 2$ corner $A_2$ of the random alternating matrix $A\in\Alt_3(F)$ with $\SN^\text{alt}(A)$ fixed, where the distribution is invariant under $\GL_3(\mathfrak{o})$ simultaneously on both sides.  By multiplying a power of $\pi$ to $A$, there is no loss we  study the case $\SN^{\alt}(A)=(0)$, Then $\SN^{\alt}(A_2)=(m)$ has distribution (here $t=1/q$)

$$\mathbf{P}_{2<3}^{\alt}((m)\mid(0))=\frac{Q_{(m)}(t;t^2)P_{(m)}(1;t^2)}{\Pi_{t^2}(t;1)}=\begin{cases}
1/(1+t+t^2) & m=0\\
(1-t^2)t^m/(1+t+t^2) & m>0\\
\end{cases}.$$ 
\end{example}

In the Hermitian case, things are less clear. The reason we have left the statements in terms of Hall-Littlewood polynomials is that the explicit formulas one gets via \Cref{prop:hl_principal_formulas} are not nearly as simple as, for example, the GUE corners process. There is significant cancellation in the sum over $\kappa$ in \Cref{thm:Corner process}, as the following example shows, though we have not tried to understand this in the general case.

\begin{example}
Consider the $1\times 1$ corner $a_{11}$ (i.e., the entry on the left upper side) of the random Hermitian matrix $A\in\Her_2(F)$ with $\SN^{\her}(A)=\lambda$ fixed, where the distribution is invariant under $\GL_2(\mathfrak{o})$ simultaneously on both sides. By multiplying a power of $\pi$ to $A$, there is no loss we study the case $\SN^{\her}(A)=(m,0)$, where $m\ge 0$. In this case, $v(a_{11})=n$ is distributed with respect to the law (here $t=1/q$)

\begin{equation}\label{eq: 1x1 corner of 2x2 Hermitian}\mathbf{P}_{1<2}^{\her}((n)\mid(m))=\sum_{\kappa\ge 0} \frac{P_{(m)/(\kappa)}(1;-t) Q_{(n)/(\kappa)}(-1;-t)P_{(n)}(-t;-t)}{P_{(m)}(1,-t;-t)\Pi_{-t}(1;t)}\end{equation}
If $m=0$, this is the Hermitian case of \Cref{thm:Corner of the invertible matrix}, and we have
$$\text{RHS of \eqref{eq: 1x1 corner of 2x2 Hermitian}}=\begin{cases}
(1-t)/(1+t^2) & n=0\\
(1-t^2)t^n/(1+t^2) & n>0\\
\end{cases}$$
If $m>0$ is odd, then we have
$$
\text{RHS of \eqref{eq: 1x1 corner of 2x2 Hermitian}}=\begin{cases}
1/(1+t^2) & n=0\\
(1-t^2)t^n/(1+t^2) & 0<n<m, n\text{ even}\\
t^{m+1}/(1+t^2) & n=m
\end{cases}
$$
If $m>0$ is even, then we have
$$
\text{RHS of \eqref{eq: 1x1 corner of 2x2 Hermitian}}=\begin{cases}
1/(1+t^2) & n=0\\
(1-t^2)t^n/(1+t^2) & 0<n<m, n\text{ even}\\
(1-t-t^2)t^m/(1+t^2) & n=m\\
(1-t^2)t^n/(1+t^2)& n>m
\end{cases}
$$
The cases that we do not write down have zero probability.
\end{example}

For the Hermitian case, in contrast to the alternating case, we are not aware of a direct sampling algorithm or more concrete probabilistic description. However, linear-algebraic computations can help explain at least basic features of the support of the distribution in the above example. Writing $A=B\diag(\pi^{m},1)B^*$, where $B=\begin{pmatrix} b_{11} 
& b_{12} \\ b_{21} 
& b_{22}\end{pmatrix}\in\GL_2(\mathfrak{o})$ is Haar distributed,  we have $a_{11}=\pi^m\Nm(b_{11})+\Nm(b_{12})$ where $\Nm: F \to F_0$ is the norm map. From this it is clear, for example, that $v(a_{11})$ cannot be an odd number less than $m$, because the valuation of $\Nm(b_{12})$ is always even. It would be interesting to find a more concrete explicit form or sampling algorithm for the probabilities in the Hermitian case.



\appendix 

\section{Comments on positive characteristic}\label{sec:char_appendix}

As mentioned in \Cref{footnote:pos_char} in the Introduction, we believe that all results stated in this paper remain valid over any non-archimedean local field with finite residue field, provided the characteristic of the local field is not $2$ in the Hermitian case. Moreover, we believe that the proofs should be essentially the same. However, we rely on many results from the literature which were proven in characteristic zero, as far as we can tell for reasons of motivation rather than technical reasons. Our interest is in $p$-adic fields, so we have not carefully checked that everything goes through in positive characteristic. At the same time, it seemed potentially useful to offer a roadmap to carefully verifying them, for later readers who may be interested in the positive characteristic case. That is the purpose of this Appendix.

The only results of this paper which use characteristic $0$ are \Cref{thm: alt_Hironaka} (alternating) and \Cref{thm: her_Hironaka} (Hermitian), which are essentially contained in Hironaka's papers on spherical functions. For the alternating case we use \cite{Hironaka}, and for the Hermitian case we use \cite{Hironaka_Hermitian}, which completes the chain of previous papers \cite{Hironaka_Hermitian_and_Symmetric_I}, \cite{Hironaka_Hermitian_and_Symmetric_II}, and \cite{Hironaka_Hermitian_and_Symmetric_III}. All of these papers assume characteristic zero, but as far as we have checked, the assumption is not used. However, these papers rely on other results, some of which are proven in arbitrary characteristic, but some proven only in characteristic zero. The latter are:
\begin{enumerate}
\item Hironaka defines the integral $\zeta(x;s)$ in \cite[(2.1)]{Hironaka} (alternating case) and \cite[(1.2)]{Hironaka_Hermitian} (Hermitian case). Both works claim that these are rational functions, and thus have analytic continuation to the whole complex plane, see the top of page 484 of \cite{Hironaka} (alternating) and \cite[Remark 1.1]{Hironaka_Hermitian} (Hermitian). These claims are based on Denef's study of Igusa zeta functions, which requires the field to have characteristic zero, so we can apply Heisuke Hironaka's reduction of singularities: see \cite{denef1984rationality}, \cite{denef1983evaluation}. Nevertheless, we believe that the positive characteristic reduction of singularities still holds in our particular case.


\item The proof of \cite[Lemma 1.8]{Hironaka_Hermitian} cites \cite[Lemma 2.3]{sato1989functional}, which deals with Fourier transform of $p$-adic prime powers. We believe that careful readers can check sentence by sentence that the proof in \cite{sato1989functional} does not require the characteristic to be zero. 

\item Page 208 of \cite{Hironaka_Hermitian_and_Symmetric_I} cites \cite{satake1963theory}, which claims that the Satake isomorphism gives an algebra isomorphism of the Hecke algebra onto the ring of symmetric Laurent polynomials. As shown in \cite[Chapter V]{Macdonald}, this also holds for the function field case.
\end{enumerate}

\newcommand{\etalchar}[1]{$^{#1}$}

\end{document}